\documentclass[11pt]{article}
\def\hybrid{\topmargin 0pt \oddsidemargin 0pt
        \headheight 0pt \headsep 0pt
        \textwidth 160true mm % US paper
        \textheight 231true mm % US paper
        \marginparwidth 0.0in
        \parskip 0pt plus 1pt \jot = 1.5ex}
\catcode`\@=11
\def\marginnote#1{}

\newcount\hour
\newcount\minute
\newtoks\amorpm
\hour=\time\divide\hour by60 \minute=\time{\multiply\hour by60
\global\advance\minute by-\hour}
\edef\standardtime{{\ifnum\hour<12
\global\amorpm={am}%
        \else\global\amorpm={pm}\advance\hour by-12 \fi
        \ifnum\hour=0 \hour=12 \fi
        \number\hour:\ifnum\minute<10
        0\fi\number\minute\the\amorpm}}
\edef\militarytime{\number\hour:\ifnum\minute<10
0\fi\number\minute}

\def\draftlabel#1{{\@bsphack\if@filesw {\let\thepage\relax
   \xdef\@gtempa{\write\@auxout{\string
      \newlabel{#1}{{\@currentlabel}{\thepage}}}}}\@gtempa
   \if@nobreak \ifvmode\nobreak\fi\fi\fi\@esphack}
        \gdef\@eqnlabel{#1}}
\def\@eqnlabel{}
\def\@vacuum{}
\def\draftmarginnote#1{\marginpar{\raggedright\scriptsize\tt#1}}

\def\draft{\oddsidemargin -.5truein
        \def\@oddfoot{\sl preliminary draft \hfil
        \rm\thepage\hfil\sl\today\quad\militarytime}
        \let\@evenfoot\@oddfoot \overfullrule 3pt
        \let\label=\draftlabel
        \let\marginnote=\draftmarginnote
   \def\@eqnnum{(\theequation)
   \rlap{\kern\marginparsep\tt\@eqnlabel}%
\global\let\@eqnlabel\@vacuum} } %%%%%%%%%%%
\relax

\hybrid

\input{amssym.def}
\input{amssym}
\usepackage{eufrak}

\newtheorem{theorem}{Theorem}

\newtheorem{corollary}[theorem]{Corollary}
\newtheorem{definition}[theorem]{Definition}

\newtheorem{proposition}[theorem]{Proposition}
\newtheorem{remark}[theorem]{Remark}

\def\beq{\begin{equation}}
\def\be{\begin{equation}}
\def\eeq{\end{equation}}
\def\ee{\end{equation}}
\def\K{{\Bbb K}}
\def\C{{\Bbb C}}
\def\R{{\Bbb R}}
\def\Z{{\Bbb Z}}
\def\gg{{\frak g}}
\def\gh{{\frak g}_{\h}}

\def\uq{U_q(sl(2))}
\def\Uq{U_q(sl(n))}
\def\Uqq{U_q(\gg)}
\def\gg{\mbox{$\frak g$}}
\def\ot{\otimes}
\def\h{{\hbar}}
\def\qq{q^{-1}}
\def\vv{V^{\ot 2}}
\def\Id{{\rm Id}}
\def\Mat{{\rm Mat}}
\def\Ren{R_{\End}}
\def\span{{\rm span}}
\def\Sym{{\rm Sym}}
\def\Im{{\rm Im}}
\def\Tr{{\rm Tr}}
\def\ad{{\rm ad}}
\def\End{{\rm End}}

\def\det{{\rm det}}

\def\Ind{{\rm Ind}}
\def\diag{{\rm diag}}
\def\Vect{{\rm Vect}}
\def\De{{\Delta}}
\def\de{{\delta}}
\def\Gr{{\rm Gr}}
\def\Cas{{\rm Cas}}
\def\LL{{\cal L}}
\def\TT{{\cal T}}
\def\al{\alpha}
\def\bw{{\bigwedge}}

\def\O{{\cal O}}
\def\Om{{\Omega}}
\def\lq{{\cal L}(R_q)}
\def\lhq{{\cal L}({R_q,\h})}
\def\loq{{\cal L}({R_q,1})}
\def\sloq{{\cal SL}({R_q,1})}
\def\slq{{\cal SL}({R_q})}

\def\lrqh{{\cal L}({R_q,\h})}

\def\slhq{{\cal SL}({R_q,\h})}
\def\slhqc{{\cal SL}^C({R_q,\h})}
\def\slqc{{\cal SL}^C({R_q})}
\def\sloqc{{\cal SL}^C({R_q,1})}
\def\loqc{{\cal L}^C({R_q,1})}
\def\hh{{\frak h}}

\def\Lie{{\rm Lie}}
\def\la{{\lambda}}
\def\ggg{{\gg}_{\Gamma}}
\def\L{{\cal L}}
\def\SL{{\cal SL}}
\def\ddk{{d_{RK}}}

\def\ii{{\textbf{i}}}
\def\tlhq{\tilde{\cal L}(R_q,\hbar)}
\def\tlq{\tilde{\cal L}(R_q)}

\def\tip{\tilde{I}_+}
\def\tim{\tilde{I}_-}
\def\krrq{\K_q[\R^4]}
\def\krq{\K_q[\R^3]}
\def\Mw{\rm Mw}
\def\Om{\Omega}
\def\om{\omega}
\def\parho{\partial_\rho}
\def\kr{\K[\R^3]}
\def\krr{\K[\R^4]}
\def\Dir{\rm Dir}
\def\rhoq{\rho_q}
\def\BBB{{\cal B}_q}
\def\HHH{{\cal H}_q}
\def\CCC{{\cal C}_q}
\def\khq{\K_q[H^2]}
\def\twl{{\tilde\bw_q}(\LL)}
\def\Trr#1{{\rm Tr\str{-1.3}}_{R^{\mbox{\scriptsize
$(#1)$}}}}
\def\str#1{\rule[#1mm]{0pt}{1mm}}

\begin{document}

\makeatletter
\renewcommand{\theequation}{{\thesection}.{\arabic{equation}}}
\@addtoreset{equation}{section}
\makeatother

\title{Braided affine geometry and $q$-analogs of wave operators}
\author{\rule{0pt}{7mm} Dimitri
Gurevich\thanks{gurevich@univ-valenciennes.fr}\\
{\small\it LAMAV, Universit\'e de Valenciennes,
59313 Valenciennes, France}\\
\rule{0pt}{7mm} Pavel Saponov\thanks{Pavel.Saponov@ihep.ru}\\
{\small\it Division of Theoretical Physics, IHEP, 142281 Protvino,
Moscow region, Russia}}

\maketitle

\begin{abstract}
The main goal of this review is to compare different approaches to constructing geometry
associated with a Hecke type braiding (in particular, with that related to the quantum group
$U_q(sl(n))$). We make an emphasis on affine braided geometry related to the so-called
Reflection Equation Algebra (REA). All objects of such type geometry  are defined in
the spirit of affine algebraic geometry via polynomial relations on generators.

We begin with comparing the Poisson counterparts of "quantum
varieties" and describe different approaches to their
quantization. Also, we exhibit two approaches to introducing
$q$-analogs of vector bundles and  defining the Chern-Connes
index for them on quantum spheres. In accordance with the
Serre-Swan approach, the $q$-vector bundles are treated as
finitely generated projective modules over the corresponding
quantum algebras.

Besides, we describe the basic properties of the REA used in this construction and
compare different ways of defining $q$-analogs of partial derivatives and differentials
on the REA and algebras close to them. In particular, we present a way of introducing a
$q$-differential calculus via Koszul type complexes. The elements of the $q$-calculus
are applied to defining $q$-analogs of some relativistic wave operators.
\end{abstract}

\section{Introduction}

By a braided geometry we mean a sort of noncommutative
geometry related to a {\it braiding} which is defined as follows.
Let $V$ be a finite dimensional vector space over the ground field $\K$
(of complex numbers $\C$ or real numbers $\R$). An operator
\be
R:\vv\to\vv \label{solYB}
\ee
is called {\it a braiding} provided  it satisfies the following
relation on the space $V^{\otimes 3}$
\be
R_{12}R_{23}R_{12} = R_{23}R_{12}R_{23}.
\label{YBE}
\ee
Here we use the standard notations $R_{12} = R\otimes \Id$ and
$R_{23} = \Id\otimes R$,  where  $\Id$ is the identity operator on the
space $V$. In fact, relation (\ref{YBE}) has the meaning of a very
special representation of the Artin braid group. Besides, this relation is
equivalent to the quantum Yang-Baxter equation and often is called
the Yang-Baxter equation too.

Let us give a few examples of solutions of the Yang-Baxter equation.
The first examle is the classical flip $\sigma$ which transposes any two elements
$\sigma(x\ot y)=y\ot x$, $x,y\in V$, is a solution. The second example is
related to a $\Z_2$-graded vector space $V=V_{\bar{0}}\oplus V_{\bar{1}}$
where the flip is replaced by its super-analog $\sigma(x\ot y)=(-1)^{\overline{x}
\,\overline{y}} y\ot x$. Here $x$ and $y$ are homogeneous elements of $V$
and $\overline{x},\,\,\overline{y} \in \Z_2$ are their parities. Emphasize, that
all super-flips (the classical flips included) are involutive, i.e. $\sigma^2=\Id$.

Other known examples come from the quantum groups (QG) $\Uqq$ (see \cite{CP}).
Consider a finite dimensional $\Uqq$-module $V$ and take the image of the universal
$R$-matrix in $\End(V^{\otimes 2})$. Then, the product of this image and the
flip $\sigma$ gives a solution of (\ref{YBE}).

If $\gg=sl(n)$ and $V$ is the first fundamental $U_q(sl(n))$-module ($q\in \K$ is a
fixed non-zero number), then the corresponding braiding $R$ satisfies the second
degree equation
\be
(R-q\Id)(R+\qq\Id)=0.
\label{Heck-cond}
\ee
In this case, the representation of the group
algebra of the Artin braid group is reduced to the representation of the Hecke
 algebra. For this reason a braiding $R$ satisfying the additional condition
(\ref{Heck-cond}) is called {\em the Hecke symmetry}.

If an algebra $\gg$ belongs to the series $B_n,\,C_n$ or $D_n$, then the
corresponding braiding $R:V^{\ot 2}\to V^{\ot 2}$, where $V$ is also the
first fundamental $\gg$-module, satisfies a third degree equation. We call it
{\em a Birman-Murakami-Wenzl (or BMW) symmetry}.
These BMW symmetries as well as the aforementioned Hecke symmetries are
deformations of the classical flips. By using the so-called gluing procedure
(see \cite{G2, MM, GPS3}) it is possible to construct the Hecke symmetries
which are deformations of the super-flips. In order to point out the braidings
and symmetries which are deformations of the classical flips we call them
{\em quasiclassical}. Note, that there exists a big family of other braidings
which are not deformations of either classical or super-flips.

Given a braiding $R$ (quasiclassical or not), the following very natural question
arises: which associative algebras can be connected with it? The simplest examples
are $q$-analogs of the symmetric $\Sym(V)$ and skew-symmetric $\bigwedge(V)$
algebras on the space $V$ endowed with a Hecke symmetry (\ref{solYB}). They are
respectively defined as follows
\be
 \Sym_q(V)=T(V)/\langle \Im (q\,\Id-R)\rangle,\qquad
\bw_q(V)=T(V)/\langle \Im (\qq\,\Id+ R)\rangle.
\label{odin}
\ee
Hereafter $T(V)$ stands for the free tensor algebra of a given
space $V$ and $\langle S \rangle$ is the two-sided ideal generated
by a subset $S\subset T(V)$. As follows from results of \cite{G2},
these algebras have {\em a good deformation property}. For a
 quasiclassical Hecke symmetry $R$ this means that for a
generic $q$ and all positive integers $k$, we have
\be
\dim \Sym^{(k)}_q(V)=\dim
\Sym^{(k)}(V),\qquad \dim \bw_q^{k} (V)=\dim \bw^{k}(V)
\label{good-dp}
\ee
where $\Sym^{(k)}_q(V)$ and $\bw_q^{k}(V)$
are the $k$-th order homogeneous components of the algebras
$\Sym_q(V)$ and $\bw_q(V)$, respectively.

If $R$ is a BMW symmetry,  $q$-analogs of symmetric and
skew-symmetric algebras of the space $V$ can be introduced as
well\footnote{Observe that, in general, for an arbitrary braiding $R:\vv\to\vv$ the
"$R$-analogs" of the symmetric $\Sym(V)$ and skew-symmetric
$\bigwedge(V)$ algebras are not defined. The definition
of these algebras as the quotients $T(V)/\langle \Im
(\Id-R)\rangle$ and $T(V)/\langle \Im (\Id+ R)\rangle$) respectively is
acceptable only if 1 and -1 are the eigenvalues of $R$ (otherwise
these quotients are trivial). But, in general, even under this condition the good
deformation propriety of the algebras involved (assuming $R$ to
be quasiclassical) is not ensured.}.

Two other well known examples belong to the class of the so-called
quantum matrix algebras.
(All algebras below are assumed to be unital.) These are the RTT algebra and the
Reflection Equation Algebra (REA). The RTT algebra is an
associative algebra generated by  formal indeterminates
$T_i^{\,j}$, $1\leq i,j \leq n=\dim V$, subject to the system of
relations \cite{FRT} \be R_{12}T_1T_2 -T_1T_2R_{12} = 0,
\label{RTT} \ee which is the compact notation for the matrix
equation
$$
R(T\ot \Id)(\Id\ot T)-(T\ot \Id)(\Id\ot
T)R=0.
$$
Here $T=\|T_i^{\,j}\|$ is the $n\times n$ quantum matrix with
noncommutative entries $T_i^{\,j}$.

The REA is another  associative algebra with formal
generators $L_i^{\,j}$, $1\leq i,j \leq n=\dim V$
subject to the relations
\be
R_{12}L_1R_{12}L_1 - L_1R_{12}L_1R_{12} = 0,\quad L_1=L\ot\Id
\label{RE}
\ee
where $L=\|L_i^{\,j}\|$.
Upon replacing $L_1$ in this relation by $L_2=\Id\ot L$ we get a
similar algebra which will be referred to as {\it the REA of second type}.

Besides these two examples, there are other quantum matrix
algebras associated with pairs of compatible
braidings. We discuss them in section \ref{sec:3}.

Also, note that the Reflection Equation (with a parameter)
appeared in connection with the theory of integrable
models with a boundary in \cite{C} (see also \cite{KS}).

If $R$ is a quasiclassical Hecke symmetry, then the RTT algebra
and the REA have the good deformation property
and can be treated as two different $q$-analogs of the symmetric
algebra $\Sym(\End(V))$, where $\End(V)$
stands for the vector space of endomorphisms of the space
$V$. In this case we denote the algebras (\ref{RTT})
and (\ref{RE}) as $\Sym_q(\TT)$ and $\Sym_q(\LL)$ respectively.
Here $\TT=\span(T_i^{\,j})^n_{i,j=1}$ and $\LL=\span(L_i^{\,j})^n_{i,j=1}$
are the linear spans of the corresponding generators. The
$q$-analog $\bw_q(\TT)$ (resp., $\bw_q(\LL)$) of the skew-symmetric
algebra $\bw(\End(V))$ with the good deformation property can
also be associated to any quasiclassical Hecke symmetry
(see section \ref{sec:7}).

Being equipped with the usual matrix coproduct, the algebras
$\Sym_q(\TT)$ and $\Sym_q(\LL)$ become bialgebras
(in the latter case {\em braided}). Besides, if $R$ is an
{\em even}\footnote{This means that it is skew-invertible (see
section \ref{sec:3}) and the skew-symmetric algebra $\bw_q(V)$ is
finite-dimensional.} Hecke symmetry, then in each of these
algebras there is a $q$-analog of the determinant of the matrix $T$
or $L$ which is the group-like element with respect
to the matrix coproduct. Let $\det_q T$ (resp., ${\rm Det}_q L$) be
this analog in the algebra $\Sym_q(\TT)$ (resp.,
$\Sym_q(\LL)$).

Assuming $\det_q T$ to be central (${\rm Det}_q
L$ is always central) we can define the quotients
$$\Sym_q(\TT)/\langle \det_q T-1\rangle\quad {\rm and} \quad  \Sym_q(\LL)/\langle
{\rm Det}_q L-1\rangle. $$
These quotients are the
Hopf algebras (in the latter case the {\em braided} Hopf algebra)
and can be treated as two different deformations of
the algebra $\K[SL(n)]$, provided $R$ is a quasiclassical Hecke symmetry.
The Hopf structure in the above quotient of RTT
algebra first appeared in papers of the
Leningrad mathematical school (see \cite{FRT} and the references
therein) and afterwards was formalized by V.~Drinfeld
\cite{Dr1}. The braided Hopf structure in the above quotient of the REA was discovered
by S.~Majid (see \cite{M4} and the references therein).

Note that if $G$ is a group from the series $B_n$, $C_n$ or
$D_n$, then there exist similar deformations of the algebra $\K[G]$.
They can be also realized as appropriate quotients of the RTT
algebra and of the REA respectively. In the sequel, the
notation $\K_q[G]$ stands for the quantum deformation of the
algebra $\K[G]$ which is the mentioned quotient of the
RTT algebra.

In the present paper we review different ways of introducing quantum
(braided) analogs of coordinate algebra of an affine regular variety
and compare the roles of the RTT algebra and
the REA in braided geometry. Also, we
 exhibit a regular way  of defining $q$-analogs of some relativistic wave
operators (Laplace, Maxwell, Dirac ones) on a $q$-analog of
the Minkowski space algebra.

We are mainly interested in a braided version of  {\em affine} algebraic geometry.
This means that all algebras we are dealing with are introduced by means of
some polynomial relations on their generators. The coefficients of these relations
analytically depend on the deformation parameter $q$ which can be specialized.
The RTT algebra, the REA, their quotients mentioned above, and their skew-symmetric
counterparts are examples of such algebras.

Other examples are provided by the so-called $q$-(quantum, braided) varieties.
The corresponding algebras arise from a quantization of commutative algebras
of functions on some classical varieties.
A typical example of a classical variety to be quantized is an orbit $\O\subset \gg^*$ of a
semisimple element $a\in \gg^*$ (a {\em semisimple orbit} for shot) where $\gg$ is a simple
Lie algebra (our main example is $sl(n)$) corresponding to a complex connected Lie group $G$.
Let $\K[\O]$ be the coordinate algebra of the affine variety $\O\subset \gg^*$.

Another way to realize the above coordinate algebra is based on the fact that the orbit
$\O$ is isomorphic (as a $G$-set) to a coset $G/G_a$, where $G_a\subset G$ is the
stabilizer of the element $a$. Via the map dual to the projection $G\to \O$ we can realize
the space of functions on $\O$ as a subalgebra of $\K[G]$ consisting of functions
such that $f(xb)=f(x)$, where $ x\in G$, $b\in G_a$. Abusing the notations we denote this
subalgebra $\K[G/G_a]$.

In general, upon deforming the algebras $\K[\O]$ and $\K[G/G_a]$, we get non-equivalent
quantum algebras. Moreover, if $G/G_a$ is not a symmetric homogeneous space, the
algebra $\K[G/G_a]$ has a large family of quantum deformations. Any algebra $A$ of this
family is a covariant $\Uqq$-module and its product (denoted $\circ$) is coordinated with
$\Uqq$ action in the following sense
\be
X(a\circ b)=\circ (X_1(a)\otimes X_2(b)),\qquad \forall\, a,b \in A,
\quad X\in \Uqq \label{coo}
\ee
where $\De(X)=X_1\ot X_2$ is the Sweedler's notation for the coproduct in the QG $\Uqq$.

In the next section we consider semiclassical counterparts of the
above quantum algebras on semisimple
orbits\footnote{We call the corresponding Poisson brackets
$G$-covariant ones. We reserve the term
"a Poisson-Lie bracket" (which is often used for these brackets) for
the linear Poisson bracket on the space $\gg^*$.}.
Following \cite{DGS} we show that  on a generic
semisimple orbit $\O\subset \gg^*$ there is
 a family of $G$-covariant non-equivalent Poisson structures. This
family is parameterized by a variety of dimension equal to the rank
of $G$ excluding some subvarieties of smaller
dimensions. However, in general these Poisson structures can be
quantized only formally (i.e. the products in the corresponding
quantum algebras are represented by series whose
convergence is disregarded). One of such structures, the so-called
reduced Sklyanin bracket, can be quantized in terms of "quantum
cosets". More precisely, the resulting quantum algebra (denoted $\K_q[G/G_a]$)
is treated as a subalgebra of the algebra $\K_q[G]$ similarly to the classical
pattern.

Nevertheless, if $G=SL(n)$, then on any semisimple $G$-orbit there is a Poisson
pencil whose quantization can be described in the spirit of affine algebraic
geometry. Moreover, it is a restriction of a Poisson pencil defined on the
vector space $gl(n)^*$. The latter Poisson pencil
\be
\{\,\,,\,\}_{a,b}=a\{\,\,,\,\}_{PL}+b\{\,\,,\,\}_{REA}
\label{PP}
\ee
is generated by the linear Poisson-Lie bracket $\{\,\,,\,\}_{PL}$ and the Poisson
counterpart of the standard\footnote{In  the sequel we also use the term
{\em standard} for the Hecke symmetry (and other objects) related to the QG $\Uq$ .}
REA denoted $\{\,\,,\,\}_{REA}$. It is known that the bracket $\{\,\,,\,\}_{PL}$ can be
restricted to any orbit in $\gg^*$ for any Lie algebra $\gg$. We call this restricted
bracket the Kirillov-Kostant-Souriau (KKS) one and denote it   $\{\,\,,\,\}_{KKS}$.

The fact that the bracket $\{\,\,,\,\}_{REA}$ can be also restricted to any orbit in $gl(n)^*$
was proved by J. Donin \cite{D}. We would like to emphasize a great contribution to this
area made by our friend J. Donin passed away 4 years ago. Being an expert
in deformation theory, he studied many types of Poisson brackets (those on homogeneous
spaces included) and their quantizations.

Thus, the whole pencil (\ref{PP}) can be restricted to any orbit in $gl(n)^*$, and consequently,
to any  orbit in $sl(n)^*$. Observe that the quantum counterpart of the Poisson pencil (\ref{PP})
can be defined by polynomial relations which are a slight modification of those for the standard REA.
Let us describe the quantum algebra explicitly.

Given a skew-invertible symmetry $R$, we call the modified REA (mREA) the algebra
defined by the following system of quadratic-linear relations on its generators $L_i^{\,j}$
\be
R_{12}L_1 R_{12}L_1-L_1 R_{12} L_1 R_{12}-\h\,(R_{12}L_1-
L_1 R_{12})=0,
\label{mREA0}
\ee
where $L = \|L_i^{\,j}\|$. In what follows we  denote the mREA by $\lrqh$.

Note that for the standard braiding $R$ the algebra $\lrqh$ transforms into the enveloping algebra
$U(gl(n)_\h)$ as $q\to 1$. (Hereafter $\gh$ stands for the Lie algebra which differs
from the Lie algebra $\gg$ by the factor $\h$ in the Lie brackets.) As a consequence,
in the standard case the Poisson pencil (\ref{PP}) turns out to be the semiclassical counterpart of
the two parameter algebra $\lhq$.

Comparing  relations (\ref{RE}) and (\ref{mREA0}) we can also see that
for the standard Hecke symmetry $R$ the following holds $\Sym_q(\LL)=\LL(R_q,0)$.
In what follows the algebra $\LL(R_q,0)$ will be also denoted $\lq$.
Note that the algebra $\lq$ is actually isomorphic to $\lhq$ if $q\not=\pm 1$ (see section
\ref{sec:3}). But since the isomorphism breaks at $q=\pm 1$ we prefer to distinguish these
algebras and to use different notation for them.

Now, go back to the restriction of the Poisson pencil
(\ref{PP}) to a given semisimple orbit $\O$. We want to represent the
result of its quantization  as a
proper quotient of the two-parameter algebra $\lhq$.
To this end, we have to find "the coordinate ring of the corresponding
 quantum orbit". This problem is rather subtle if the orbit $\O$ is not
generic, i.e., the eigenvalues of the corresponding element $a\in
\O$ are not pairwise distinct. As was observed in
\cite{DM}, the degeneracy of the eigenvalues disappears after
quantization (see also \cite{GS2}, remark 18). We do not
consider such orbits here. Our main examples --- the hyperboloids
(spheres) --- are the generic orbits and, at the same time, these
orbits are symmetric due to their low dimension.

Note that any symmetric orbit $\O$ possesses an additional property: all
$G$-covariant brackets on such an orbit $\O$ originate from the restriction of the Poisson
pencil (\ref{PP}) to it. Thus, in this case any quantum coset
algebra $\K_q[G/G_a]$ can be also realized
as a quotient of the mREA. For example, the quantum sphere
algebra can be introduced as a quantum coset
$\K_q[SU(2)/H]$ and as a quotient of the mREA (see section
\ref{sec:6}). Note that the popular Podle\`s quantum sphere
algebra \cite{P1} is just such a quotient but written in terms of
generators different from ours and
endowed with an involution (conjugation) possessing classical properties.
However, contrary to the quantum hyperboloid algebra,
the Podle\`s quantum sphere cannot be realized as a real algebra
even if the parameter $q$ is real. We discuss the
problem of an appropriate definition of an involution in a quantum
algebra in sections \ref{sec:5} and \ref{sec:6}.

Nevertheless, even on symmetric orbits different ways of
introducing quantum algebras give rise to different types of
"quantum geometry". In section \ref{sec:6} we compare these
approaches on an example of a quantum sphere
(hyperboloid). In particular, we describe two ways of constructing
quantum analogs of line bundles and computing the
Chern-Connes index for them. Recall that according to the Serre-Swan
approach any vector bundle over a regular affine algebraic variety
or a smooth compact one can be realized as a projective module
over its coordinate algebra (all
projective modules are assumed to be finitely
generated). The Chern-Connes index $\Ind(\pi, e)$ is introduced via
a pairing of a representation $\pi$ of a given noncommutative
algebra $A$ and a projective $A$-module $M\cong e A^{\oplus n}$
(or $M\cong A^{\oplus n}e$) where
$e\in \Mat_n(A)$ is the corresponding idempotent. This pairing is defined as follows
\be
\Ind(\pi, e)=\Tr(\pi(\Tr(e))).
\label{idx}
\ee

In \cite{HM} a family of idempotents over the algebra $\K_q[SU(2)
/H]$ was constructed. Besides, one of these
idempotents was paired with an infinite dimensional representation
of the quantum sphere algebra taken from \cite{MNW}.
The quantum index thus calculated equals the value of the classical
index on the corresponding line bundle.

Another approach was suggested in \cite{GLS2} (also, see \cite{GS2}) where a
family of idempotents over a $q$-hyperboloid algebra was constructed, the $q$-hyperboloid
being realized as a quotient of the algebra $\lhq$. In
contrast with \cite{HM}  these idempotents were constructed via
a series of the Cayley-Hamilton identities valid for some matrices with entries from
the mREA (in particular, for the matrix $L$ entering the definition of the mREA).
Such type identities were found in \cite{GPS1, GPS2}
for the quantum matrix algebras associated with a large class of
general type Hecke symmetries.
However, only  for the REA (modified or not)  the  matrix powers
in the Cayley-Hamilton identities have the usual sense (see section \ref{sec:3}).

The second ingredient coming in the index formula is a
representation. The problem of constructing representations of
 the REA (mainly, in the standard case) was
considered in a series of papers (see, for example, \cite{K, Mu1,
KSS, DKM, DS, S}). Note that in the standard case the REA has
some specific properties simplifying the problem (see section \ref{sec:5}).
A general approach to constructing the finite dimensional,
$R$-invariant\footnote{In the $\Uq$ case this mean that the
maps $\lhq\to\End(V)$ are $\Uq$-morphisms.} (or equivariant)
representations was suggested in \cite{GPS3}. In that
paper a quasitensor category (called  Schur-Weyl)  of the REA representations
was defined for the REA associated with a skew-invertible
general type Hecke symmetry $R$.

An important peculiarity of this category is a modification
of the notion of the trace. For the quantum matrices the usual trace
must be replaced by the categorical (or quantum) trace $\Tr_R$
(see sections \ref{sec:3} and \ref{sec:4} for detail). In particular,  this
new trace enables us to define
the procedure of {\it sl-reduction} for the REA and its
representations. In a sense, this procedure is analogous to
the classical passage from $U(gl(n))$ to $U(sl(n))$ which is a
motivation for the term ``$sl$-reduction''.

Given a representation of the algebra $\lhq$ belonging to the
Schur-Weyl category and a projective module from the
aforementioned family, the $q$-index is defined in \cite{GLS2} by
a formula analogous to (\ref{idx}) but both usual traces coming in
it are replaced by the quantum trace $\Tr_R$. As a result, the
$q$-index equals a $q$-integer. This is an intrinsic property of
the braided affine geometry: all the numerical characteristics
(dimensions, indexes, etc.) of its objects (algebras, modules,
etc.) become $q$-numbers.

Another basic feature of the braided affine geometry is a modification of
the notion of a Lie algebra and a vector field. In the
 $\Uqq$ case the problem of defining a quantum (braided) Lie algebra
 can be formulated as
follows. We look for a deformation of the Lie bracket
$[\,\,,\,]:\gg^{\ot 2}\to \gg$ such that the deformed bracket $[\,\,,\,]
_q$ is a $\Uqq$-covariant map (we assume the space
$\gg$ to be endowed with a $\Uqq$ action which is a deformation
of the usual adjoint one) and the corresponding {\em
enveloping algebra} is a quadratic-linear one and it possesses a {\it good
deformation property}. This means that it is canonically isomorphic
to its associated graded algebra\footnote{We refer the reader to
\cite{PP} where a general form of the Jacobi identity and that of the PBW theorem
ensuring such an isomorphism are presented. Note that this Jacobi identity
has nothing in common with that from section \ref{sec:5} which enables us to define the
adjoint representation.} and the latter one has a good
deformation property in the sense of definition based on formula
(\ref{good-dp}). Also, we are interested in finding $q$-analogs of axioms of usual
(or super-) Lie algebras.

There are known numerous attempts
\cite{W,DGHZ,DGG,LS,M2,GM} to define the  quantum (braided)
analogs of a Lie
algebra (often without requiring the good deformation property for
the "enveloping algebra"). It turned out \cite{GPS3} that for
simple Lie algebras the braided deformations with the desired
properties exist only for the Lie algebra $\gg=sl(n)$. Also,
such a deformation can be defined for the general linear algebra $\gg=gl(n)$.
Moreover, a braided analog of the Lie algebra $gl(n)$ can be
associated with any skew-invertible Hecke symmetry $R$ and the role of the
corresponding enveloping algebra is played by the mREA
related to the symmetry $R$. Note that although the mREA had been known for a long
time, a "braided Lie algebra" was extracted from this algebra only in \cite{GPS3}.

By using this braided Lie algebra we can define a $q$-analog of the
adjoint action on the space $\LL$ in the natural way:
$\ad_q x(y)=[x,y]_q$, $\forall\,x,y\in \LL$. In the $\Uq$ case it is a
deformation of the usual adjoint action.
In order to define  braided analogs of  vector fields arising from
the usual adjoint action
we should extend this $q$-adjoint action to the higher components
of the algebra $\lq$ (or of $\slq$ if we deal with an $sl$-reduced algebra).
For usual or super-algebras this can be done by means of the coproduct acting on
the generators in the additive way: $\Delta(X) =  X\ot 1+1\ot X$. This leads to the
classical (or super-) Leibnitz rule for the vector fields. The above coproduct (properly
extended to higher degree elements) is still valid for the mREA related to an involutive
symmetry $R: R^2 = \Id$.

However, if a symmetry $R$ coming in the definition of the REA is not involutive, this
method of constructing the "braided vector fields" fails. Fortunately, in the algebra
$\loq$ (we put $\h=1$ for the sake of concreteness, since algebras $\lrqh$ are isomorphic
for all $\hbar\not=0$) there exists a coproduct which endows it with a {\em braided bialgebra
structure}. Note that in the $\Uq$ case it is also a deformation of the classical coproduct.
Thus, using this coproduct we can define the braided vector fields on the algebra $\lq$ or
its $sl$-reduced counterpart $\slq$ which are analogs of the vector fields arising from the
usual adjoint action.

We call these braided vector fields {\em tangent}  since their classical counterparts are
tangent to orbits in $sl(n)^*$. In the case $n=2$ the braided tangent vector fields are subject
to the relation analogous to $x\,X+y\,Y+z\,Z=0$ which is valid for the infinitesimal
rotations $X,\, Y,\, Z$ in the space $\R^3$, $x,y,z$ being $\R^3$ coordinates. Note that
the vector fields $X,Y,Z$ have the meaning of the angular momentum components.

In the present paper we use a different method of defining the tangent braided
vector fields. This method is based on a conjecture that any element $f\in\lq$ or $f\in\slq$
has a canonical (i.e. completely $q$-symmetrized) form. If it is so, the action of a braided
vector field on the element $f$ can be defined without using any form of the Leibnitz rule.
Namely, it suffices to apply the braided vector field to the first factors of the summands
constituting the canonical form of $f$ (with a subsequent renormalization).  This method
is close to the technique  used  in \cite{G3} for constructing  Koszul type complexes (see
section \ref{sec:7}). Note that these complexes are usually employed in the theory of quadratic
algebras \cite{Ma, PP}. In a similar manner, using the aforementioned conjecture we
introduce braided analogs of partial derivatives.

In the last section we use the braided analogs of the partial derivatives and the tangent
vector fields in order to define $q$-analogs of basic wave operators on the $q$-Minkowski
space and $q$-hyperboloid algebras. A $q$-analog of the Minkowski space algebra was
introduced  in the early 90's in  \cite{CSSW1, CSSW2, SWZ, OSWZ}. Initially, this algebra
was defined via $q$-analogs of spinors. Lately, it was treated as a particular case of the
REA \cite{M3, Me1, MMe, AKR}.

In section \ref{sec:8} we analyze other possible candidates for the role of the $q$-Minkowski space
algebra assuming them to be quadratic and $\uq$-covariant\footnote{So, accordingly to the
common viewpoint the complex QG $\uq$ is treated to be the quantum analog of the Lorenz
group. In this connection we would like to mention an original approach of \cite{Dob} where a
$q$-Minkowski space is $SU_q(2,2)$-covariant.}. We introduce a {\em truncated} REA $\krrq$
which differs slightly from (\ref{RE}) and in the $\uq$ case we treat it as the $q$-Minkowski
space algebra.

Also, $q$-analogs of the Laplace\footnote{Note that we employ the term "Laplace operator" in
 a large sense by admitting that the metric coming in its definition can be indefinite.} and Dirac
operators on the $q$-Minkowski space and quantum sphere algebras were considered in a
number of papers \cite{Me2, BK, P4, PS} and others. In fact, both operators arise from the
quadratic Casimir element $\Cas=\Tr_q L^2$ which is central in the algebra $\krrq$. Passing
to the $sl$-reduced algebra $\krq$ we get the reduced Casimir element $\Cas_{sl}$. Then we
obtain the $q$-Laplace operator on the algebra $\krq$ (resp., $\krrq$) by replacing the
algebra generators in the Casimir element $\Cas_{sl}$ (resp., $\Cas$) by appropriate
$q$-derivatives. In order
to get the $q$-Laplace operator on the $q$-sphere ($q$-hyperboloid) we
replace the generators of the algebra $\krq$ in $\Cas_{sl}$ by
braided analogs of some tangent vector fields. The element
$\Cas_{sl}$ is also involved in the
construction of $q$-Dirac operator on the algebras in question (see
section \ref{sec:8}).

Particular cases of such tangent vector fields are
 infinitesimal (hyperbolic) rotations whose q-analogs are often introduced via their
relations with QG $\uq$ in the spirit of \cite{FRT} or \cite{LS} (we discuss such an approach in
section \ref{sec:4}). In our construction of $q$-derivatives or braided tangent vector fields we use a
technique of Koszul type complexes which is applicable in a much more general context.  An
advantage of this approach becomes more evident when we pass to the quantum sphere
(hyperboloid) on which the analogs of tangent and cotangent bundles are realized as one-sided
projective modules. These modules play the crucial role in our defining the $q$-Maxwell
operator on a quantum hyperboloid.

The first attempt to construct an algebra of differential  forms on  quantum spheres is due
to P. Podle\`s \cite{P2, P3}. However, his approach, based on the Leibnitz rule,
leads to an algebra with non-classical dimensions on a generic quantum sphere. There is known
another approach to constructing the differential calculus on a quantum sphere which
uses $q$-analogs of complex coordinates $z$ and $\bar{z}$ in the spirit of
{\em projective} geometry (see for instance \cite{CHZ, SSV}).

Restricting ourselves to the {\em braided affine} geometry we disregard this approach
as well as that based on $C^*$-algebras which  involves  a big amount of functional
analysis.  We refer the reader to  the paper \cite{A} which
reviews some aspects of quantum geometry
essentially  based on the  RTT algebra and using
different tools of  noncommutative geometry due to A. Connes.
In contrast,  our approach is completely algebraical:
all "varieties" are braided affine, all representations
are finite dimensional. Mainly dealing with the REA algebra,
playing the crucial role  in braided  affine geometry, we  exhibit recent
 results and constructions which are covered neither by the review
 \cite{M3} nor by the monograph \cite{M4}.

To complete the Introduction, we would like to mention
 the approach based on the so-called
 Heisenberg double and going back to \cite{W}. It deals with the algebra
 $\K_q[GL(n)]\ot \bw_q(\LL)$  treated as a $q$-differential
 algebra on the group $GL(n)$. The difficulties of such a
 differential calculus
 are analyzed in numerous papers. We refer the reader to the recent
 paper \cite{IP2} where an application of this approach to the
 so-called quantum top is exhibited.

Our paper is organized as follows. In the next section we describe
semiclassical (i.e. Poisson) counterparts of quantum algebras related to
standard Hecke symmetries. In section \ref{sec:3} we recall a general method of
introducing the quantum matrix algebras (the REA and RTT algebras
included) which is based on the notion of  compatible braidings.
Also, there we reproduce the Cayley-Hamilton identity valid for the
quantum matrix algebras related to the Hecke symmetries of general
type and we discuss some special properties of the REA.

In section \ref{sec:4}
we compare the methods of constructing the (m)REA representation theory
in general and in the standard cases. Also, we
exhibit the coproduct which plays the central role in constructing
the aforementioned Schur-Weyl category.
In section \ref{sec:5} we present a treatment of the mREA as an
enveloping algebra of a {\em braided Lie algebra} and
discuss the problem of defining an involution  (conjugation) in this
"enveloping algebra".

In section \ref{sec:6} we compare different
ways of defining $q$-analogs of line bundles on
a quantum sphere (hyperboloid) and
the Chern-Connes index for them.
In section \ref{sec:7} we introduce some elements of
differential calculus on quantum algebras based on the Koszul type
complexes. In section \ref{sec:8} we apply this technique in order to define
the braided
analogs of partial derivatives and other vector fields. Then we introduce
the $q$-analogs of basic wave operators on the $q$-Minkowski
space and $q$-hyperboloid space algebras.

\section{Poisson counterparts of quantum varieties}
\label{sec:2}

In this section we consider Poisson counterparts of quantum
varieties. By a quantum variety algebra
we mean a $\Uqq$-covariant algebra which is a deformation of
(function algebra of) a usual homogeneous
$G$-space. Here $G$ is a classical matrix connected complex group (so,
$\K=\C$), $\gg$ is its Lie algebra,
and $\Uqq$ is the corresponding quantum group. Besides, we
assume that a triangular decomposition
of the Lie algebra $\gg$ is fixed.

We consider a homogeneous $G$-space of the form
$G/G_{\Gamma}$ where $G_{\Gamma}$ is {\em the Levi
subgroup}
 corresponding to a subset $\Gamma$ of the set of simple
 positive roots. The Lie algebra $\ggg=\Lie(G_{\Gamma})$
(called {\em the Levi subalgebra}) is generated by the Cartan
subalgebra $\hh\subset \gg$ and by the root vectors
$E_{\pm\alpha}$ corresponding to the roots $\pm \alpha$ for all
$\alpha\in \Gamma$. Besides, we assume
$\langle E_{\alpha},E_{-\alpha} \rangle=1$ where
$\langle\,\,,\,\rangle$ is the pairing defined by the Killing form.

The coset $G/G_{\Gamma}$ is isomorphic (as a $G$-set) to an
orbit $\O\subset \gg^*$ of a semisimple element
of the space $\gg^*$ ({\em a semisimple orbit} in short).
Consequently, we  identify the algebra
$\K[\O]$ and a subalgebra  in the algebra  $\K[G]$.

Consider an example of the above isomorphism for the matrix
group $G=SL(2)$
$$
L\in SL(2) \Leftrightarrow L=\left(\begin{array}{cc}
a&b\\
c&d \end{array}\right),
\quad ad-bc=1,\quad a,b,c,d\in \K.
$$
In this simple case the only Levy subgroup different from $SL(2)$
itself is the Cartan subgroup $H$ consisting of all
diagonal matrices with the unit determinant. We put $\sigma(L)=
(x,h,y)$ where $(x=ab,\, h=-ad-bc,\,y=-cd)$. Note that
$\sigma(L)=\sigma (Lg)$ for $g\in H$. Thus, the functions
$x,\,y,\,z$ are defined on the coset $SL(2)/H$. With respect to
the left $SL(2)$-action $L\to gL$, $ g\in SL(2)$, the space of
vectors $(x,h,y)$ becomes a spin 1 $SL(2)$-module. We
treat the vector $(x,h,y)$ as an element of $sl(2)^*\cong sl(2)$.
Since the quantity $h^2+4xy= (ad-bc)^2=1$ is stable under
the left $SL(2)$-action, the image of the map $\sigma$ belongs to
a hyperboloid passing through the point $(0,1,0)$.

Now, consider the Sklyanin bracket $\{\,\,,\,\}_{\K[G]}$ defined
on the space $\K[G]$ as follows
\be
\{f,g\}_{\K[G]}=\circ \left(\rho_l^{\ot 2}(r_-)(f\ot g)-\rho_r^{\ot 2}(r_-)
(f\ot g)\right), \quad f,g \in \K[G] \label{Skl} \ee
where $\circ$ is the ordinary pointwise product in the algebra
$\K[G]$, $r_-=(r_{12}-r_{21})\in \gg^{\ot 2}$ is the
skew-symmetrized classical $r$-matrix (up to a factor 1/2), and $\rho_l$
(resp., $\rho_r$) is the representation $\gg\to \Vect(G)$
of the algebra $\gg$ by the right-invariant (resp., left-invariant)
vector fields. Namely,
$$
\rho_l(X)f(a)=\partial_t f(e^{tX}\,a)|_{t=0}, \qquad \rho_r(X)f(a)=
\partial_t f(a\,e^{-tX})|_{t=0}
\quad \forall X\in \gg, \; a\in G.
$$

We assume the classical $r$-matrix $r$ to be chosen in such a
way that
$$
r_-= \sum_{\al\in \Om^+} (E_{\al}\ot E_{-\al}-E_{-\al}\ot E_{\al})
$$
where $\Om^+$ is the set of positive roots and the
symmetrized part $r_+=(r_{12}+r_{21})\in \gg^{\ot 2}$
is the {\em split Casimir element}
$$
r_+= \sum_{\al\in \Om^+}( E_{\al}\ot E_{-\al}+E_{-\al}\ot E_{\al})+...
$$
where we omit the terms belonging to ${\frak h}^{\ot 2}$.

It is well known that the Sklyanin bracket is a Poisson one and it is
compatible with the standard matrix coproduct in the
space $\K[G]$ in the sense of formula (\ref{below}) where we put $M=G$.

If $M$ is a homogeneous $G$-space, the space
$\K[M]$ can be equipped with a coaction
$\De:\K[M]\to \K[M]\ot \K[G]$. We say that a Poisson bracket
$\{\,\,,\,\}_{\K[M]}$ defined in the space $\K[M]$ is
$G$-{\em covariant} if
\be
\De\{f,g\}_{\K[M]}=\{\De(f),\De(g)\}_{\K[M]\ot \K[G]},\quad \forall \label{below}
f,g\in\K[M]
\ee
where the bracket in the right hand side is defined via the product
of the Poisson structures.

A Poisson bracket on the space $\K[M]$ is called $G$-{\em invariant} if the
operator $f(x)\to f(gx)$, $g\in G$
commutes with the bracket.
It is easy to see that the  sum of a  $G$-invariant bracket and a
$G$-covariant one is again a $G$-covariant bracket.

\begin{remark}{\rm
Hereafter we use the notion of a "bracket" in a large sense. By a bracket we mean an
operator
 $\{\,\,,\,\}: \K[M]^{\otimes 2}\to \K[M]$ which is  bilinear, skew-symmetric
and satisfying the Leibnitz rule (the Jacobi
identity is not required).
The notion of $G$-covariance (invariance) can be
extended to such ``brackets".  } \end{remark}

Consider an example of $G$-covariant Poisson bracket on a given
homogeneous $G$-space $\O$, namely, the
Sklyanin bracket reduced to the space $\K[\O]$. We denote this
bracket $\{\,\,,\,\}^{\O}$. The fact that this bracket
is well defined follows from \cite{DGS} (see also considerations
below).

Write the Sklyanin bracket as a difference of two terms
\be
\{\,\,,\,\}_{\K[G]}= \{\,\,,\,\}_{left}-\{\,\,,\,\}_{right}
\label{Sk-br-dif}
\ee
where
$$
\{f,g\}_{left}=\circ \rho_l^{\ot 2}(r_-)(f\ot g),\qquad
\{f,g\}_{right}=\circ \rho_r^{\ot 2}(r_-)(f\ot g).
$$

Being reduced to the space $\O$, the bracket $\{f,g\}_{left}$
becomes\footnote{If $r$ is an arbitrary element of $\bw^2(\gg)$ and $\gg$ is
represented in the space of functions on a homogeneous space $M$
by vector fields $\rho:\gg \to \Vect(M)$, the
expression $\{f,g\}_r= \circ \rho^{\ot 2}(r)(f\ot g)$ is often called the
$r$-matrix bracket. Thus, the bracket
$\{f,g\}^{\O}_{left}$ is just of this type.}
$$
\{f,g\}^{\O}_{left}= \sum_{\al\in \Om^+} \left(E_{\al}(f) E_{-\al}(g)-
E_{-\al}(f) E_{\al}(g)\right)
$$
where the elements $E_{\pm \al}$ are treated to be vector fields
naturally defined on the homogeneous space
$G/G_{\Gamma}$ which is isomorphic to $\O$.
As for the reduced bracket $\{\,\,,\,\}^\O_{right}$, it becomes
$G$-invariant. Nevertheless, in general the above
left and right brackets are not Poisson ones, since none of them
satisfies the Jacobi identity.

As was shown in \cite{GP}, on a semisimple orbit $\O\subset
\gg^*$ the bracket $\{\,\,,\,\}^{\O}_{left}$ becomes a Poisson
one iff the orbit $\O$ is symmetric\footnote{In \cite{GP} the orbits
possessing this property were called the $R$-matrix type
orbits. As was shown there, if $\O$ is an $R$-matrix type orbit,
then it is semisimple or nilpotent. For nilpotent orbits
a necessary and sufficient condition to be of the $R$-matrix type was
also found.}. In the case $G=SL(n)$ a semisimple orbit is
symmetric iff it is the orbit of a diagonal matrix with two different
eigenvalues. It is not difficult to see that
$\{\,\,,\,\}_{right}^{\O}$ and $\{\,\,,\,\}^{\O}_{left}$ become Poisson
brackets simultaneously. Indeed, the Schouten bracket
of $\{\,\,,\,\}^\O$ contains three terms: the Schouten bracket of
$\{\,\,,\,\}^{\O}_{left}$ with itself, that of
$\{\,\,,\,\}_{right}^{\O}$ with itself and that of these two brackets.
Since the latter Schouten bracket always vanishes
the former Schouten brackets vanish simultaneously (see
\cite{DGS} for detail).

Moreover, if $\O$ is a symmetric orbit, the bracket $\{f,g\}_{right}
^{\O}$ equals (up to a factor) the Kirillov-Kostant-Souriau
(KKS) one $\{f,g\}^\O_{KKS}$ which is the restriction of the linear
Poisson-Lie bracket defined on $\gg^*$ to the orbit $\O$.
This follows from the fact that in this case the space $\K[\O]$ is
multiplicity free  and therefore any two $G$-invariant
brackets are proportional to each other.

Finally, on a symmetric orbit $\O$ there exists a Poisson pencil
$$
\{\,\,,\,\}^\O_{a,b}=a\, \{\,\,,\,\}^\O_{KKS}+b\{\,\,,\,\}^{\O}_{left}.
$$
In \cite{DGK} an attempt was undertaken to quantize this Poisson
pencil on the ${\Bbb{CP}}^n$ type orbits (i.e.
those of minimal dimension) via the generalized Verma modules.
Another way of quantizing such a Poisson pencil in the
spirit of the affine algebraic geometry uses quotienting
 the REA. We exhibit this way below.

Let us go back to the example above. Consider the reduced
Sklyanin bracket on the hyperboloid
$\O=\{(x,h,y)|h^2+4xy=1\}$. Since this orbit is symmetric, the
bracket $\{f,g\}_{right}^{\O}$ equals (up to a factor)
the KKS one. In order to find this factor we compute
$$
\{x,y\}_{right}^{\O}=\{ab,-cd\}_{right}={-}(a^2d^2-b^2c^2)={-}(ad-bc)
(ad+bc)=h.
$$
Thus, it equals the KKS bracket.

Compute now the bracket $\{f,g\}^{\O}_{left}= (X(f) Y(g)-Y(f) X(g))
$ where $X,H,Y$ are the standard generators
of $sl(2)$ represented by the infinitesimal hyperbolic rotations:
$$
X=h\partial_y-2x\partial_h,\quad H=2x\partial_x-2y\partial_y,\quad
Y=-h\partial_x+2y\partial_h.
$$
We have
$$
\{h,x\}^{\O}_{left}=2hx,\quad \{h,y\}^{\O}_{left}=-2hy,\quad \{x,y\}^{\O}
_{left}=h^2.
$$
Thus, the Sklyanin bracket reduced to the hyperboloid $\O$ is
$$
\{h,x\}^{\O}=2(h-1)x,\quad \{h,y\}^{\O}=-2(h-1)y,\quad
\{x,y\}^{\O}=h(h-1).
$$
As for the above Poisson pencil we have
$$
\{h,x\}^\O_{a,b}=a(2x)+b(2hx),\quad \{h,y\}^\O_{a,b}=a(-2y)
+b(-2hy),\quad \{x,y\}^\O_{a,b}=ah+bh^2.
$$
It is easy to see that this bracket is Poisson on the whole space
$sl(2)^*$ too.

Our calculations are analogous to those from \cite{Sh} (see
appendix) where the compact form of the group, i.e. $SU(2)$,
was considered, but a concrete form of a given complex group is
somewhat pointless. In \cite{KRR} it was shown that the
Sklyanin bracket reduced to a semisimple orbit $\O$ is compatible
with the KKS one iff the orbit is symmetric. This entails
that if an orbit $\O$ is not symmetric then the brackets $\{\,\,,\,\}
_{right}^{\O}$ and $\{\,\,,\,\}^\O_{KKS}$ are not compatible
(i.e. their Schouten bracket does not vanish).

Note that the authors of \cite{KRR} also deal with compact forms of
groups and corresponding orbits. However, as we said
above this does not affect the result.

Nevertheless, if a given orbit $\O\subset \gg^*$ is not symmetric,
the family of $G$-covariant Poisson brackets becomes
larger.  Let $\O$ be a generic
orbit. This means that $G_{\Gamma}$ equals the Cartan subgroup
$H$. Let us describe the family of $G$-covariant Poisson
structures following \cite{DGS}.

Consider an element
$$
v=\sum_{\al\in\Om^+} c(\al)(E_{\al}\ot E_{-\al}-E_{-\al}\ot E_{\al})
$$
where $c(\al)$ is a $\K$-valued function on the set $\Om^+$ such
that if $\al+\beta\not=0$ then
\be
c(\al+\beta)={c(\al)c(\beta)+1 \over c(\al)+c(\beta)}.
\label{c1}
\ee
Here we assume that $c(\al)+c(\beta)\not=0$ for any $\al, \beta \in
\Om^+$ such that $\al+\beta\not=0$. This condition
forbids some low dimensional subvarieties. Now, we modify
the Sklyanin bracket (\ref{Sk-br-dif}) on $G$ by
replacing the bracket $\{\,\,,\,\}_{right}$ for the following one
$$
\{f,g\}_{right,c}=\circ (\rho_r^{\ot 2} (v)(f\ot g)).
$$

The bracket $\{f,g\}_{right,c}$ can be also reduced to the
homogeneous space $\O\cong G/H$ and the bracket
$$
\{\,\,,\,\}_{c}^{\O}=\{\,\,,\,\}^{\O}_{left}- \{\,\,,\,\}_{right, c}^{\O}
$$
is a Poisson one as was shown in \cite{DGS}.

Note that if
$\la:\Om^+\to \K$ is a linear function then $c(\la)= \coth(\la(\al))$ is
a solution of (\ref{c1}). Thus, the variety of $G$-covariant Poisson
brackets $\{\,\,,\,\}_{c}^{\O}$ can be identified with a vector space of
dimension $\dim \hh$ with exception of the aforementioned
low dimensional subvarieties. Observe that the reduced Sklyanin
bracket corresponds to the case $c(\al)=1$ for all $\al\in \Om^+$.

If an orbit $\O$ is not generic, i.e. $G_{\Gamma}\not=H$,
then different functions $\la$ can give rise to the same
Poisson structure on $\O$, i.e. in this case the map $c\to
\{\,\,,\,\}_{c}^{\O}$ is not injective.

Emphasize, that the bracket $\{\,\,,\,\}_{left}- \{\,\,,\,\}_{right, c}$ is not
Poisson one on $G$ (it becomes Poisson bracket
only after the reduction to $\O$). But it can be treated in terms of
the so-called classical dynamical $r$-matrices
(see \cite{L} and the references therein). Also, note that another
but equivalent approach to the classification problem
of $G$-covariant Poisson structures on homogeneous
$G$-spaces was developed in \cite{Ka}.

A way of formal quantization of the Poisson brackets $\{\,\,,\,\}_{c}
^{\O}$ was suggested in \cite{DGS} where as it is usual in the
framework of formal quantization, the problem of convergency of
series defining the deformed product is disregarded. It would
be interesting to express the bracket $\{\,\,,\,\}_{c}^{\O}$ in terms
of the coordinate algebra $\K[\O]$ and
to quantize it explicitly by presenting the result in the spirit of affine
algebraic geometry or to show that it is not possible.

Meanwhile, on the whole vector space $gl(n)^*$ there exists a
$GL(n)$-covariant Poisson structure which is the semiclassical
counterpart of the REA (for this reason we denote the
corresponding bracket $\{\,\,,\,\}_{REA}$). Moreover, as shown in
\cite{D}
the bracket $\{\,\,,\,\}_{REA}$ can be restricted to any orbit
$\O\subset gl(n)^*$ giving rise to the $GL(n)$-covariant Poisson
structure on the coordinate algebra $\K[\O]$\footnote{In fact there
are two structures: one of them is considered below,
another one is associated with the second type REA in a similar way.}.

Besides, the bracket $\{\,\,,\,\}_{REA}$ is compatible with the linear
Poisson-Lie bracket $\{\,\,,\,\}_{PL}$ also defined on $gl(n)^*$. By direct
computations it is easy to see, that the the semiclassical counterpart of
the two parameter mREA $\lhq$ is just the Poisson pencil (\ref{PP})
generated by these two brackets. Thus,  this Poisson pencil can be quantized
 explicitly in the spirit of affine algebraic geometry.
 By treating  its quantum counterpart via the algebra $\lhq$  we avoid
 the convergency problem of the deformation quantization scheme.
The restrictions of the above Poisson pencil to semisimple orbits can be also
quantized in a similar manner. As a result of such a quantization we get
 quotients of the mREA discussed below.

The explicit form of the bracket $\{\,\,,\,\}_{REA}$ reads (see \cite{Mu1, GPS3})
\be
\{f,g\}_{REA}= \circ\, r_+^{\rm l,r}(f\ot g) -
\circ\, r_+^{\rm r,l}(f\ot g) - \circ \,r_-^{\rm ad,ad}(f\ot g),\qquad
\forall\, f,g\in {\Bbb K}[gl(m)^*].
\label{brr}
\ee
Here
$$
r_-=\sum_{i<j} e_i^j\ot e_j^i-e^i_j\ot e^j_i,\qquad
r_+=\sum_{i,j} e_i^j\ot e_j^i,
$$
the elements $e_i^j$ form the standard basis of the Lie algebra
$gl(n)$
$$
[e_i^j,e_k^r] = \delta_k^j\,e_i^r - \delta_i^r\,e_k^j.
$$
and the element $r_+\in gl(n)^{\ot 2}$ is the $gl(n)$  {\em split Casimir element}.
Superscripts indicate the type of the action. Namely,
$$
r_-^{\rm ad,ad}(f\ot g)=\sum_{i<j}\left({\rm ad}\,{e_i^j}(f)\ot {\rm ad}\,
{e_j^i} (g)-
{\rm ad}\,{e_j^i}(f)\ot {\rm ad}\, {e_i^j} (g)\right),
$$
$$
r_+^{\rm l,r}(f\ot g)=\sum_{i,j}\rho_l(e_i^j)f \otimes \rho_r(e^i_j)
g,\qquad
r_+^{\rm r,l}(f\ot g)=\sum_{i,j}\rho_l(e_i^j)g \otimes \rho_r(e^i_j)f.
$$
Here the notation ${\rm ad}\,{e_i^j}$ stands for the vector field
acting on $l_k^r\in \K[GL(n)]$ by the rule
${\rm ad}\,{e_i^j}(l_k^r)=\delta_k^j\,l_i^r - \delta_i^r\,l_k^j$, whereas
the notations $\rho_l,\,\, \rho_r$ have
the same meaning as above.

Let us list the basic properties of the bracket $\{\,\,,\,\}_{REA}$.
First of all, it is $GL(n)$-covariant. More precisely,
it is a linear combination of the $r$-matrix bracket (the third term in
(\ref{brr})) and $GL(n)$-invariant bracket (the two
first terms in (\ref{brr})). Second, it can be restricted to the
subspace $sl(n)^*$ since for all $f\in \K[GL(n)]$ we have
$\{f,\ell\}_{REA}=0$, where $\ell=\sum l_i^i$. The restricted bracket
$\{\,\,,\,\}_{sl-REA}$ is $SL(n)$-covariant.

\begin{remark}{\rm
The first two terms in (\ref{brr}) can be interpreted in the following way \cite{G3}.
The space $sl(n)^{\ot 2}$ considered as an $sl(n)$-module can be
decomposed into a direct sum of irreducible
$sl(n)$-modules. For $n>2$ this decomposition contains two
components isomorphic to the algebra $sl(n)$. One of these
components belongs to the skew-symmetric part
$\bigwedge^2(sl(n))$, the other one belongs to the symmetric part
$\Sym^{(2)}(sl(n))$.
So, there exists unique (up to a factor) nontrivial $sl(n)$-morphism
$\sigma: \bigwedge^2(sl(n))\to
\Sym^{(2)}(sl(n))$ (it identifies the corresponding $sl(n)$-components
and kills all the others). On extending the morphism
$\sigma$ up to the algebra $\Sym(sl(n))$ via the Leibnitz rule, we
get an $SL(n)$-invariant bracket. (Note that for a Lie algebra
$\gg$ of the series $B_n$, $C_n$ or $D_n$ no similar bracket exists
since the space $\gg^{\ot 2}$ is multiplicity free.)
At some specific value of the factor of the morphism $\sigma$, the
corresponding bracket coincides with that defined by
the two first terms of (\ref{brr}) restricted to $sl(n)^*$.}
\end{remark}

The bracket $\{\,\,,\,\}_{sl-REA}$ can be restricted to any
orbit $\O\subset sl(n)^*$ (see \cite{D}). If such an orbit is
semisimple, it is closed. Let $I_{\O}$ be an ideal of functions
vanishing on such an orbit $\O$. Then $\{f,g\}_{sl-REA}\in I_{\O}$ for any
$f\in \K[sl(n)]^*\cong \Sym(sl(n))$ and $g\in I_{\O}$. The
quantization of the Poisson pencil generated by the KKS bracket on $\O$ and the
bracket $\{\,\,,\,\}_{sl-REA}$ restricted to the orbit can be done via a
proper quotienting the standard mREA $\lhq$ in the spirit of affine
algebraic geometry. Thus, we get the braided coordinate algebra
$\K_q[\O]$ which is a deformation of the classical algebra $\K[\O]$
whereas
a quantization of the restricted bracket $\{\,\,,\,\}_{sl-REA}$ alone
can be realized as a quotient of the algebra $\lq$.

If the  orbit $\O$ is generic, the ideal $I_{\O}$ is generated by the
elements $\Tr L^k-a_k$ with appropriate $a_k\in \K$,
$2\le k\le n$ (since $\Tr(L)=0$). Here $L$ is the classical analog of
the matrix $L$ coming in the definition of the algebra $\lq$,
i.e. it satisfies the defining relation of the REA but with $R$
replaced by the usual flip. Then we get the algebra $\K_q[\O]$ as a
quotient of the algebra $\lhq$ or $\lq$ over an ideal generated by the elements
$\Tr_R (L^k)-a_k$, $2\le k\le n$, where $\Tr_R$ is quantum trace
introduced below. If an orbit $\O$ is not generic the problem of finding
the "quantum coordinate algebra" $\K_q[\O]$ is more subtle (see \cite{DM}).

Emphasize that for other simple Lie algebras there exist orbits on
which the bracket $\{\,\,,\,\}_{c}^{\O}$ is compatible
with the KKS one for some $c$ (such orbits were called {\em
good} in \cite{DGS}). However, these brackets and the
corresponding
pencils were quantized in \cite{DGS} via the deformation
quantization scheme only. It is not known whether
it is possible to define the corresponding quantum algebra by a finite
set of polynomial relations among their generators (at least for the matrix
groups).

Let us give a short summary of the considerations above. If an
orbit $\O\subset sl(n)^*$ is symmetric, the family of $SL(n)$-covariant
brackets is just the Poisson pencil generated by the KKS bracket
and the $r$-matrix one. Quantization of any Poisson bracket from
this pencil can be always realized as a proper quotient of the
algebra $\lhq$. Thus, we get a
{\em braided affine variety}. If a semisimple orbit $\O$ is not
symmetric, the quantization method depends on a given
bracket on $\K[\O]$. However, the previous method is still valid,
provided such a bracket is a linear combination of the KKS bracket
and that $\{\,\,,\,\}_{sl-REA}$.

\section{Quantum matrix algebras}
\label{sec:3}

The REA and the RTT algebra defined above are
particular examples of the so-called {\it quantum matrix algebras}. First
examples of these appeared in works \cite{FM} and \cite{H}. A general definition and extensive studying
of the algebraic structure were given in \cite{IOP}.

The definition of the quantum matrix algebra is based on the pair of
the compatible braidings
which are defined below.
\begin{definition}
\label{def:com-p}
{\rm
Let $V$ be a finite dimensional vector space over the ground field
$\Bbb K$,
$\dim V = N$, and let $R,F\in {\rm End(V^{\otimes 2})}$
be two invertible operators.
An ordered pair $\{R,F\}$ is called a pair of compatible braidings if
the following conditions are
satisfied
\begin{itemize}
\item
Both operators $R$ and $F$ obey the quantum Yang-Baxter
equation (\ref{YBE}).
\item
The operators $R$ and $F$ satisfy the compatibility conditions
\be
R_{12}F_{23}F_{12} = F_{23}F_{12}R_{23}\, ,
\qquad
F_{12}F_{23}R_{12} = R_{23}F_{12}F_{23}\, .
\ee
\end{itemize}
}
\end{definition}

Given a pair of compatible braidings $\{R,F\}$ we shall additionally
assume
the braidings to be {strictly skew-invertible}.
\begin{definition}\label{def:sk-inv}
{\rm
An operator $R$ is said to be {\em skew-invertible} if there exists an operator
$\Psi^R\in {\rm End}(V^{\otimes 2})$ such that
\be
{\rm Tr}_{(2)}R_{12}\Psi_{23}^R = P_{13} =
{\rm Tr}_{(2)}\Psi_{12}^R R_{23}\, ,
\label{def:psi}
\ee
where the subscript in the notation of the trace shows the index
of the space $V$,
where the trace is applied (the enumeration of the factors
in the tensor product is taken as follows
$V^{\otimes k} = V_1\otimes V_2\otimes\dots \otimes V_k$).
The skew-invertible operator $R$ is {\it strictly skew-invertible} if the
operator
\be
C={\rm Tr}_{(2)}\Psi^R_{12}
\label{def:C}
\ee
 is invertible.}
\end{definition}
Let us note, that for a strictly skew-invertible braiding $R$ the
operator
\be
B={\rm Tr}_{(1)}\Psi^R_{12}
\label{def:B}
\ee
is also invertible \cite{O}.

As a direct consequence of definitions of $B$ and $C$ we find
\be
\Tr_{(1)}B_1R_{12} = {\rm Id} =\Tr_{(2)}C_2R_{12}.
\label{BRC}
\ee
The operators $B$ and $C$ play the crucial role for the structure
of the
quantum matrix algebras and, in particular, for the representation
theory of REA.

With the matrix $C$ one defines the {\em $R$-trace}
operation ${\rm Tr}_R:{\rm Mat}_N(W)\rightarrow W$
\be
{\rm Tr}_R(X) = \sum_{i,j=1}^NC_i^jX_j^i,\qquad X\in
{\rm Mat}_N(W),
\label{def:R-tr}
\ee
where $W$ is any linear space.
\begin{remark}
{\rm Note, that if we work with the second form of
REA, the matrix $C$ in the above definition of the $R$-trace should be changed
for the matrix $B$.}
\end{remark}

Now we are ready to give a general definition of a quantum matrix
algebra.
\begin{definition}\label{def:QMA}
{\bf \cite{IOP}} \hspace{3pt}{\rm
Given a compatible pair $\{R,F\}$ of strictly skew-invertible
braidings,
the {\em quantum matrix algebra } (QMA) ${\cal M}(R,F)$
is defined to be a unital associative algebra generated
by $N^2$ components of the matrix $M=\|M_i^{\,j}\|_{1\le i,j\le N}$ subject
to the relations
\be
\label{def:QM}
R_{12}M_{\overline 1}M_{\overline 2} = M_{\overline 1}
M_{\overline 2} R_{12}\, ,
\ee
where we use the following notation
\be
M_{\overline 1} = M_1 = M\otimes {\rm Id}, \quad M_{\overline 2} =
F_{12} M_{\overline 1} F^{-1}_{12}
\label{copy}
\ee
for the copies of the matrix $M$.}
\end{definition}

The defining relations (\ref{def:QM}) then imply the same type
relations for
consecutive pairs of the copies of $M$ (see \cite{IOP})
\be
R_{k,k+1}\,M_{\overline k}M_{\overline{k+1}} =
M_{\overline k}M_{\overline{k+1}}\, R_{k,k+1}\,,
\nonumber
\ee
where $M_{\overline {k+1}}=F_{k,k+1}M_{\overline k}F^{-1}_{k,k+1}
$.

\begin{remark}{\rm As can be easily seen, the pairs $\{R, P\}$ ($P$
being
a permutation matrix) and
$\{R,R\}$ are compatible pairs of braidings in the sense of
Definition
\ref{def:com-p}. The first of them defines the RTT algebra while
the second
one gives rise to the REA. Besides these evident examples, there exist
large multiparametric families of compatible braidings
leading to
other quantum matrix algebras.}
\end{remark}

From now on we shall be interested in a subfamily of QMAs
characterized by the
property that in the pair $\{R,F\}$ the first component $R$ is a
skew-invertible
Hecke symmetry. Recall, that it means that skew-invertible matrix
$R$ additionally obeys
 the quadratic Hecke condition
\be
(R-q\, {\rm Id})(R+q^{-1}\,{\rm Id})=0\,, \quad q\in {\Bbb C}\setminus
0\,,
\label{Hec}
\ee
where the numerical parameter $q$ is either equal to $1$ or is not
a root of unity
of any order: $q^k\not=1$, $\forall \,k\in {\Bbb Z}_+$.

With any Hecke symmetry $R$ one can associate an ordered pair
of integers $(m|n)$ called the
{\it bi-rank}\footnote{This notion was defined
in \cite{GPS3} as follows.
The Poincar\'e-Hilbert series of the algebra $\bw_q(V)$
where $V$ is a vector space endowed with a Hecke
symmetry $R:\vv\to \vv$, is always
a rational function. Let us assume it to be uncancellable.
Then $m$ (resp., $n$) is the degree of its numerator
(resp., denominator).} of the symmetry $R$. A
Hecke symmetry with
the bi-rank $(m|n)$ can be considered as a generalization of
the super-flip on
the vector super-space $V_{(m|n)}$. Therefore, we  call this
Hecke symmetry (and
the corresponding QMA) $GL(m|n)$ type braiding. Note,
that any skew-invertible Hecke symmetry of the bi-rank $(m|n)$ is
automatically strictly
skew-invertible since it can be proved (see \cite{GPS3}, corollary
8) that
\be
B\cdot C = q^{2(n-m)}{\rm Id}\,.
\label{B.C}
\ee

A Hecke symmetry $R$ realizes local (or $R$-matrix)
representations of the $A_{n-1}$ series
Hecke algebras ${\cal H}_n(q)$. The detailed treatment of the structure
of the Hecke algebras and
their $R$-matrix representations can be found in the review
\cite{OP1} and references therein.

Describe now the main properties of a general $GL(m|n)$ type
QMA. First of all we can extract the commutative {\it characteristic
subalgebra} Char($M$) which is formed by the linear span of the
QMA elements of the following form
\be
\label{char}
x(h_k) = \Trr{1\dots k}(M_{\overline 1}\dots M_{\overline
k}\,\rho_R(h_k))\quad k =1,2,\dots ,
\ee
where $h_k$ runs over all elements of the Hecke algebra ${\cal H}
_k(q)$ and
$\rho_R: {\cal H}_k(q)\rightarrow {\rm End}(V^{\otimes k})$ is the
$R$-matrix representation
of the Hecke algebra associated to a given Hecke symmetry $R$.
 The symbol $\Trr{1\dots k}$
means that we apply  the $R$-trace over the spaces from the first to the
$k$-th ones.

Among  elements of the characteristic subalgebra we distinguish
the set of "power sums"
$p_k(M)$ and Schur symmetric functions $s_\lambda(M)$ defined
respectively by the
relations
\begin{equation}
p_0(M): = {\rm Tr}_R({\rm Id}), \quad
p_k(M)= \Trr{1\dots k}(M_{\overline 1}\dots
M_{\overline k}\,R_{k-1}\dots R_1)\,,
\label{st-sm}
\end{equation}
where $R_i=R_{i,i+1}$ and
\be
s_0(M) = 1,\quad
s_\lambda(M) = \Trr{1\dots k}(M_{\overline 1}\dots
M_{\overline k}\, \rho_R(E^\lambda_\alpha))\,.
\label{f-sh}
\ee
Here $\lambda\vdash k$ is an arbitrary partition of the given
integer $k\in {\Bbb Z}_+$
and $E^\lambda_\alpha$ stands for a primitive idempotent of the
Hecke algebra ${\cal H}_k(q)$
parameterized by the Young tableau $(\lambda;\alpha)$
corresponding to the partition
$\lambda$. The right hand side of (\ref{f-sh}) does not actually
depend on the index
$\alpha$ of the primitive idempotent. It can be proved that the
elements $p_k(M)$
generate the characteristic subalgebra \cite{IOP} and the elements
$s_\lambda(M)$ form
its linear basis \cite{GPS1}.

Besides, we define {\em the $k$-th powers} $M^{\overline k}$ of
the matrix $M$
by the following rule:
\be
M^{\overline 0}= {\rm Id},\quad
M^{\overline 1}= M_1,\quad M^{\overline k}=
\Trr{2\dots k}(M_{\overline 1}\dots M_{\overline k}
\,R_{k-1} \dots R_1), \quad k=2,3,\dots .
\label{m-k}
\ee

Comparing the definitions (\ref{st-sm}) and (\ref{m-k}) we find the
relation
$p_k(M) = {\rm Tr}_R(M^{\overline k})$ which is similar to the
classical one. But in general, the
usual matrix product of quantum matrices does not possess the
classical property, that is
$$
M^{\overline k}\cdot M^{\overline p}\not=M^{\overline {k+p}}.
$$
Nevertheless, it is possible to define an associative multiplication
$\star$
of quantum matrices which
satisfies the following remarkable properties (see \cite{OP2})
$$
M^{\overline{k+1}} = M\star M^{\overline{k}} = M^{\overline{k}}\star
M,\qquad ({\rm Id}\,x(M))\star M^{\overline{k}} = M^{\overline{k}}\star
(x(M)\,{\rm Id})\,,
$$
where $x(M)$ is an arbitrary element of the characteristic
subalgebra.

The definition (\ref{m-k}) is justified by the fundamental property of
any
$GL(m|n)$ type QMA. Namely, in each such algebra there
exists a
polynomial identity in powers $M^{\overline k}$ of order $m+n$
generalizing the classical Cayley-Hamilton identity valid for the numerical
matrices. The coefficients of the "quantum Cayley-Hamilton
polynomial" are some linear combinations of the
Schur functions $s_\lambda(M)$. Moreover, the quantum
Cayley-Hamilton polynomial can be presented in a factorized form
as the $\star$-product of two polynomials of the orders $m$ and $n$.
This fact enables us to introduce the notion of the "spectrum" of the
quantum (super) matrix and distinguish the "odd" and "even" eigenvalues.
The detailed treatment of these questions with complete proofs is given
in \cite{GPS1, GPS2}.

Now we give some explicit formulas for the case of ${\cal M}(R,R)$
$GL(m|n)$ type QMA (hereafter we restrict ourselves to this case).
 As we already mentioned above, this is
nothing but the REA. In this case the complicated $\star$-product of
quantum matrices reduces to the usual matrix product and the
characteristic subalgebra ${\rm Char}(L)$ is {\it central} in the REA,
that is
$$
L^{\overline k}\equiv L^k= L \cdot L^{k-1} = L^{k-1}\cdot L\,,
\qquad x(L)L = L x(L),\quad \forall\,x(L)\in {\rm Char}(L).
$$

On introducing  a special notation $[m|n]^k_r$ for the partition
$((n+1)^k,n^{m-k},r)\vdash (mn+k+r)$ we can write the Cayley-Hamilton
identity as follows \cite{GPS1}
$$
\sum_{i=0}^{m+n}L^{m+n-i}\sum_{k=\max\{0,i-n\}}^{\min\{i,m\}}
(-1)^kq^{2k-i}s_{[m|n]^k_{i-k}}(L) \equiv 0.
$$
This identity can be presented in a remarkable factorized form.
Multiplying
the above formula by $s_{[m|n]}(L)$, we can rewrite the result
as follows
\cite{GPS2}
$$
\Big(\sum_{k=0}^m (-q)^k\, L^{m-k}
s_{[m|n]^k_0}(L)
\Big) \cdot \Big(\sum_{r=0}^n q^{-r}\, L^{n-r}
s_{[m|n]_r^0}(L)\Big)\equiv 0\, .
\label{factor-ch}
$$
Hereafter the notation $[m|n]=[m|n]^0_0$ stand for the $m\times n$ rectangle.

It is useful to introduce the parametrization for the
Schur functions as the homomorphic map from the characteristic
subalgebra ${\rm Char}(L)$ into the algebra ${\Bbb K}[\mu,\nu]$
of polynomials in commuting variables $\mu= \{\mu_i \}_{1\le i\le m}
$
and $\nu= \{\nu_j \}_{1\le j\le n}$ (see \cite{GPS2})
\begin{eqnarray}
\textstyle\frac{s_{[m|n]^k}(L)}{s_{[m|n]}(L)}
&\mapsto &
{\textstyle \frac{s_{[m|n]^k}(\mu,\nu)}{s_{[m|n]}(\mu,\nu)}}
\, =\!\!\!\!
\sum_{1\le i_1<\dots <i_m\le m}\!\!\!\! q^{-k}\mu_{i_1}\dots
\mu_{i_k} =e_k(q^{-1}\mu)\, , \quad 1\le k\le m\, ,
\label{def:mu}
\\[1mm]
\textstyle \frac{s_{[m|n]_r}(L)}{s_{[m|n]}(L)}
&\mapsto&
{\textstyle \frac{s_{[m|n]_r}(\mu,\nu)}{s_{[m|n]}(\mu,\nu)}}
\, =\!\!\!\!
\sum_{1\le j_1<\dots <j_r\le n}\!\!\!\!(-q)^r\nu_{j_1}\dots\nu_{j_r}
=e_r(-q\nu)\, ,\quad 1\le r\le n\, .
\label{def:nu}
\end{eqnarray}
Here $e_k(\cdot)$ denotes the elementary symmetric polynomial
in finitely many variables --- the arguments of $e_k(\cdot)$.

Note that  this
parametrization is noncontradictory if we demand the Schur
function $s_{[m|n]}(L)$ to be invertible and the ratios of the Schur functions
in the left hand side of (\ref{def:mu}) and (\ref{def:nu})
to be algebraically independent. Note, that the indeterminates
$\mu_i$ and $\nu_j$ can be treated as elements of an algebraic
extension
of the field of fractions of the central characteristic
subalgebra\footnote{As
follows from the results of \cite{GPS2}, the characteristic
subalgebra of REA
considered as a ring is an integer domain and the field of fractions
can be
correctly defined.}.

With this parametrization the Cayley-Hamilton identity takes the
form
$$
(s_{[m|n]}(L))^2 \prod_{i=1}^m
(L - \mu_i {\rm Id}) \cdot
\prod_{j=1}^n
( L - \nu_j {\rm Id})
\equiv 0\, .
\label{factor2-ch}
$$
Recall, that since we are now dealing with the REA
all  matrix products have the usual sense.
The above totally factorized form of the Cayley-Hamilton identity
justifies an interpretation of the indeterminates
$\{\mu_i\}$ and $\{\nu_j\}$ as, respectively, {\em "even"}
and {\it "odd"} eigenvalues of the quantum super-matrix $L$.
\begin{remark}{\rm
As was shown in \cite{GPS2}, the map (\ref{def:mu})--(\ref{def:nu})
allows us to parameterize all the elements of the characteristic
subalgebra
in terms of $\mu$ and $\nu$. In particular,
$$
s_{[m|n]}(L)\, \mapsto\, s_{[m|n]}(\mu,\nu) = \prod_{i=1}^m\prod_{j=
1}^n
\left(q^{-1}\mu_i - q\nu_j\right).
$$
Therefore, the invertibility of the Schur function $s_{[m|n]}(L)$ in
terms of the
spectrum of $L$ is equivalent to invertibility of all factors
$\left(q^{-1}\mu_i - q\nu_j\right)$.
}
\end{remark}

Now we  describe "$sl$-reduction" of the REA and interrelations between this algebra and the RTT one.
Let us rewrite  the defining relations of the REA $\lq$ in the following form
\be
R_{12}L_1R_{12}L_1R_{12}^{-1} - L_1R_{12}L_1 = 0\quad
{\rm or}\quad L_1R_{12}L_1R_{12}^{-1} - R_{12}^{-1}L_1R_{12}L_1= 0\,.
\label{def:REA}
\ee
Here as usual $L= \|L^{\,j}_i\|_{1\le i,j\le N}$ ($N=\dim V$) and $R$ is a skew-invertible Hecke
symmetry of the bi-rank $(m|n)$.

One can easily prove the centrality of elements
$p_k(L) = {\rm Tr}_R(L^k)$ which generates the characteristic
subalgebra of the REA.
Indeed, consecutively multiplying the defining relation
(\ref{RE}) by the matrix $L_1$ from the left we get the following relations
$$
R_{12}L_1R_{12}L_1^k -L_1^k R_{12}L_1R_{12} = 0,
\quad \Leftrightarrow \quad L_1R_{12}L_1^kR_{12}^{-1} -
R_{12}^{-1}L_1^k R_{12}L_1= 0,\quad
\forall\,k\in {\Bbb Z}_+.
$$
Then, by applying the $R$-trace $\Trr{2}$ in the second space from the last
relation and using the
property
$$
\Trr{2}(R_{12}^{\pm}X_1R_{12}^{\mp}) = {\rm Id}_1{\rm Tr}_R(X)
$$
we find $p_k(L)L-Lp_k(L) = 0$. Emphasize once more, that this proof uses only the
skew-invertibility
of $R$ and is valid for any $GL(m|n)$ type REA.

Now, apply a linear shift to generators of a $GL(m|n)$ type REA:
\be
L^{\,j}_i\rightarrow L^{\,j}_i - \frac{\hbar}{q-q^{-1}}\delta^j_i\, 1_{\cal L}
\label{shift-1}
\ee
with a nonzero numerical parameter $\hbar$ (we retain the same notations for the
new generators and for the corresponding matrix $L=\|L^{\,j}_i\|$). Using the Hecke
condition on $R$, for the new generators we get relations (\ref{mREA0}) of the
mREA ${\cal L}(R_q,\hbar)$
$$
R_{12}L_1R_{12}L_1 - L_1R_{12}L_1R_{12} = \hbar (R_{12}L_1
- L_1R_{12}).
$$
Recall that the semiclassical  counterpart of the standard mREA is a Poisson pencil  considered in section
\ref{sec:2}.

If the bi-rank components $m\not=n$ , then the $R$-trace of the
unit matrix is nonzero in virtue of the following fomula
(see \cite{GPS3})
$$
{\rm Tr}_R({\rm Id}) = {\rm Tr}(C) = q^{n-m}(m-n)_q,\qquad k_q=
\frac{q^k-q^{-k}}{q-q^{-1}}.
$$
In this case we are able to make an "$sl$-reduction" of the mREA (or REA).
To this end we pass from the set of generators $\{L^{\,j}_i\}$ to
another set in which
the central element $\ell = p_1(L) = {\rm Tr}_R(L)$ is chosen as
one of the generators.
Let us extract the $R$-traceless part of the matrix $L$:
\be
F=L-\frac{\ell}{{\rm Tr}_R({\rm Id})}\, {\rm Id},\quad \ell= {\rm Tr}
_R(L).
\label{shift}
\ee
Then the defining relations for the mREA can be rewritten in terms
of generators
$F^{\,j}_i$ and $\ell$ as follows
\be
R_{12}F_1R_{12}F_1 - F_1R_{12}F_1R_{12} = \Bigl(\hbar 1_{\cal L}-
\frac{ (q-q^{-1})}{{\rm Tr}_R({\rm Id})}\,\ell\Bigr)
 (R_{12}F_1 - F_1R_{12}),\quad F\ell=\ell F.
\label{F-l}
\ee
Here the generators $F_i^{\,j}$ are linearly dependent since ${\rm
Tr}_R(F) = 0$.

We see, that contrary to the classical case, the set of generators
$F_i^{\,j}$ does not generate a subalgebra in mREA.
But the centrality of $\ell$ allows us to consider a quotient algebra
\be
{\cal SL}(R_q,\hbar)={\cal L}(R_q,\hbar)/\langle\ell\rangle.
\label{SL-red}
\ee
The defining relations for the ${\cal SL}(R_q,\hbar)$ generators $F_i^j$ (we keep the same
notation for them) can be obtained from (\ref{F-l})
by setting $\ell = 0$. The passage to the quotient ${\cal SL}
(R_q,\hbar)$ (resp. $\slq$)
will be referred to as an $sl$-reduction of the mREA ${\cal L}
(R_q,\hbar)$ (resp., $\lq$).

At last, we define {\it the truncated} mREA ${\tilde {\cal L}}
(R_q,\hbar)$
as an associative algebra generated by the elements $F_i^{\,j}$
and $\ell$
subject to the rules:
$$
R_{12}F_1R_{12}F_1 - F_1R_{12}F_1R_{12} = \hbar
 (R_{12}F_1 - F_1R_{12}),\quad F\ell=\ell F,\quad {\rm Tr}_R(F) =
 0.
$$

We complete this section by discussing a covariance property of the REA with respect to
the adjoint coaction of the RTT algebra. The coaction is based on the Hopf algebra structure
which can be defined in the RTT algebra under the additional condition that the element
$\det_qT$ is central (see Introduction). This condition is necessary for definition of the antipodal
map. Note, that in the standard case this coaction allows us to
endow the REA with a covariant $U_q(sl(N))$-module structure since the quantum group is restricted
dual to the standard RTT algebra quotiented by the condition $\det_qT =1$.

The coproduct in the RTT algebra reads  $\Delta(T) = T\stackrel{.}{\otimes}T$ \cite{FRT}.
Hereafter, the notation $A\stackrel{.}{\otimes}B$ means the  matrix product
of two equal size matrices $A$ and $B$ where their matrix elements are not multiplied
but tensorized. Thus, matrix elements of the matrix $A\stackrel{.}{\otimes}B$ are
$$
(A\stackrel{.}{\otimes}B)_i^{\,j} = A_i^{\,k}\otimes B_k^{\,j}.
$$

Consequently, all the algebras considered above (${\cal L }(R_q)$, ${\cal L }
(R_q,\hbar)$,
${\tilde {\cal L }}(R_q,\hbar)$ and ${\cal SL }(R_q,\hbar)$) can be given
a structure of the (right) adjoint comodule over the RTT algebra.
The  adjoint coaction $\delta_r$ on the generators can be presented  as follows
$$
\delta_r( A_i^{\,j}) =
\sum_{k,p=1}^NA_k^{\,p}\otimes T_i^{\,k}S(T_p^{\,j}),\quad{\rm or}
\quad
\delta_r( A) = TAS(T)\,,
$$
where $A$ stands for $L$ or $F$ and $S(\cdot)$ denotes the antipodal map.
The products in the above algebras are covariant with respect to this coaction, i.e.,  the
multiplication in these algebras commutes with the coaction.

Thus, in the standard case all these algebras  can be endowed with an action of the QG
$\Uq$ so that their products are coordinated with this action in the sense of formula (\ref{coo}).

\section{Elements of the mREA representation theory}
\label{sec:4}

In this section we reproduce some elements of representation theory of the REA related
to a skew-invertible Hecke symmetry of general type and compare our approach with that arising from
\cite{FRT} and \cite{LS} which is valid for the standard REA.

Our method of constructing the representation theory of the REA is based on
treating the space $\LL = {\rm span}(L_i^{\, j})$ as a space of endomorphisms of the basic (fundamental) space
$V$. Indeed, as was shown in \cite{GPS3}, the dual spaces $V$ and $V^*$ possess the
basis sets $\{x_i\}$ and $\{x^j\}$ in which the right and left pairing are uniquely (up to a nonzero factor)
defined by the formulae
\be
V\ot V^*\to \K:\quad x_i\ot x^j \mapsto \de_i^j\qquad {\rm and}\qquad
V^*\ot V\to \K:\quad x^j\ot x_i \mapsto B_i^{\,j},
\label{l-r-pairing}
\ee
where $\|B_i^j\|$ is the matrix corresponding to the operator (\ref{def:B}) in the basis $\{x_i\}$.
Recall that this matrix is invertible for any skew-invertible Hecke symmetry. In order to construct
the Schur-Weyl category of the REA representations we need only the skew-invertibility
of a Hecke symmetry.

The crucial point is that the map $\rho_V: \LL\to \End(V)$ defined by the rule
\be
\rho_V(L_i^j)\triangleright x_k= x_i B^j_k
\label{B-rep}
\ee
realizes an irreducible representation of the algebra $\loq$ (the symbol $\triangleright$ stands
for the action of a linear operator). This action together with the pairing $\langle x^j, x_i\rangle=
B_i^{\,j}$ is the motivation for treatment of $L_i^j$ as the element $x_i\ot x^j$.

If $R$ is the usual flip (therefore $B_i^j=\de_i^j$) we get the basic (fundamental) representation
of the algebra $U(gl(n))$. If $R$ is a super-flip then $B$ is the parity operator: $B(x)=(-1)^{\bar x} x$ where
$x$ is a homogeneous element and ${\bar x}$ is its parity.

Let us consider an example, namely the REA related to the $\uq$ Hecke symmetry. In the basis
$\{x_i\otimes x_j\}\in V^{\otimes 2}$ we get the following matrix
$$
R = \left(
\matrix{
q&0&0&0\cr
0&q-q^{-1}&1&0\cr
0&1&0&0\cr
0&0&0&q}\right).
$$
Then the matrix $B$ becomes $\diag(\qq, q^{-3})$. Let
$L=\left(\begin{array}{cc}
a&b\\
c&d
\end{array}\right)$
be the matrix of the mREA generators. Explicitly, the corresponding mREA $\lhq$ is
defined by the following system
\be
\begin{array}{l@{\hspace{20mm}}l}
q ab - q^{-1}ba =\h\, b &q( bc - cb) -(q-\qq) a(d-a)=\h\,(a-d)\\
q ca - q^{-1}ac =\h\, c & q(cd - dc) - (q-\qq)c a=-\h\, c\\
ad-da=0 & q(db - bd) - (q-\qq) a b=-\h\, b.
\end{array}
\label{sys1}
\ee

We present this system in terms of another set of generators $\{\ell,\,h,\,b,\,c\}$,
where $\ell=q^2\Tr_R L=q^{-1} a+q^{1} d,$ and $h=a-d$. The multiplication rules (\ref{sys1})
take the form
 \be
\begin{array}{l@{\hspace{20mm}}l}
q^{2}hb-bh +(q-\qq) \ell b =2_q\,\h\, b & b\ell - \ell b=0\\
q^{2}ch-hc +(q-\qq) \ell c =2_q\,\h\, c& c\ell - \ell c=0\\
2_q\, q\left(bc-cb\right)+\left(q^{2}-1\right)h^{2} +(q-\qq) \ell h=
2_q\,\h\, h & h \ell - \ell h = 0.
\end{array}
\label{sys2}
\ee

Observe that the element $\ell$ is central in this algebra. In what follows we also
consider the quotient $\slhq=\lhq/ \langle \ell \rangle$.
Thus, this algebra is generated by three elements $b,\,h,\,c$
subject to the system
\be
q^2hb-bh=2_q\,\h\, b,\quad q^2c h-hc=2_q\,\h\, c,\quad 2_q q (bc-
cb)+(q^2-1) h^2=2_q\,\h\, h.\label{sys3} \ee
Note that at $q=1$ this algebra $\slhq$ coincides with $U(sl(2)_{\h})$.

The explicit form of the representation (\ref{B-rep}) reads
$$
\rho_V(a)=\left(\begin{array}{cc}
\qq&0\\
0&0
\end{array}\right),\quad \rho_V(b)=\left(\begin{array}{cc}
0&q^{-3}\\
0&0
\end{array}\right),\quad \rho_V(c)=\left(\begin{array}{cc}
0&0\\
q^{-1}&0
\end{array}\right),\quad \rho_V(d)=\left(\begin{array}{cc}
0&0\\
0&q^{-3}
\end{array}\right).$$

In order to get the corresponding representation of the $sl$-reduced algebra (\ref{sys3}) we
use the general recipe suggested in \cite{S}. Let $\rho_U: {\cal L}(R_q,1)\to \End(U)$ be
an mREA representation and the central element $\Tr_R(L)$ is represented by a scalar operator
$$
\rho_U(\Tr_R(L)) = \chi_1\,{\rm Id}_U
$$
where $\chi_1 = \chi(\Tr_R(L))$ is the value of a character $\chi: Z(\loq)\rightarrow {\Bbb K}$ of the center
$Z(\loq)$.

Then the straightforward calculation shows that the ${\cal SL}(R_q,1)$ generators $F_i^{\,j}$
are represented in the space $\End(U)$ by the following linear operators
\be
\bar\rho_U(F_i^{\,j}) = \frac{1}{\omega}\,\Bigl(\rho_U(L_i^{\,j}) -
\delta_i^{\,j}\,\frac{\chi_1}{\Tr C}\,{\rm Id}_U\Bigr), \quad
\omega = 1- (q-q^{-1})\frac{\chi_1}{\Tr C}.
\label{sl-red-rep}
\ee

Applying formula (\ref{sl-red-rep}) to the above representation $\rho_V$ we get
the representation of the algebra $\sloq$ (\ref{sys3})
\be
\bar\rho_V(h)=w\,\left(\begin{array}{cc}
q&0\\
0&-q^{-1}
\end{array}\right),\quad
\bar\rho_V(b)=w\,\left(\begin{array}{cc}
0&q^{-1}\\
0&0
\end{array}\right),
\quad \bar\rho_V(c)=w\,\left(\begin{array}{cc}
0&0\\
q&0
\end{array}\right),
\quad w={q^2+1 \over q^4+1}. \label{taut} \ee

Let us go back to a general case and introduce another basic
(contragradient) representation $\rho_{V^*}: {\cal L}(R_q,1)
\rightarrow {\rm End}(V^*)$ of the algebra $\loq$ defined by the formula
\be
\rho_{V^*}(L_i^{\,j})\triangleright x^k = -x^rR_{ri}^{\,kj}.
\label{rep-dua}
\ee
Note that in the classical case (when $R$ is the usual flip)
we have  $\rho_{V^*}(a)=-\rho_V(a)^*$. Assuming $R$ to be involutive and
introducing an involution via $R$ as discussed in the next section we can get such
a contragradient representation. However, if $R$ is not involutive this method
fails.

In the case in question this representation is
$$
\rho_{V^*}(a)=\left(\begin{array}{cc}
-q&0\\
0&0 \end{array}\right),\,
\rho_{V^*}(b)=\left(\begin{array}{cc}
0&0\\
-1&0 \end{array}\right),\,
\rho_{V^*}(c)=\left(\begin{array}{cc}
0&-1\\
0&0 \end{array}\right),\,
\rho_{V^*}(d)=\left(\begin{array}{cc}
q^{-1}-q&0\\
0&-q \end{array}\right).
$$

As was mentioned above, in \cite{GPS3} there was constructed a rigid quasitensor
Schur-Weyl category $SW(V)$ of vector spaces generated by $V$ and $V^*$ such that each
object of this category can be endowed with
an $\loq$-module structure. The term ''rigid`` means that with any object $U$ this category
contains its dual. Recall that an object $U^*$ is dual to $U$ if there exists nontrivial pairings
$U\otimes U^*\to {\Bbb K}$ and $U^*\otimes U\to {\Bbb K}$ which are the categorical morphisms.
A peculiarity of the case related to an even $R$ is that (by assuming
$\det_q T$ to be central, see Introduction)
the dual space to $V$ can be identified with a subspace in some tensor power of the space
$V$ while for the general  case the dual space $V^*$ must be introduced independently.

Using the pairing (\ref{l-r-pairing}) we introduce an important morphism ${\rm tr}_R: {\cal L}\rightarrow
{\Bbb K}$ which will be referred to as a categorical trace
\be
{\rm tr}_R(L_i^{\, j}) = \delta_i^j.
\label{tr-r-cat}
\ee
\begin{remark}{\rm
The fact that the above definition of the ${\rm tr}_R$ gives a categorical
morphism follows from the identification of $L_i^{\,j}$ and $x_i\otimes x^j$. Now the categorical
trace just equals to the pairing
$$
{\rm tr}_R(L_i^{\, j})  = \langle x_i,x^j\rangle
$$
which is a morphism (see \cite{GPS3}).
 }
\end{remark}

Anyway, in a general case  the initial skew-invertible braiding $R$
can be extended in a unique way to the braidings $R_{U,W}$
for any two objects  $U$ and $W$  of the Schur-Weyl category.
Besides, all morphisms of the category are assumed to be natural
or functorial as defined in \cite{T}. This means that given two morphisms
 $f:U\rightarrow U'$ and $g:W\rightarrow W'$ the relation
$$
(g\ot f)\circ R_{U,W}=R_{U',W'}\circ (f\ot g),
$$
(where the symbol $\circ$ denotes the composition of the maps)
is assumed to be fulfilled. By putting $W'=W$ and $g=\Id$  we get a
condition on a morphism $f$ which means that $f$ is covariant.

As we said above the space $\LL$ is identified with  the product $V\ot V^*$.
Consequently, this space is an object of the Schur-Weyl category.

\begin{definition} Given an object $U$ of the Schur-Weyl category, a
representations $\rho_U: {\cal L}(R_q,1)\rightarrow \End(U)$ is called R-invariant
or equivariant if its restriction to the space $\LL$ is a categorical morphism.
\end{definition}

For example, the above representations $\rho_V$ and $\rho_{V^*}$ are equivariant.

Our next gaol is to exhibit a way of multiplying equivariant representations of the algebra $\loq$.
To this end we need the braiding
$$
R_{\rm End}: \End(V)^{\otimes 2}\to \End(V)^{\otimes 2},
$$
which is a particular but very important case of the aforementioned braidings $R_{U,W}$.

Its explicit form in the basis $\{L_1\stackrel{.}{\otimes} L_2\}$ of the space $\LL^{\ot 2}$
reads as follows \cite{GPS3}
\be
R_{\rm End}(L_1\stackrel{.}{\otimes} L_2)= \Tr_{(0)}(R_{10}L_1R^{-1}_{10}\stackrel{.}{\otimes}
L_1R_{10}\Psi_{02})P_{12},
\label{rend1}
\ee
where the operator $\Psi$ is defined in (\ref{def:psi}). (The notation $\stackrel{.}{\otimes}$ was introduced
in the previous section. For instance, the entries of the matrix $L_1\stackrel{.}{\otimes} L_2$ where
$L_1=L\ot \Id,\,\,L_2=\Id\ot L$ are $L_i^j\ot L_k^l$.)

In the tensor power $\LL^{\ot k}$ there is a basis different from that $\{L_{i_1}^{j_1}\ot L_{i_2}^{j_2}\ot...\ot L_{i_k}^{j_k}\}$.
This new basis is composed from the matrix elements of the tensor product
$L_{\overline 1}\stackrel{.}{\otimes}L_{\overline 2}\stackrel{.} {\otimes}\dots \stackrel{.}{\otimes}
L_{\overline k}$ and it turns out to be more convenient for our calculations. Here in accordance
with the general definition (\ref{copy}) we set
\be
L_{\overline 1} = L\otimes {\rm Id}^{\otimes (k-1)},\quad
L_{\overline {p+1}} =
R_pL_{\overline p}R_{p}^{-1},\quad p\le k-1.
\label{L-copies}
\ee
Using the skew-inverse operator $\Psi$ one can prove that the matrix
elements of $L_{\overline 1}\stackrel{.}{\otimes}L_{\overline 2}
\stackrel{.}{\otimes}\dots \stackrel{.}{\otimes} L_{\overline k}$ are in
one-to-one correspondence with those of
$L_{1}\stackrel{.}{\otimes}L_{2} \stackrel{.}{\otimes}\dots
\stackrel{.}{\otimes} L_{k}$ and thus form the basis in
$\LL^{\ot k}\cong ({\rm End}(V))^{\otimes k}$.

Making use of (\ref{rend1}) we find that
\be
R_{{\rm End}}(L_{\overline 1}
\stackrel{.}{\otimes} L_{\overline 2}) =
L_{\overline 2}\stackrel{.}{\otimes} L_{\overline 1},
\label{R-end}
\ee
and, as a consequence,
$$
R_{{\rm End}}(L_{\overline k}
\stackrel{.}{\otimes} L_{\overline p}) =
L_{\overline p}\stackrel{.}{\otimes} L_{\overline k},
\quad \forall p> k,\quad k,p\in {\Bbb Z}_+.
$$

Note that each homogeneous component of the associated graded algebra $\Gr\, \lhq$ can be identified with
an object of the category $SW(V)$. So, the braiding $\Ren$ can be extended onto the whole algebra
$\Gr\, \lhq$ and, consequently, to the algebra $\lhq$ (in particular, to $\loq$). Besides, the braidings transposing
elements of $\lhq$ and elements of arbitrary objects of the category $SW(V)$ are also defined (see \cite{GPS3}
for detail). Below the symbol ${\sf R}(a \ot b)$ stands for the general form of braiding transposing $a$ and
$b$ which are elements of objects of the category $SW(V)$ (in particular, one or both of them can belong to
$\lhq$).

Our next goal is to introduce a coproduct  which endows $\loq$ with a braided bi-algebra structure.
This enables us to multiply its equivariant representations. In order to describe
the coproduct and the corresponding braided bi-algebra structure we need the following definition.

\begin{definition}{\rm
Given any skew-invertible Hecke symmetry $R$,
introduce  a  braided associative
algebra $\mbox{\bf L$(R_q)$}$  by the data:
\begin{enumerate}
\item As a vector space over the field $\Bbb K$ the algebra
$\mbox{\bf L$(R_q)$}$ is isomorphic to the tensor product of two
copies of mREA ${\cal L}(R_q,1)$
$$
\mbox{\bf L$(R_q)$} \cong {\cal L}(R_q,1) \otimes {\cal L}(R_q,1)\,.
$$
\item The product $* : \mbox{\bf L$(R_q)$}^{\otimes
2}\rightarrow \mbox{\bf L$(R_q)$}$ is defined by the rule
\be
(a_1\otimes b_1)* (a_2\otimes b_2)=a_1 a^\prime_2 \otimes
b^\prime_1 b_2\,,\qquad
\forall a_i\otimes b_i \in \mbox{\bf L$(R_q)$}\,,
\label{br-pr}
\ee
where $a_1a^\prime_2$ and $b_1b^\prime_2$ are the
products of elements of mREA and
\be
a^\prime_2\otimes b^\prime_1= {\sf R}(b_1\otimes a_2)\,.
\label{ti-el}
\ee
\end{enumerate}
}
\end{definition}
Note that the associativity of the $*$-product is proved in \cite{GPS3}.

Let us define a coproduct  $\Delta: {\cal L}(R_q,1)\rightarrow \mbox{\bf L$(R_q)$}$
as a linear map of the following form:
\be
\begin{array}{l}
\Delta(e_{\cal L})= e_{\cal L}\otimes e_{\cal L}\\
\rule{0pt}{6mm} \Delta(L_i^{\,j}) = L_i^{\,j}\otimes e_{\cal L}+e_{\cal L}
\otimes L_i^{\,j} - (q-q^{-1})\sum_k L_i^{\,k}\otimes L_k^{\,j}\\
\rule{0pt}{6mm}
\Delta(ab)=\Delta(a)*\Delta(b)\qquad\forall\,a,b\in {\cal L}(R_q,1)\,.
\end{array}
\label{copr}
\ee
In addition to (\ref{copr}), we introduce a linear map
$\varepsilon :{\cal L}(R_q,1) \rightarrow {\Bbb K}$
\be
\varepsilon(e_{\cal L})= 1\qquad
\varepsilon(L_i^{\,j})= 0\qquad
\varepsilon(ab)= \varepsilon(a)\varepsilon(b)
\quad \forall\,a,b\in {\cal L}(R_q,1)\,.
\label{coed}
\ee

\begin{proposition} {\bf (\cite{GPS3}}
\label{pro:10} The maps $\Delta$ and $\varepsilon$ endow the algebra $\loq$ with
a bi-algebra structure. This means that the relations
$$
(\Id\otimes \Delta)\Delta = (\Delta\otimes \Id)\Delta ,\qquad
(\Id\otimes \varepsilon)\Delta = \Id =(\varepsilon\otimes \Id)\Delta
$$
hold and moreover $\De$ is a homomorphism of the algebra $\loq$ into
\mbox{\bf L$(R_q)$}.
\end{proposition}

\begin{remark}
{\rm
The above braided bialgebra structure was deduced from the braided bialgebra structure
in $\lq$  discovered by S. Majid (see \cite{M3, M4}). A passage from the latter structure to
the former one can be given by a shift of generators (\ref{shift-1}). Nevertheless, namely,
in the form (\ref{copr}) the coproduct involved in this construction is very useful for defining
"braided vector fields" (see section \ref{sec:8}). Also, note that for $q=1$ this coproduct takes
the well known additive form  on generators of a  Lie algebra or its generalized analog discussed
in the next section.}
\end{remark}

Now, we are able to define a product of two representations
of the algebra $\loq$.

\begin{proposition} {\bf \cite{GPS3}}\hspace{3pt}
Given two equivariant mREA modules $U$ and $W$, let
$\rho_U:{\cal L}(R_q,1)\rightarrow \End(U)$ and
$\rho_W:{\cal L}(R_q,1)\rightarrow \End(W)$ be
the corresponding equivariant representations.
Consider a map
$\rho_{U\otimes W}:\mbox{\bf L$(R_q)$}\rightarrow \End(U\otimes
W)$
defined by the following rule
{\rm
\be
\rho_{U\otimes W}(a\otimes b)\triangleright (u\otimes w) =
(\rho_U(a)\triangleright u^\prime)\otimes (\rho_W(b^\prime)
\triangleright w)\,, \qquad \forall\,a\otimes b\in \mbox{\bf L$(R_q)$},
\; \forall\, u\in U,\;\forall\,w\in W,
\label{rep-ll}
\ee }
where
$$
u^\prime \otimes b^\prime={\sf R}(b\otimes u)\,.
$$
Then the map {\rm (\ref{rep-ll})} defines a representation of the braided
algebra $\mbox{\bf L$(R_q)$}$ in the space $U\otimes W$.
\end{proposition}

This proposition implies the following corollary.

\begin{corollary}
Let $U$ and $W$ be two ${\cal L}(R_q,1)$-modules with
equivariant
representations $\rho_U$ and $\rho_W$. Then the map
${\cal L}(R_q,1)\rightarrow \End(U\otimes W)$
given by
 {\rm
\be a\mapsto \rho_{U\otimes W}(\Delta(a))\,,
\qquad \forall\,a\in {\cal L}(R_q,1)\,,
\label{rep-prod}
\ee}
where the coproduct $\Delta$ and the map $\rho_{U\otimes W}$
are given
respectively by formulae {\rm (\ref{copr})} and {\rm (\ref{rep-ll})},
is an equivariant representation of ${\cal L}(R_q,1)$.
\end{corollary}

Using this result we can endow any tensor products of $V$ and $V^*$
(as well as subspaces of these tensor products, extracted by Young projectors)
with a $\loq$-module structure so that the corresponding representation
is equivariant (see \cite{GPS3}).

As an important example we consider an analog of the "adjoint"
representation $\rho_{V\otimes V^*}$, arising from those $\rho_V$ and $\rho_{V^*}$.
To this end we need an explicit form of the braiding  $R_{{\rm End},V}:{\rm End}(V)
\otimes V \rightarrow V\otimes {\rm End}(V)$ since it comes in construction of the
product of representations.

The matrix of this braiding turns out to be \cite{GPS3}:
\be
R_{{\rm End},V}(L_1\stackrel{.}{\otimes} x_2) = x_2\stackrel{.}{\otimes}
{\rm Tr}_{(0)}(R_{10}L_0\Psi_{02})P_{12}.
\label{br-m}
\ee
The action
$$
\rho_{V\otimes V^*}(L_{i_1}^{\,j_1})
\triangleright L_{i_2}^{\,j_2} = (\rho_V\otimes \rho_{V^*})\circ
\Bigl(\Delta(L_{i_1}^{\,j_1})\otimes( x_{i_2}\otimes x^{\,j_2})\Bigr)
$$
or, in matrix form, $\rho_{V\otimes V^*}(L_1)\triangleright L_2$
can be easily calculated with the use of formulae (\ref{B-rep}),
(\ref{copr}),
(\ref{rep-dua}), and (\ref{br-m}).  The result is given by a
somewhat complicated expression:
$$
\rho_{V\otimes V^*}(L_1)\triangleright L_2 =
(q-q^{-1})L_1I_2+L_1B_2P_{12}-{\rm Tr}_{(0)}(R_{10}L_1R_{10}
\Psi_{02})P_{12}.
$$

Nevertheless, in the basis of quantum matrix copies (\ref{L-copies}) the
adjoint action has a nice and simple form, namely
\be
\rho_{V\otimes V^*}(L_{\bar 1})\triangleright L_{\bar 2} =
L_1R_{12} - R_{12}L_1.
\label{adj-rep}
\ee
This example shows again that basis of copies (\ref{L-copies}) is very useful tool.

Now, we restrict ourselves to the standard case and explain a way of constricting
a family of representations of the corresponding (m)REA arising from methods
\cite{FRT} and \cite{LS}. Note that the approach  of \cite{LS}  was used in \cite{DS}
for constructing representation theory of the braided Lie algebra $sl(2)_q$ (see section \ref{sec:5}).

The authors of \cite{FRT} suggested a way of finding the center of the QG $\Uqq$.
This center can be easily calculated via the formula $\Tr_R L^k$
where $L$ is a matrix with entries belonging to the QG. Lately it was understood that
this matrix is subject to the defining relations of the corresponding REA.
The matrix $L$ was defined in \cite{FRT} via the formula
\be
L=S(L^-)L^+,
\label{REA-emb}
\ee
where  $L^\pm$ are the matrices
composed of the quantum group generators and $S(\cdot)$ is the antipode.

Let us consider the simplest example related to the QG $\uq$. In this case
with the use of  the basis $\{\ell,h,b,c\}$ discussed at the beginning of this section we
get the following parametrization of the mREA generators via the generators
$E$, $F$ and $H$ of the QG $\uq$
\be
h = -\frac{\lambda}{q}(qEF-q^{-1}FE),\quad
b = \lambda q^{-H}E,\quad c =\lambda q^{-H-2}F,\quad
\ell  = q^{2H-1}+q^{-2H-3}+ \frac{\lambda^2}{q^2}FE,
\label{ex:emb}
\ee
where we set for shortness $\lambda = q-q^{-1}$. The generators $E$, $F$ and $H$ are
normalized as follows
$$
q^HE = qEq^{H},\qquad q^HF = q^{-1}Fq^{H},\qquad [E,F] = \frac{q^{2H}-q^{-2H}}{q-q^{-1}}.
$$

However, the generators $\ell,h,b,c$ are not independent in the above parametrization.
In \cite{FRT} the problem of finding relations between these generators
was not considered. Such relations in the $\Uq$ case were found in \cite{LS}.
In that paper for the QG $\Uq$ a family of generators similar to those above was constructed.
In the simplest case ($n=2$) these generators $X_0$, $X_{\pm}$ and $C$ satisfy the system
of relations \cite{LS}
\begin{eqnarray}
 && q^2X_0 X_+-X_+ X_0=q C X_+ \nonumber\\
&& q^{-2}X_0 X_--X_- X_0=-q^{-1} C X_- \nonumber\\
&&X_+X_- - X_-X_+ +(q^2 - q^{-2})X_0^2 = (q+q^{-1})CX_0\nonumber\\
&& CX_m= X_mC,\quad (m=0,\pm),\label{sl-LS}
\end{eqnarray}
and
\be
C^2-(q-q^{-1})^2\Bigl(X_0^2 + \frac{qX_-X_+ + q^{-1}X_+X_-}{q+q^{-1}}\Bigr) = 1.
\label{qv-dep}
\ee
It is easy to see that upon setting
$$
X_0 = \frac{qC}{2_q}\,h,\qquad X_+ = qCb, \qquad X_- = qCc
$$
the first three lines of (\ref{sl-LS}) transform precisely into relations (\ref{sys3}) with $\h=1$.
Thus, the algebra ${\cal SL}(R_q,1)$ is embedded in the QG $\uq$ extended by the element $C^{-1}$.

Since $C$ is central, it becomes a nontrivial scalar in each finite dimensional representation of the
QG $\uq$. This enables one to get a representation of the algebra $\slhq$ once a $U_q(sl(2))$ representation
is given. Note that such a representation of the algebra $\slhq$ can be treated as that of
$\lhq$ but with the image of $\ell$ depending on the images of $b$, $c$ and $h$ (due to (\ref{qv-dep})).
Thus, we can get a subfamily of all finite dimensional representations of the corresponding REA constructed
by the general method. The same is valid for all QG $\Uq$.

In conclusion we want to emphasize that our method of constructing the REA representations does not give
rise to non-equivariant ones. For example, the standard REA has one dimensional representations, which cannot be
obtained either by our method or by that of \cite{LS}. We refer the reader to \cite{K, Mu1} for such type
representations. Also note that
our approach is valid for a generic $q$ whereas the method of \cite{LS} allows to consider
the case of special $q$ (roots of unity). This case was considered in
\cite{DS}. Besides, this method enables one to construct Verma type modules for the
standard algebra $\slhq$ (also, see \cite{LS}).
Note that in general we cannot define such type modules.

However, by comparing our approach and that from \cite{LS} we would like to stress again that
the latter one cannot be applied to the REA connected with nonstandard braidings since the corresponding
QG like objects are not known for the general type braidings.

\section{Braided Lie algebras}
\label{sec:5}

In this section we discuss the role of the mREA in the definition of
the quantum (braided) analogs
of a Lie algebra. We are mainly interested in the Lie algebra type
objects similar to $gl(n)$ and $sl(n)$.
For such objects the problem can be formulated as follows. Let
$R:\vv\to\vv$ be a skew-invertible
braiding. Consider the space $\End(V)$ equipped with the braiding
$\Ren:\End(V)^{\ot 2}\to \End(V)^{\ot 2}$ (\ref{R-end}).
We want to introduce a braided analog of the Lie bracket in the
space $\End(V)$, which is an
$\Ren$-invariant operator
$$
[\,\,,\,]:\quad\End(V)^{\ot 2}\to \End(V)
$$
and to define the corresponding {\em enveloping algebra}
$U(\End(V))$ with the good deformation property.
Actually, we first introduce an analog of the enveloping algebra
$U(gl(n))$ and then define a $gl(n)$ type
bracket.

Before going into detail we would like to make some historical
comments. The first generalization of a
Lie algebra was a super-Lie algebra. We were informed by M.
Gerstenhaber and J. Stasheff, to whom
we express our profound gratitude, that a super-version of the
Jacobi identity was introduced in \cite{N}
(also, in the 40-50's  super-algebra type objects
were popular in connection with
the Whitehead product). Lately, the super-algebras were
intensively studied by physicists in the frameworks
of models with fermion-boson symmetries.

The basic example is a super-Lie algebra $gl(m|n)$ defined in the space
$\End(V)$ where $V=V_{(m|n)}=V_{\overline{0}}\oplus V_{\overline{1}}$ is
a $\Z_2$ graded space and
$\dim V_{\overline{0}}=m$, $\dim V_{\overline{1}}=n$. A braiding
$R$ is a super-flip
$R(x\ot y)= (-1)^{{\overline x}{\overline y}}y\ot x$.

Attempts to introduce the Lie algebra type objects graded by other
finite commutative group
were undertaken in the late 70's \cite{Sch}.

The next generalization of the Lie algebra notion was related to
the involutive symmetries $R:\;R^2={\rm Id}$.
In \cite{G1} the following notion was introduced.
\begin{definition}
\label{def:gLie}
{\rm
The data
$$
(V,\, R:\vv\to\vv,\, [\,\,,\,]:\vv \to V)
$$
where $R$ is an involutive symmetry is called a {\em generalized
Lie algebra} if the axioms below hold true:
\begin{enumerate}
\item $[\,\,,\,]\,R(X\ot Y)=-[X,Y]$;
\item
$[\,\,,\,]\,[\,\,,\,]_{12}(I+R_{12}R_{23}+R_{23}R_{12}) (X\ot Y\ot Z)=
0$;
\item $R [\,\,,\,]_{23} (X\ot Y \ot Z)= [\,\,,\,]_{12}
R_{23}R_{12}(X\ot Y\ot Z)$
\end{enumerate}
for all $X,Y,Z\in V$.
}
\end{definition}
Note that the generalized Jacobi identity (the axiom 2) can be also
presented in one of the following
equivalent forms
\begin{itemize}
\item[2.a]
$[\,\,,\,]\,[\,\,,\,]_{23}(I+R_{12}R_{23}+R_{23}R_{12})(X\ot Y\ot Z)=
0$;
\item[2.b]
$[\,\,,\,]\,[\,\,,\,]_{12}(X\otimes(Y\ot Z-R(Y\ot Z)))=[X,[Y,Z]]$;
\item[2.c]
$[\,\,,\,]\,[\,\,,\,]_{23}((X\ot Y-R(X\ot Y))\otimes Z)=[[X,Y],Z]$.
\end{itemize}

Let us denote the generalized Lie algebra originated from a vector
space $V$ and a braiding $R$
as $\gg(V,R)$. Its enveloping algebra $U(\gg(V,R))$ can be naturally
defined as a quotient of the free tensor algebra $T(\gg(V,R))$:
$$
U(\gg(V,R))=T(\gg(V,R))/\langle X\ot Y-R(X\ot Y)-[X,Y]
\rangle\,\quad \forall\,X,Y\in V.
$$
The enveloping algebra $U(\gg(V,R))$ has the good deformation
property. In particular, due to the PBW-like
theorem \cite{PP} the associated graded algebra ${\rm Gr}\,U(\gg(V,R))$ is canonically
isomorphic to the algebra
$$
\Sym(\gg(V,R))=T(\gg(V,R))/\langle{\rm Im}({\rm Id}-R)\rangle.
$$

Besides, it becomes a {\em braided Hopf algebra}, being equipped
with an appropriate coproduct $\De$
and an antipode $S$. These operators have the classical form
$\De(X)=X\ot 1+1\ot X$ and $S(X)=-X$
on elements $X\in V$ and can be naturally extended to the
whole algebra with the use of the operator $R$ (see \cite{G2}).

The axioms of Definition \ref{def:gLie} are satisfied by the
following data
\be
(\End(V),\, \Ren:\End(V)^{\ot 2}\to \End(V)^{\ot 2},\, [\,\,,\,]: \End(V)
^{\ot 2}\to \End(V))
\label{glVR}
\ee
where $\Ren$ is an extension of a skew-invertible involutive
symmetry $R:\vv\to\vv$ to the
space $\End(V)^{\otimes 2}$ (see section \ref{sec:4}) and the
bracket is defined by
\be
[\,\,,\,]=\circ(\Id-\Ren).
\label{skobka}
\ee
Here $\circ:\End(V)^{\ot 2}\to \End(V)$ is the usual product in the
algebra $\End(V)$. For the
generalized Lie algebra (\ref{glVR}) we use the notation
$gl(V,R)$.

Now, we suppose $R=R_q$ to be a Hecke symmetry. Consider
the corresponding mREA $\loq$. As was shown in \cite{GPS3},
the algebra $\loq$ possesses the good deformation property.

Now we define a {\em braided Lie algebra} such that its enveloping
algebra coincides with the
mREA $\loq$. To this end we rewrite the mREA multiplication rules
(\ref{mREA0}) in the following
form
\be
L_{\overline 1}\,L_{\overline 2} - R^{-1}_{12}\, L_{\overline 1}\,
L_{\overline 2}\, R_{12} = L_1\,R_{12} - R_{12}\,L_1\,,
\label{mREA1}
\ee
where the matrices $L_{\overline k}$ are defined in (\ref{L-copies}).
Then, on the linear space $\L=\span(L_i^{\,j})\cong \End(V)$ we
introduce a Lie type
bracket $[\,\,,\,]:\L^{\ot 2}\to \L$ setting by definition
\be
 [L_{\overline 1},L_{\overline 2}]= L_1\,R_{12} - R_{12}\,L_1.
\label{REA-br}
\ee
This bracket can be written in a form, similar to (\ref{skobka}). For
this purpose we define
the linear operator $Q:\;{\cal L}^{\otimes 2}\to {\cal L}^{\otimes 2}$:
\be
Q(L_{\overline 1}\stackrel{.}{\otimes}L_{\overline 2}) = R^{-1}_{12}\,
L_{\overline 1}
\stackrel{.}{\otimes}L_{\overline 2}\, R_{12}.
\label{Q-op}
\ee

Then, as can be easily seen from (\ref{mREA1}), the bracket
(\ref{REA-br}) reads
\be [\,\,,\,]=\circ (\Id-Q) \label{bra} \ee
where $\circ$ has the same meaning as above. Note that the
multiplication of the operators
$L_i^{\,j}$ can be calculated on the base of (\ref{B-rep})
\be
L_{i_1}^{\,j_1}\circ L_{i_2}^{\,j_2} = B_{i_2}^{\,j_1}L_{i_1}^{\,j_2}\,.
\label{prod}
\ee
One of  important distinction  between  the bracket (\ref{REA-br}) and that
(\ref{skobka}) of the generalized
Lie algebra $gl(V,R)$ consists in the fact that the operator $Q$
differs from the operator
$\Ren$ whereas if $R$ is an involutive symmetry these operators
are equal to each other.

In terms of the operator $Q$ we can express an operator
${\cal S}_q: {\cal L}^{\otimes 2}\to {\cal L}^{\otimes 2}$
possessing the property
$$
{\cal L}^{(2)}(R_q)=\Im\, {\cal S}_q,
$$
where ${\cal L}^{(2)}(R_q)\subset {\cal L}^{\otimes 2}$ is the second
order homogeneous component of the REA ${\cal L}(R_q)$.
Such an operator ${\cal S}_q$ can be defined by the following formula (see \cite{GPS3})
$$
{\cal S}_q = \frac{1}{2_q^2}\,((q^2+q^{-2})\,{\rm Id} +Q+Q^{-1}).
$$

Thus, the operator ${\cal S}_q$ can be treated as a total
$q$-symmetrizer on the space ${\cal L}^{\otimes 2}$.
\begin{theorem}
The bracket {\rm (\ref{REA-br})} possesses the following
properties:
\begin{itemize}
\item[{\rm 1.}] The bracket is skew-symmetric in the following
sense
{\rm
\be
[\,\,,\,]\,{\cal S}_q(L_{\overline 1}\stackrel{.}{\otimes}
L_{\overline 2}) = 0\,,
\label{q-skew}
\ee}
where ${\cal S}_q$ is the symmetrizer introduced above.
\item[{\rm 2.}] The bracket satisfies the generalized Jacobi identity
of the form
$$
[\,\,,\,]\,[\,\,,\,]_{23}((X\ot Y-Q(X\ot Y))\otimes Z)=[[X,Y],Z]\qquad
\forall\, X,Y,Z\in{\cal L},
$$
otherwise stated, the adjoint action defined by
{\rm
\be
\ad\,L_i^{\,j}(L_k^{\,l})
=[L_i^{\,j},L_k^{\,l}]
\label{ad-rep}
\ee
}
is a representation of the algebra $\loq$.
\item[{\rm 3.}] The bracket is $\Ren$-invariant. This means that
{\rm
\be
\begin{array}{l}
R_{\End}[\,\,,\,]_{23} =
[\,\,,\,]_{12}(R_{\End})_{23}(R_{\End})_{12},\\
\rule{0pt}{7mm}
R_{\End}[\,\,,\,]_{12} =
[\,\,,\,]_{23}(R_{\End})_{12}(R_{\End})_{23}\,,
\end{array}
\label{R-inv}
\ee
}
where the both sides of the above equalities are operators acting
in ${\cal L}^{\otimes 3}$.
\end{itemize}
\end{theorem}

{\bf Proof.} To prove the first statement of the theorem, we note
that the operator $Q$
has the minimal polynomial of the form \cite{GPS3}
$$
(Q+q^2{\rm Id})(Q+q^{-2}{\rm Id})(Q - {\rm Id})=0
$$
which in turn gives us the following expression for the inverse operator
$Q^{-1}$
$$
Q^{-1} = Q^2 + (q^2+q^{-2}-1)\,Q - (q^2+q^{-2}-1)\,{\rm Id}.
$$
Now, the property (\ref{q-skew}) follows directly from the relation
(\ref{bra}) and explicit forms of
${\cal S}_q$ and $Q^{-1}$ written above.

The second claim was actually proved in section \ref{sec:4},
relation (\ref{adj-rep}), where the
construction (\ref{ad-rep}) was found from the general approach to
the mREA representation theory.
But it can be also easily verified by direct calculations. To this end we rewrite the
''adjoint action'' (\ref{ad-rep}) in the
basis of matrix copies $L_{\bar k}$. Then this action is given by
relation (\ref{REA-br})
$$
L_{\bar 1}\triangleright L_{\bar 2} = L_1R_{12} - R_{12}L_1 =
(L_{\bar 1} - L_{\bar 2})R_{12},
$$
where the symbol $\triangleright$ stands for the adjoint action.

This entails that
\be
L_{\bar 2}\triangleright L_{\bar 3} =  (L_{\bar 2} - L_{\bar 3})R_{23}.
\label{REA-br2}
\ee
Now, by  writing the mREA commutation relation as
$$
R_{12}L_{\bar 1}L_{\bar 2} - L_{\bar 1}L_{\bar 2}R_{12} = L_{\bar
2} - L_{\bar 1}
$$
we apply  the both sides of this equality to the element $L_{\bar 3}$
and verify that the results are equal to each other. Indeed, we have
$$
R_{12}L_{\bar 1}\triangleright(L_{\bar 2}\triangleright L_{\bar 3}) =
(L_{\bar 2} - L_{\bar 1})R_{23} - R_{12}R_{23}R_{12}^{-1}(L_{\bar
2} - L_{\bar 1}),
$$
$$
L_{\bar 1}\triangleright(L_{\bar 2}R_{12}\triangleright L_{\bar 3}) =
L_1R_{12}^{-1}R_{23}^{-1}R_{12} - R_{12}^{-1}L_1R_{23}^{-1}
R_{12} -
L_1R_{23}^{-1} + R_{23}^{-1}R_{12}^{-1}L_1R_{12}.
$$
and
$$
(L_{\bar 2} - L_{\bar 1})\triangleright L_{\bar 3} = L_{\bar 2}R_{23} -
R_{23}L_{\bar 2}
- R_{23} L_1R_{12}^{-1}R_{23}^{-1} + R_{23}R_{12}^{-1}
L_1R_{23}^{-1}.
$$
Then, taking the difference of the first two expressions and
applying the Hecke
condition to connect $R$ and $R^{-1}$, we convince ourselves that
the difference
coincide with the third expression above, that is
$$
R_{12}(L_{\bar 1}\triangleright)( L_{\bar 2}\triangleright)  -
(L_{\bar 1}\triangleright )(L_{\bar 2}\triangleright) R_{12} =
(L_{\bar 2}\triangleright) - (L_{\bar 1}\triangleright).
$$

The proof of the third claim of the theorem is a matter of a trivial
calculation on the base
of (\ref{R-end}), (\ref{REA-br}) and (\ref{REA-br2}). Consider, for
example, the action of
the first relation (\ref{R-inv}) on the basis element
$L_{\bar 1}\stackrel{.}{\otimes}L_{\bar 2}\stackrel{.}{\otimes}
L_{\bar 3}$ (we shall omit the symbols $\stackrel{.}{\otimes}$ for simplicity).
The right hand side
of the relation  leads to the following chain of transformations:
$$
L_{\bar 1}L_{\bar 2}L_{\bar 3}\stackrel{(\Ren)_{12}}{\longrightarrow}
L_{\bar 2}L_{\bar 1}L_{\bar 3}\stackrel{(\Ren)_{23}}{\longrightarrow}
L_{\bar 2}L_{\bar 3}L_{\bar 1}\stackrel{[\,,\,]_{12}}{\longrightarrow}
(L_{\bar 2} - L_{\bar 3})R_{23}L_{\bar1},
$$
while the left hand side gives
$$
L_{\bar 1}L_{\bar 2}L_{\bar 3}\stackrel{[\,,\,]_{23}}{\longrightarrow}
L_{\bar 1}(L_{\bar 2} - L_{\bar 3})R_{23}\stackrel{\Ren}
{\longrightarrow}
(L_{\bar 2} - L_{\bar 3}) L_{\bar 1}R_{23}
$$
which coincides with the above result for the right hand side since
$L_{\bar 1}$
commutes with the matrix $R_{23}$.\hfill\rule{6.5pt}{6.5pt}
\medskip

We keep the notation $gl(V,R)$ for the following data
$$(\LL=\span(L_i^j),\, Q:\LL^{\ot 2}\to \LL^{\ot 2},\, [\,,\,]=\circ(\Id-Q))$$
which is  similar to that  above related to an involutive $R$.
We call this data the {\em braided Lie algebra} of {\it gl}-type.

Now, consider its $sl$-reduction. For this
purpose we pass to the mREA generators $F_i^{\,j}$ and $\ell$
introduced in
(\ref{shift}). Recall that this passage requires $\Tr(C)\not=0$.
The commutation relations of the new mREA generators are given
by (\ref{F-l}).
It is not difficult to calculate the adjoint action in terms of $F$ and
$\ell$:
\begin{eqnarray}
&&\ad\,\ell(\ell) = 0,\qquad
\ad\,\ell(F_1) =
-(q-q^{-1})\,\Tr(C)\,F_1
\nonumber\\
&&\ad\,F_{\overline 1}(\ell) = 0\,,
\nonumber\\
&&\ad\,F_{\overline 1}(F_{\overline2})
 = F_1 R_{12} - R_{12}F_1 + (q-q^{-1}) R_{12}F_1 R_{12}^{-1}\,.
\label{sl-ad}
\end{eqnarray}

As was defined in section \ref{sec:3}, the $sl$-reduction of the mREA
${\cal L}(R_q,1)$
consists in passing to the quotient algebra ${\cal SL}(R_q,1)$
(\ref{SL-red}). The explicit
commutation relations of the ${\cal SL}(R_q,1)$ algebra read
\be
R_{12}F_1 R_{12}F_1 - F_1 R_{12}F_1 R_{12} =
(R_{12}F_1 - F_1 R_{12}),\qquad \Tr_R(F) = 0.
\label{sl-comm}
\ee
Let us consider the space
$${\cal SL} = \span\{F_i^{\,j}\}\subset {\cal L} = \span\{L_i^{\,j}\}$$
formed by the traceless elements with respect to the categorical
trace ${\rm tr}_R$  (\ref{tr-r-cat}):
$$
\SL = \{X\in {\cal L}\,|\,{\rm tr}_R(X) = 0\}.
$$
The restriction $Q\mapsto Q_{sl}:{\cal SL}^{\ot 2}\to {\cal SL}^{\ot 2}
$ can also be naturally defined. It is easy to see that the categorical
trace of the bracket (\ref{REA-br}) is zero, so the bracket can
be restricted to the subspace ${\cal SL}$. However, the
corresponding ``adjoint action'' of the restricted
bracket following from (\ref{sl-comm})
$$
\ad_{sl}\,F_{\overline 1}(F_{\overline 2})=[F_{\bar
1},F_{\bar2 }]_{sl} =
F_1 R_{12} - R_{12}F_1
$$
does not define a representation of $\SL$ because it does not
contain the last term of formula (\ref{sl-ad}).

Thus, if similarly to the braided Lie algebra $gl(V,R)$ we define
the algebra $sl(V,R)$ as the following data
$$
(\SL=\span(F_i^{\,j}), Q_{sl}:\SL^{\ot 2}\to \SL^{\ot 2} , [\,\,,\,\,]
_{sl}:\SL^{\ot 2}\to \SL)
$$
where $[\,\,,\,]_{sl}$ is the bracket (\ref{bra}) restricted to $\SL$,
we can see that  the Jacobi identity in the form
valid for the former algebra fails for
the algebra $sl(V,R)$.

A particular case of such a braided Lie algebra related to
$\Uq$ was introduced in \cite{LS} (see the previous section).
 A similar construction was suggested in \cite{DGHZ}.

In our approach we need no object of QG type. Moreover, our
construction is valid for Hecke symmetries of general type. In
particular, it embraces $q$-deformations
of Lie super-algebras $gl(m|n)$.

Quantum Lie algebras related to the QG $\Uqq$ of other series were
introduced in \cite{DGG}. Let us explain the main idea of that construction.
Fix a Lie algebra $\gg$ and consider a $\Uqq$-covariant analog $[\,\,,\,]_q$
of the classical Lie bracket $[\,\,,\,]:\gg^{\ot 2}\to\gg$. For this end we decompose
the $\Uqq$-module $\gg^{\ot 2}$ into a direct sum of irreducible $\Uqq$-modules.
Here we assume the space $\gg$ to be endowed  with a $\Uqq$ action deforming
the usual adjoint action. By means of the QG coproduct this action can be extended
onto $\gg^{\ot 2}$. In contrast with the previous case the
$\Uqq$-module $\gg^{\ot 2}$
is multiplicity free. So, the bracket $[\,\,,\,]_q$ can be defined in a
unique (up to a factor) way
similarly to their classical counterparts but in the category of
$\Uqq$-modules. Thus, a "quantum (braided) Lie algebra" bracket can
be introduced.

However, we know no reasonable Jacobi identity which could be
written for such a "quantum Lie algebra". Moreover, the
above "enveloping algebra" does not possesses the good
deformation property. The point is that even the "symmetric
algebra" corresponding to such an "enveloping algebra" is not a
deformation of its classical counterpart.

As was mentioned in \cite{GPS3}, we think that an axiomatic introduction of
a generalized (quantum, braided) Lie algebra
is possible iff the corresponding
braiding is involutive. Nevertheless, in some papers such algebras
related to non-involutive braidings
are introduced by the same (or very close) axiom system.
However, this way of introducing the Lie algebra type objects is not justified
by exhibiting meaningful examples.

To conclude this section, we would like to discuss the problem of
defining a conjugation (involution) in a generalized (quantum,
braided) Lie algebra.
Let us first assume that $R:\vv\to\vv$ is a skew-invertible real
involutive symmetry (so, $\K=\R$). Also, suppose that there exists a
nondegenerated pairing $\vv\to \R$ which is $R$-invariant.

Then the space $\End(V)$ can be identified with $\vv$ and this
identification is a categorical morphism in terminology
 of section \ref{sec:4}. Let us introduce an involution in $\End(V)$ as the
 image of the operator $R:\vv\to\vv$.
By passing to the complexification of the algebra $\End(V)$
we complete this operator with the complex conjugation of
numerical coefficients. The final operator
$$
*:\End(V)\to \End(V),\quad x\to x^*
$$
is called {\em conjugation}. It can be naturally extended to the
enveloping algebra $U(gl(V,R))$ via the relation
$(x\circ y)^*=\circ (*\ot *)R(x\ot y)$ where $\circ$ is the product in
$U(gl(V,R))$.

This conjugation is involutive, $\Ren$-invariant and it is coordinated
with the bracket of generalized Lie algebra $gl(V,R)$ in the
following sense
\be
[x,y]^*=-[x^*, y^*]. \label{invol} \ee
Note that this relation is universal: it does not depend on
 a skew-invertible involutive symmetry $R$. The latter relation entails
that the family of elements which are skew-symmetric with respect to the conjugation $x^*=-x$
is a subalgebra. Also, the map $x\to -x^*$ is an isomorphism of the generalized Lie algebra
$gl(V,R)$.

Note that if $R$ is the usual twist and the pairing is Euclidean the subalgebra of skew-symmetric
elements is just $u(n)$.

Now, let $V$ be a vector space endowed with a Hecke symmetry. We are interested
in a problem of classification of all involutions (conjugations) in the algebra $\End(V)$
which are $R$-invariant and verifying (\ref{invol}). In the $U_q(sl(2))$ case it can be shown by
 a direct calculation
that the only $R$-invariant linear operators on the space $\End(V)$ are scalar on each
irreducible component in the decomposition $\End(V) = \K\oplus {\cal SL}$ where ${\cal SL}$
is the subspace of the traceless elements.
So, the only   involution  on ${\cal SL}$ which is R-invariant operator
reads $x\to x^* = -x$. It is out of
interest and is not a deformation of the involution on the algebra $sl(2,{\Bbb C})$ giving
rise to the algebra $su(2)$. We would like to emphasize that using in algebras in question involutions
which are not categorical morphisms is not motivated by their "braided nature" and leads to some
shortcoming (see the next section).

\section{Quantum sphere: different approaches}
\label{sec:6}

In this section we consider different approaches to introducing
quantum sphere (hyperboloid) and developing some aspects of
geometry on it.
Historically, the RTT algebra was the first defined quantum
matrix algebra. So, quantum homogeneous space algebras (and in the first
turn, the quantum sphere) were initially defined by imitating the classical definition of
homogeneous $G$-spaces as cosets $G/H$.

However, as explained in section \ref{sec:2} on a symmetric orbit there exists another way
of defining the corresponding "quantum variety" via quotienting the standard (m)REA.
Let us compare these ways on the example of a quantum sphere.
Though in this case the both ways lead to equivalent results, they
yield completely different methods of adjacent geometry, in particular,
those of constructing projective modules over $q$-spheres.
Now, describe a quantum sphere explicitly.

First, describe the quantum algebra $\K_q[SL(2)]$ which is a
quantum counterpart
 of the Sklyanin bracket on the group $SL(2)$.
Computing the defining relations for the entries of the
matrix $T=\left(\begin{array}{cc}
a&b\\
c&d
\end{array}\right)$ we get
$$ab=qba,\,\, ac=qca,\, \,bd=qdb,\, \,bc=cb,\,\, cd=qdc,$$
$$ ad-da=(q-\qq)bc,\,\, ad-qbc=da-\qq bc=1.$$
Consider the subalgebra of this algebra consisting of elements
invariant with respect to the coaction
$$\De: \left(\begin{array}{cc}
a&b\\
c&d
\end{array}\right)\to \left(\begin{array}{cc}
a\ot z&b\ot z^{-1}\\
c\ot z&d\ot z^{-1}
\end{array}\right)$$
where $z$ is a formal invertible indeterminate. We shall denote this subalgebra
$\K_q[SL(2)/H]$. It is a $q$-deformation of the algebra of functions $\K[SL(2)/H]$.

The algebra $\K_q[SL(2)]$ can be equipped with the following conjugation operator
\be
a^*=d, \quad b^*=-qc,\quad c^*=-\qq b, \quad d^*=a,
\label{conj-op}
\ee
which is assumed to be antilinear and subject to the usual condition
$$(x\,y)^*=y^*\, x^*,\quad \forall x,y\in \K_q[SL(2)].$$
Consequently, it is involutive.
Here we assume $\K=\C$ and $q$ to be real.

Note that the Hecke symmetry coming in the definition of the
algebra $\K_q[SL(2)]$ is a particular case of the so-called
braiding of real type as defined in \cite{M4} (Definition 4.2.15). In
our normalization of a Hecke symmetry $R$ this condition reads
${\overline R}^{{\sf T}}=R$. If $R$ is such a Hecke symmetry
it is possible to introduce an involution in the corresponding RTT
algebra in a similar way.

It is not difficult to see that the subalgebra $\K_q[SL(2)/H]$ is
closed with respect to the conjugation (\ref{conj-op}). The algebra
$\K_q[SL(2)/H]$ equipped with this
involution will be denoted $\K_q[SU(2)/H]$. Namely, this algebra is
usually considered as one of avatars of a quantum sphere.

Being equipped with the above conjugation (\ref{conj-op}), the algebra
$\K_q[SL(2)]$ is treated to be a quantum counterpart of the algebra
$\K[SU(2)]$. We denote it as $\K_q[SU(2)]$. All its $*$-representations
(i.e. those respecting the $*$-operator) in a Hilbert space were classified
in \cite{VS}, where it was shown that there is a series of one-dimensional
representations and an infinite dimensional one. They can be restricted to
the subalgebra $\K_q[SL(2)/H]$.

Now, we consider a way of defining the quantum sphere (hyperboloid)
via the REA. To this end consider the mREA $\lhq$ related to the
standard Hecke symmetry. This algebra is explicitly given by the system (\ref{sys1})
or (\ref{sys2}) in the generators $\ell, h, b, c$. Also, consider its quotient
 $\slhq$ defined by the system (\ref{sys3}). We shall refer to the algebra $\slhq$
as {\em $q$-noncommu\-tative} if $\h\not=0$ and as {\em $q$-commutative} if $\h=0$.
Note, that since the algebras $\lhq$ and $\lq$ are isomorphic for $q\not=\pm 1$
we avoid to call them in a similar manner by keeping this terminology for their $sl$-quotients
which are not isomorphic to each other for any value of $q$.

Consider the central element
$$
\Cas=q^2\,\Tr_q L^2=\qq a^2+\qq bc+q cb+q d^2={\ell^2\over
2_q}+ \qq\, bc+{h^2 \over 2_q}+q\, cb\in \lhq.
$$
Its image in the algebra $\slhq$ reads
$$ \Cas_{sl}=\qq\, bc+{h^2 \over 2_q}+q\, cb.$$
It is also central in this algebra. So, it is natural to introduce the
following quotient
$$
\slhqc=\slhq/\langle \Cas_{sl}-C \rangle
$$
where $C\in \K$ is a number. We assume that $C\not=0$.

In what follows the elements $\Cas$ and $\Cas_{sl}$ are called
the {\em $q$-Casimir elements} in the algebras $\lhq$ and $\slhq$
respectively. Note that the center of the algebra $\lhq$ is generated
by the elements $\Cas_{sl}$ and $\ell$.

The algebra $\slhqc$ is called the {\em quantum hyperboloid} since at
$q=1$ we get just a usual hyperboloid (one- or two-sheeted in
dependence on $C\in \K = \R$). By a {\it quantum sphere} one often
means this algebra but considered over the field of complex numbers ($\K=\C$)
and endowed with a conjugation (involution) defined on the generators as follows
$$
\ell^*=\ell,\qquad h^* = h,\qquad  b^*=c,\qquad c^*=b.
$$
Onto the whole algebra this involution is extended via the classical properties.

Nevertheless, such an involution does not allow us to define a
quantum sphere as a real algebra. This is a shortcoming of the involution
which is not a categorical morphism whereas a quantum hyperboloid can be treated
as a real algebra if $q\in \R$.

Considering irreducible representations $V_k,\,\,\dim V_k=k+1$
of the algebra $\sloq$ we can see that the constant
$C$ coming in the definition of the algebra
$\sloqc$ depends on $k$. The situation is similar to the case of
the algebra $U(sl(2))$ in which the value of the Casimir element
depends on $k$ as well. By direct computation we get \cite{S}
$$
C(k)=C(V_k)=q^{-2}\,\frac{k_q(k+2)_q}{(k+1)_q}\,
\frac{\bigl((k+2)_q+k_q\bigr)}{((k+2)_q-k_q)^2}.
$$

As for the algebra $\slqc=\SL^C(R_q, 0)$, besides one-dimensional
representations, it possesses two Verma type representations
which differ from one another by the sign. Finding these
representations is left to the reader.

Now, we want to discuss different ways of quantizing vector
bundles on the sphere. According to the Serre-Swan approach any
vector bundle on an affine algebraic or a smooth variety can be realized as a
finitely generated projective module over the coordinate algebra of the variety.
As follows from \cite{R}, along with a formal deformation $A\rightarrow A_\h$ of
a commutative algebra $A$ any projective $A$-module $M$ can be also
deformed into an $A_\h$-module $M_\h$.
Otherwise stated, the idempotent $e(M)$ corresponding to the
module $M$ can be formally deformed into an idempotent
$e(M_\h)$ with entries
belonging to the algebra $A_\h$.

Nevertheless,
we are dealing with a non-formal deformation and our deformation
parameters can be specialized. So, we want to get explicit
expressions for idempotents over a quantum sphere (hyperboloid).
Construction of such idempotents in the framework of the first approach was
done in \cite{HM}. In that paper a series of idempotents $e_{\pm 1},\,e_{\pm 2},...$
was constructed which define the left $\K_q[SU(2)/H]$ modules
$\K_q[SU(2)/H]^{\oplus (|k|+1)} e_k$. Let us reproduce two of these idempotents:
$$
e_{-1}=\left(\begin{array}{cc}
ad&-\qq ab\\
cd&-\qq cb
\end{array}\right)=\left(\begin{array}{c}
a\\
c
\end{array}\right)(d, -\qq b),\quad e_{1}=\left(\begin{array}{cc}
da&-q dc\\
ba&-q bc
\end{array}\right)=\left(\begin{array}{c}
d\\
b
\end{array}\right)(a, -qc).
$$
(These formulae differ from those of \cite{HM} by replacing $q\to \qq$.)

The authors of \cite{HM} used this explicit realization of "quantum line bundles"
in order to compute the Chern-Connes index in a particular case: the projective
module corresponding to $k=-1$ and a representation of the quantum sphere
constructed in \cite{MNW}. Recall, that the index is defined by the pairing
of (the class of) a representation $\pi$ of a given algebra $A$ and
(the class of) a projective $A$-module $M$ in accordance with the rule
\be
\Ind(\pi, e)=\Tr (\pi(\Tr(e))),
\label{ind}
\ee
where $e$ is an idempotent corresponding to $M$.

Another approach to constructing projective modules over braided
varieties (orbits) was suggested in \cite{GS1}, \cite{GLS2}. Also, in
\cite{GLS2} a $q$-analog of the Chern-Connes index was
introduced and computed on $q$-spheres.
The main tool employed in this approach is a series of the Cayley-
Hamilton identities which are valid for the matrix $L_{(1)}=L$ and for its higher
analogs $L_{(k)}$, $k\ge 2$. Here $L$ is the matrix of the mREA generators and
$L_{(k)}$ is the $(k+1)\times (k+1)$ dimensional matrix whose explicit form is
given in \cite{GLS2}.

By passing to the algebra $\slhq$ we get reductions of the matrices
$L_{(k)}\to F_{(k)}$ and the Cayley-Hamilton identities for them.
Being reduced to the algebra $\slhqc$
these identities
take the form $p_k(F_{(k)})=0$ where $p_k$ is a polynomial with numerical
coefficients. Assuming the roots of the
polynomial $p_k$ to be distinct, we can associate with
it $k+1$ idempotents
$e_i(k),\ \ i=1,2,..., k+1$ with entries belonging to the algebra
$\slhqc$.

Consider an example: the idempotents arising from the matrix $L$
itself. In the algebra $\slhq$ the matrix $L$ reduces to
\be
F=\left(\begin{array}{cc}
{qh\over 2_q}&b\\
c&-{\qq h\over 2_q}
\end{array}\right)=b\left(\begin{array}{cc}
0&1\\
0&0
\end{array}\right)+ h\left(\begin{array}{cc}
{q\over 2_q}&0\\
0&-{\qq \over 2_q}
\end{array}\right)+c \left(\begin{array}{cc}
0&0\\
1&0
\end{array}\right).
\label{Lmat}
\ee
This matrix satisfies the Cayley-Hamilton identity
\be
F^2-\qq\h\, F-{\Cas_{sl}\over 2_q}\,\Id=0.
\label{CHid}
\ee
(This is an $sl$-reduction of the Cayley-Hamilton identity for the initial
non-reduced form of the matrix $L$.) While we pass to the algebra $\slhqc$ in
the space $V_k$ the coefficient ${\Cas_{sl}\over 2_q}$ becomes ${C \over 2_q}$.
Assuming the roots $\mu_i,\, i=0,1$ of the equation
$$
\mu^2-\qq \h \mu-{C \over 2_q}=0
$$
to be distinct we define two idempotents
$$
e_i(1)={L-\mu_i\Id \over \mu_i-\mu_{j}},\quad i,j=0,1, \; j\not=i.
$$

The $q$-index introduced in \cite{GLS2} has a form similar to
(\ref{ind}) but the trace $\Tr$ is replaced by its braided analog $\Tr_R$
$$
{\rm Ind}_q(\pi,e) = \Tr_R(\pi(\Tr_Re)).
$$
Note that the elements $\Tr_R\,e_i(k)$ are central whereas it is not true
for $\Tr \,e_i(k)$ if $q\not=1$. Also, observe that this method of
constructing projective modules over
braided orbits is valid for the mREA related to skew-invertible
Hecke symmetries of general type
(but the case of non-even symmetries is less studied).

The following proposition was proved in \cite{GLS2}.

\begin{proposition} Let $k$ be a sufficiently large positive integer. Then,
with a proper numbering of the idempotents $e_i(m), \, 0\leq i\leq m $ we have
$$
\Ind(\pi_k,e_i(m))=(m+k-2i+1)_q
$$
where $\pi_k$ is the (right) $R$-invariant representation of
the mREA $\loq$ in the space $V_k$.
Here we assume the $R$-trace to be normalized so that
it is an additive-multiplicative functional on the corresponding
Schur-Weyl category (see \cite{GLS1}).
\end{proposition}

Emphasize that we are dealing with a $q$-noncommutative
quantum (braided) hyperboloid algebra. On setting $q=1$, we get a
noncommutative hyperboloid
algebra which is a quotient of $U(sl(2))$. However, since the value
of the Casimir element (both in the classical and quantum cases)
depends on the space where the algebra in question is
represented we are somewhat dealing with a series of (braided)
hyperboloids. Namely, in the formula
for the $q$-index the idempotent $e_i(m)$ has entries
belonging to the algebra $\loqc$ with $C=C(V_k)$.

Note that the modules
$$
M_m=e_0(m)\loqc^{\oplus (m+1)} \quad {\rm and}\quad M_{-m}=
e_m(m)\loqc^{\oplus (m+1)}, \quad m=1,2,...
$$
correspond to the line
bundles $\O(m)$ and $\O(-m)$ respectively over the projective
space ${\Bbb{CP}}^1$ (also, we put $M_0=\loqc$, this module
corresponds to the
 trivial line bundle). We refer the reader to \cite{GLS2, GS2} for
 detail.

Now, we consider a quantum analog of the cotangent vector
bundle on the hyperboloid also presented as a projective module.
First, consider the usual hyperboloid
$$H^2=\{b,h,c\in sl(2)^*\,|\, bc+{h^2\over 2}+ cb=C\not=0\}$$
The space $\bigwedge^1(H^2)$ of one-forms on it consists of all
linear combinations
$$db\,\al+dh\, \beta+dc\, \gamma\qquad \al,\, \beta,\,\gamma \in
\K[H^2]$$ modulo the submodule
$$
(db\, c+{dh\, h\over 2}+dc\, b)\varphi,\,\,\varphi\in \K[H^2].
$$
Hereafter $\K[H^2]=\K[\R^3]/\langle 2bc+{h^2\over 2}-C\rangle$
 is the coordinate ring of the hyperboloid $H^2$. So, we
realize the space $\bigwedge^1(H^2)$ as a right $\K[H^2]$-module
but since the algebra $\K[H^2]$ is commutative it can be endowed
with a two-sided module structure. It is not difficult to show that this
module is projective.

By passing to a $q$-analog $\K_q[H^2]=\slq/\langle \qq
bc+{h^2\over 2_q}+qcb-C\rangle$ of the algebra $\K[H^2]$ we
naturally define a $q$-analog of the space ${\bigwedge}^1(H^2)$ as a
quotient-module
$$
{\bigwedge}^1_q(H^2)= \K_q[H^2]^{\oplus 3}/{M'}
$$
where $M'=e'\K_q[H^2]^{\oplus 3} $ is the right $\K_q[H^2]$-module
such that the corresponding idempotent is
$$
e'={1\over C}(\qq c,\, {h\over 2_q},\, qb)^{\sf T}(b,\, h,\,
c).
$$

Otherwise stated, in the right $\K_q[H^2]$-module
${\bigwedge}^1_q(H^2)$ the relation
\be
\qq db\, c+{dh\, h\over 2_q}+q dc\, b=0
\label{diff}
\ee
is imposed. It is also possible to realize the space ${\bigwedge}^1_q(H^2)$
as a left module $\K_q[H^2]^{\oplus 3}/M''$ where $M''=
\K_q[H^2]^{\oplus 3}e'' $ and
\be
e''={1\over C}(c,\, {h},\, b)^{\sf T}(\qq b,\, {h\over
2_q},\, q c).
\label{eprim}
\ee

This means that in the space of braided differentials we impose
the relation
$$\qq b\, dc+{h\, dh\over 2_q}+q c\, db=0.$$
Similarly to (\ref{diff}) this relation follows from the $q$-Casimir: in
order to get these relations we replace the left (resp., right) factors
in each summand of the $q$-Casimir by their differentials.

Emphasize that we are only dealing with one-sided modules
without endowing them with a two-sided module structure.
In order to convert such a one-sided module into a two-sided one
we have to introduce a transposition of
"functions" and "differentials" arising from the initial braiding $R$.
However, as was shown in \cite{AG2} there is no such a transposition
which would ensure the good deformation property of the module.

In a similar way we can treat the space ${\bigwedge}^2_q(H^2)$ of
quantum 2-differentials as a one-sided $\K_q[H^2]$-module. Note
that the same approach is valid on a $q$-noncommutative hyperboloid
(sphere) or its classical analog $q=1$ (see \cite{AG1}).
Observe that the Leibnitz rule for the differential $d$ is not applicable
in all these cases. However, it is possible to construct a $q$-analog of
the classical de Rham operator by analyzing the
decomposition of the spaces ${\bigwedge}^i_q(H^2),\,\, i=1,2$ into
a direct sum of
irreducible $\uq$-module (and similarly for other algebras) and
defining this $q$-analog via the classical pattern (see \cite{AG1} for
detail).

\section{Differential calculus via Koszul complexes}
\label{sec:7}

In this section we discuss the role of Koszul type complexes in
constructing the differential calculus on the quantum matrix algebras.

Let $V$ be a vector space over the ground field $\K$, $T(V)=
\bigoplus_{k=0}^{\infty} V^{\otimes k}$
 be its free tensor algebra and $I\subset \vv$ be a vector
 subspace. Consider the quadratic algebra $A=T(V)/\langle I
 \rangle$ and introduce
 the following subspaces
 $$
I^{\cap n}=V^{\ot (n-2)}\ot I \bigcap V^{\ot (n-3)}\ot I\ot V \bigcap
V^{\ot (n-4)}\ot I\ot \vv\bigcap...\bigcap I\ot V^{\ot (n-2)}
\subset V^{\ot n},\,\, n\geq 3.
$$
Also we set by definition $I^{\cap 2} = I$, $I^{\cap 1} = V$ and $I^{\cap 0} = \K$.

Then the Koszul complex is defined by the chain of maps
$$
...\to I^{\cap n}\ot A\to I^{\cap (n-1)}\ot A \to... \to I\ot A\to V\ot
A\to A\to\K\to 0
$$
where $A\to \K$ is the counit and for $n\geq 1$ the differential
$d:I^{\cap n}\ot A\to I^{\cap (n-1)}\ot A $ reads
$$
d(y_{i_1} y_{i_2}...y_{i_n}\ot x)=y_{i_1} y_{i_2}...y_{i_{n-1}}\ot
y_{i_n} x\quad \forall x\in A,\,\, y_{i_1},..., y_{i_n}\in V.
$$

In fact, this complex splits up into disjoint complexes
$$
...\to I^{\cap n}\ot A^{(m)}\to I^{\cap (n-1)}\ot A^{(m+1)} \to... \to
I\ot A^{(m+n-2)}\to V\ot A^{(m+n-1)}\to A^{(m+n)}\to 0
$$
where $A^{(m)}=V^{\ot m}/I^{\cup m}$ and
$$
I^{\cup m }=V^{\ot (m-2)}\ot I + V^{\ot (m-3)}\ot I\ot V + V^{\ot
(m-4)}\ot I\ot \vv+...+I\ot V^{\ot (m-2)}\subset V^{\ot m},\,\, m\geq 2.
$$
Also, we put $A^{(1)}=V$ and $A^{(0)}=\K.$
The algebra $A$ is called {\em Koszul} if these complexes are
acyclic for $m+n\geq 1$.

Assume now that in the space $\vv$ there exists a subspace
$I_+$ complimentary to $I$ (which in the sequel will be denoted $I_-$)
such that for a given $m\geq 3$ the subspaces $I_-^{\cup m }
\subset V^{\ot m}$ and $I_+^{\cap m }\subset V^{\ot m}$ are
complimentary, i.e.
\be
I_-^{\cup m }\bigcap I_+^{\cap m }=\{0\}\quad {\rm and}\quad
I_-^{\cup m }\oplus I_+^{\cap m }=V^{\ot m}. \label{cond} \ee
Let $P_{+}^{(m)}:V^{\ot m}\to I_{+}^{\cap m }$ be
the projector taking $I_{-}^{\cup m }$ to 0. Then for any element
$x$ of the homogeneous component $A^{(m)}$, $m\ge 2$, we have $x=P_{+}
^{(m)}(x)\in I_+^{\cap m }$ modulo $I_{-}^{\cup m }$.
If $I_+$ (respectively $I_-$) is subspace of symmetric (respectively skew-symmetric)
elements then the projector $P_+^{(m)}$ is the operator of the complete
symmetrization of elements from $A^{(m)}$.
Presenting any element of the algebra $A$ as a sum of
homogeneous components we can symmetrize it as well.
In this case we say that the element is presented in the {\em
canonical form}.

Now, consider a family of quadratic algebras $A(\nu)=T(V)/\langle
I(\nu) \rangle$ where $I(\nu)\subset \vv$ is a subspace depending
on a parameter $\nu$.
Since we are dealing with a non-formal deformation, we assume
that $I(\nu)$ is defined by a system of finite linear combinations of
generators of
$\vv$ with coefficients analytically depending on $\nu$ in a neighborhood
$U$ of $0$. So, in $U$ the parameter $\nu$ can be specialized.
We are interested in such families of algebras that
$\dim\, A^{(m)}(\nu)=\dim\, A^{(m)}$ for any integer $m$
(at least for a generic $\nu\in U$), where $A=A(0)$  .

Assume that for a family $I_-(\nu)=I(\nu)$ there exists another one
$I_+(\nu)$ such that for any $m\geq 2$ the subspaces
$I_-^{\cup m }(\nu)$ and $I_+^{\cap m }(\nu)$ are complementary
and the projectors $P_+^{(m)}(\nu)$ analytically depend on $\nu \in
U$ except my be a
finite set of points which does not include 0. Then for a generic
$\nu\in U$ we have $\dim\, A^{(m)}(\nu)=\dim\, A^{(m)}$ for all $m$,
i.e. the algebra $A$ have the good deformation
property.

Given families $I_{\pm}(\nu)$ and therefore the projector $P_{+}
^{(2)}(\nu)$, we look for
a higher projector $P_{+}^{(m)}(\nu),$ $m\geq 3$ as a polynomial
in $(P_{+}^{(2)}(\nu))_{12},\,...,\,(P_{+}^{(2)}(\nu))_{m-1\,m}$
(the subscripts indicate the spaces in the product $V^{\ot m}$
where the operator acts) with analytical coefficients.
Having constructed such a polynomial, we can conclude that for
a generic $\nu$ $\dim\, A^{(m)}(\nu)=\dim\, A^{(m)}$ (here "a
generic $\nu$" means "all $\nu$ except for a countable set").
If such projectors exist for all $m\geq 3$ we conclude that this
property is valid for all homogeneous components of the algebra
$A(\nu)$.

If in addition the initial algebra $A$ is Koszul, then this is also true for the
algebra $A(\nu)$ for a generic $\nu$ (see \cite{PP}).

In examples below we have two families of subspaces $I_{\pm}
(\nu)\subset \vv$ and we want to show that
the both algebras
$$
A_+(\nu)=T(V)/\langle I_-(\nu) \rangle,\qquad A_-(\nu)=T(V)
/\langle I_+(\nu) \rangle
$$
have the good deformation property. In order to show this,
we should construct the projectors $P_{\pm}^{(m)}(\nu),$ $ m\geq 3$
in terms of the  operators $P_{\pm}^{(2)}(\nu)$.
The "skew-symmetrization" operators $P_{-}^{(m)}(\nu)$ are defined similarly to the projectors
$P_+^{(m)}(\nu)$, provided that the subspaces $I_+^{\cup m}(\nu)$
and $I_-^{\cap m}(\nu)$ are complementary to each other for any $m\geq 2$.
We call a couple of the subspaces
$I_{\pm}(\nu)$ {\em regular} if the subspaces $I_{\pm}^{\cup m}(\nu)$ and
$I_{\mp}^{\cap m}(\nu)$ are complementary for any $m\geq 2$.

In \cite{G2} this scheme was applied to the algebras $\Sym_q(V)$
and $\bigwedge_q(V)$ related to the Hecke symmetries. Namely,
a series of projectors $P_\pm^{(m)}(q)$ related to $I_-(q)=\Im(q\Id-R_q)$
and $I_+(q)=\Im(q^{-1}\Id+R_q)$ was constructed. As follows from \cite{G2},
these algebras have the good deformation property. Note that if a family of Hecke
symmetries is not quasiclassical the monomials $x_1^{a_1}...x_n^{a_n}$
do not form a basis in the algebra $\Sym_q(V)$. So, the method of
verifying the good deformation property of this algebra based on the
ordering the generators is not valid any more, whereas the above scheme is still
applicable.

Again, consider the algebras $\Sym_q(\TT)$ and $\lq=\Sym_q(\LL)
$ defined by formulae (\ref{RTT}) and (\ref{RE}) respectively
(we do not assume a Hecke symmetry $R$ to be quasiclassical).
This means that the corresponding subspaces $I_-(\TT)$ and $I_-(\LL)$
are determined by the left hand side of these formulae.
Now, define the complementary "symmetric" components
$I_+(\TT)$ and $I_+(\LL)$ by putting
$$
I_+({\TT})=RT_1T_2+T_1T_2 R^{-1},\qquad I_+({\LL})=
RL_1RL_1+L_1RL_1 R^{-1}
$$
and define the corresponding algebras
\be
\bw_q(\TT)=T(\TT)/\langle I_+({\TT})\rangle , \qquad \bw_q(\LL)=
T(\LL)/\langle I_+({\LL})\rangle.
\label{def:ext-alg}
\ee

Note that the projector $P^{(2)}_{+}(\LL)$ is nothing but the operator
${\cal S}_q$ discussed in section \ref{sec:5}. In the paper \cite{GPS3}  an attempt
was undertaken to construct the higher projectors $P^{(m)}_{\pm}$
via those $P^{(2)}_{\pm}$ in order to apply the above scheme. We
have only succeeded in constructing these projectors for $m=3$.
(In fact, the construction is valid for all quantum matrix algebras
associated with a compatible pair of Hecke symmetries, in particular,
for $\Sym_q(\TT)$.) Nevertheless,
as follows from \cite{Dr2} this property suffices
for concluding that for any $m\geq 4$ and a generic $q=e^{\nu}$
the dimensions of all homogeneous components
are stable (classical for a quasiclassical $R$). So, according to
\cite{Dr2} we should only control the dimensions of the third
homogeneous components. However, it would be interesting to
explicitly construct the higher projectors.

Consider one example more, where this scheme can be hopefully
applied, namely, the algebra $\slq=\lq/\langle \ell \rangle$
(an analogous quotient of the algebra $\Sym_q(\TT)$ cannot be
defined since there is no central element in $\TT$).
As we observed above the algebra $\slq$ is defined by the same
formulae as that $\lq$ but with the generators $F_i^{\,j}$ instead of
$L_i^{\,j}$. This means that the subspace $I_-({\SL})$ corresponding to this algebra
can be obtained by replacing the generators
$L_i^{\,j}$ in $I_-({\LL})$ by their traceless components $F_i^{\,j}$
defined by formula (\ref{shift}).

It is more difficult to define the result of "$sl$-reduction" of
the algebra $\bigwedge_q(\LL)$. Describe this procedure by mainly following
\cite{IP1}.

In accordance with definition (\ref{def:ext-alg}) the algebra
$\bigwedge_q(\LL)$ is generated by the elements
$L_i^{\,j}$ subject to the following multiplication rules
\be
R_{12}L_1R_{12}L_1+L_1R_{12}L_1 R_{12}^{-1} = 0.
\label{ext-alg}
\ee
We make a linear change (\ref{shift}) $L_i^{\,j}\to \{F_i^{\,j},\ell\}$
including explicitly the $R$-trace $\ell=\Tr_R(L)$
into the set of generators. Of course, as well as in the case of
mREA, we assume $\Tr_R(\Id) = \Tr(C)\not=0$.

Calculating the $R$-trace in the second space of the defining relations
(\ref{ext-alg}) we find
$$
\ell\,L+L\,\ell +\omega L^2 = 0,\qquad \omega = q-q^{-1}
$$
or
\be
\ell\,F +F \,\ell +\frac{\omega}{1+\omega(\Tr_R(\Id))^{-1}}\,F^2 = 0,
\label{ell-F}
\ee
that is the $R$-trace $\ell$ is not central in our algebra. As was
shown in \cite{IP1}, the element $\ell$ is nilpotent: $\ell^2 = 0$.

Now we can find the multiplication table of the algebra
$\bigwedge_q(\LL)$ in terms of
the generators $F$ and $\ell$. Substituting the expression of $L$
through $F$ and
$\ell$ into (\ref{ext-alg}) and using the anticommutation relation
(\ref{ell-F}) and nilpotency
of $\ell$, we get the following result
\be
\begin{array}{l}
R_{12}F_1R_{12}F_1+F_1R_{12}F_1 R_{12}^{-1} =
\kappa\,(F_1^2+R_{12}F_1R_{12}),\\
\rule{0pt}{7mm}
\ell\,F+F\,\ell = -\tau\,F^2,\\
\rule{0pt}{7mm}
\ell^2 = 0,\qquad \Tr_R(F) = 0,
\end{array}
\label{s-f}
\ee
where the numeric parameters $\kappa$ and $\tau$ read
$$
\kappa = \frac{\omega}{\Tr_R(\Id) + \omega},\quad
\tau = \kappa\,\Tr_R(\Id).
$$

Note, that the $R$-traceless elements $F_i^{\,j}$ forms a {\it
subalgebra}
of $\bigwedge_q(\LL)$ since the element $\ell$ does not enter their
multiplication table contrary to the case of mREA (see the relations
(\ref{F-l})).
So, in this case we have no need to pass to a quotient algebra in
order
to obtain the $sl$-reduction --- the traceless algebra is a subalgebra
of the initial one. Finally, we put
$$
I_+(\SL)=\span(R_{12}F_1R_{12}F_1+F_1R_{12}F_1 R_{12}^{-1} -
\kappa\,(F_1^2+R_{12}F_1R_{12}))
$$
and $\bw_q(\SL)=T(\SL)/\langle I_+(\SL) \rangle$.

Hopefully, the couple of subspaces $I_{\pm}({\SL})$ is regular. We
are able to prove this claim if the initial Hecke symmetry is of the
Temperley-Lieb type. A proof will be given in our subsequent paper.

Here we consider an example arising from the standard Hecke
symmetry. The space $I_-(q=e^{\nu})$ is determined by the left hand side
of the relations (\ref{sys3}). It is the spin 1 $\uq$-submodule of the space $\SL^{\ot 2}$
endowed with an $\uq$-action.
The corresponding algebra is just $\slq$.
The space $I_+(q)=V_2\oplus V_0$ in this case is a direct sum of
the spin 0 and spin 2 $\uq$-submodules of $\SL^{\ot 2}$.
Here
\be
\begin{array}{c}
\displaystyle
V_0={\span}( \qq b\ot c +{1 \over{2_q}}\,h\ot h + q c\ot b),\\
\rule{0pt}{5mm} V_2={\span} (b\ot b,\;q^2b\ot h+h\ot b,\; q^3b\ot
c-q h\ot h+ \qq c\ot b,\; q^2h\ot c+c\ot h,\; c\ot c).
\end{array}
\label{Iplus}
\ee
It is not difficult to see that the braiding in the space $\SL^{\ot 2}$
which is the extension of the initial Hecke symmetry is a BMW symmetry.
With the use of methods of \cite{OP2} it is possible to construct the projectors
$P_{\pm}^{(m)}$. This is the basic idea of the aforementioned proof.

\begin{remark}{\rm
Given a subset $I_-\subset \vv$ it is not clear
whether there is a subspace $I_+$ such that the
 couple $(I_-,\,I_+)$ is regular. But even if it is the case,
the complementary subspace $I_+$ is not in general unique.
Consider an example.

Let $R:\vv\to\vv,\,\,(\dim V=2)$ be an
involutive symmetry given by
$$R(x\ot x)=x\ot x,\,\,R(x\ot y)=b\, x\ot x+y\ot x,\,\, R(y\ot x)=-b\, x\ot
x+x\ot y,$$
$$R(y\ot y)=a\,b\, x\ot x-a\, x\ot y+a\, y\ot x +y\ot y$$
where $\{x, y\}$ is a basis in $V$ and $a, \, b\in \K$.
Then $I_-=\span(-b\, x^2+xy-yx)$ and $I_+=\span(x^2, xy+yx, -a\,
xy+y^2)$ (we omit the sign $\ot$). The couple $(I_-,\,I_+)$ is
regular for any
$a,\, b \in \K$.
Even if we put $b=0$, i.e. if we consider the "classical"
skew-symmetric subspace $I_-$, the space $I_+$ is not unique and
depends on $a$.}
\end{remark}

Assume $(I_-,\,I_+)$ to be a regular couple. Consider the
algebras $A_+=A=T(V)/\langle I_- \rangle$ and $A_-=T(V)/\langle
I_+ \rangle$ and associate with them two Koszul complexes
$$ d_-:A_-^{(m)}\ot A_+^{(n)}\to A_-^{(m+1)}\ot A_+^{(n-1)},\qquad
d_+:A_-^{(m)}\ot A_+^{(n)}\to A_-^{(m-1)}\ot A_+^{(n+1)}$$
where we identify $A_+^{(n)}$ and $I_+^{\cap n}$ (resp.,
$A_-^{(m)}$ and $I_-^{\cap m}$). In \cite{G2} these complexes
associated with a Hecke symmetry were called Koszul complexes
of the first kind. Now, we want to use them in order to define a
differential and "partial derivatives" on quantum matrix algebras.

Introduce a de Rham-Koszul differential $\ddk$ on the component
$A_-^{(m)}\ot A_+^{(n)}$ by setting $\ddk=n\,d_-$ (the factor $n$ is
motivated by an analogy with
the classical case). Thus, we have
$$\ddk (y\ot x_{i_1}x_{i_2}...x_{i_n})=n\, y x_{i_1}\ot
x_{i_2}...x_{i_n},\,\,y\in A_-^{(m)},\,\,x_{i_1},x_{i_2},...,x_{i_n}\in V.$$
In particular, if $m=0$ we treat this operator to be an analog of the
de Rham differential on the algebra $A_+=\bigoplus A_+^{(n)}$.

Now, define {\em braided partial derivatives} on this algebra in the
standard way. If $f\in A_+^{(n)}$ is a homogeneous element, consider its image
$\ddk (f)=\sum_i x_i\ot f_i$ where $f_i\in A_+^{(n-1)}$. Then we introduce
the braided partial derivative in $x_i$ by putting $\partial_{x_i}(f)=
f_i$. The differential $\ddk$ and the braided partial derivatives
can be extended on the whole algebra $A_+$ by linearity. Below, we use the
notation $\partial_i^j=\partial_{L_j^{\,i}}$ for the braided derivatives in the
generators of the REA. Thus, we have $\partial_i^j\, L_k^{\,l}=\de_k^j\, \de_i^l$.
Below we shall omit the term "braided".

This way of introducing the partial derivatives on quantum algebra
does not use any form of the Leibnitz rule. Emphasize that we are
only dealing with one-sided $A_+$-modules without transposing "functions" and
"differentials". Instead, we apply the de Rham-Koszul differential
and the partial derivatives to elements presented in the canonical form. So, we do not need to
verify any compatibility of the differential with the
defining relations of the algebras in question. Below we use the
same principe in order to define other "braided vector fields".

Again, let $(I_-,\,I_+)$ be a regular couple of subspaces of $\vv$.
Now, consider another complex (called the Koszul complex of the
second kind in \cite{G2}). To this end we need the space $V^*$ dual to $V$.
This means that that there exists a nondegenerated pairing
$V^*\ot V\to \K$. We extend this pairing on the spaces $(V^*)^{\ot k}
$ and $V^{\ot k}$ by the rule
$$
\langle a\ot b,\, c\ot d\rangle=\langle b,\, c \rangle\langle a,\,d
\rangle, \quad a,b\in V^*,\,\, c,d \in V
$$
and so on. Define the subspaces $I^*_-=(I_+)^\bot$ and $I^*_+=
(I_-)^\bot$ where $I^\bot\subset (V^*)^{\ot 2}$ stands for the space
orthogonal to $I\subset \vv$. Also, we put
$$
A^*_+=T(V^*)/\langle I^*_- \rangle\quad {\rm and}\quad A^*_-=
T(V^*)/\langle I^*_+ \rangle.
$$

Observe that the couple $(I_-^*, I_+^*)$ is regular and present all
elements of the algebras $A_+$ and $A_-^*$ in the canonical
form. Now introduce a differential
$$
\tilde{d}: {A_-^*}^{(m)}\ot A_+^{(n)}\to {A_-^*}^{(m-1)}\ot
A_+^{(n-1)},\,\, m,n\geq 1
$$
via the pairing $\langle\,\,,\,\rangle_{m,\, m+1}$. This means that the
last factor of $(A_-^*)^{(m)}$ and the first factor of $A_+^{(n)}$
are coupled. It is not difficult to see that $\tilde d^2 = 0$.

In what follows we  consider the operator $\de=n\, \tilde{d}$
(we have renormalized the initial differential by the same reason as
above). In \cite{G2} it was shown that the Koszul complex of the second kind
associated with a skew-invertible Hecke symmetry is acyclic for a
generic $q$. Hopefully, it is also so for the complexes associated with all
quantum matrix algebras in question.

In the next section we use the both kinds of
Koszul type complexes in order to define q-analogs of wave operators.

Now, compare our definition of the derivatives on the algebra $\lq$
with that from \cite{Me1}, \cite{Me2}. To this end, besides the
operator $Q$ (\ref{Q-op}) we also introduce $Q'$
$$
Q'(L_{\overline 1}\stackrel{.}{\ot} L_{\overline 2}) = R^{-1}_{12}\, L_{\overline
1}\stackrel{.}{\ot} L_{\overline 2}\, R_{12}^{-1}.
$$
Then we have
$$I_-(\L)=\Im(\Id-Q),\quad I_+(\L)=\Im(\Id+Q').$$

It is clear that these operators commute with each other. Also, the
both operators satisfy the quantum Yang-Baxter equation, i.e. they
are braidings.
Besides, the couples $\{Q,Q'\}$ and $\{Q',Q\}$ are compatible in
the sense of the definition of section \ref{sec:3} and $(\Id-Q)(\Id+Q')=0$.
For such a couple of operators it is possible to apply the scheme
from \cite{M1} where some aspects of the differential calculus on a
space $V$ endowed with a skew-invertible Hecke symmetry were
developed\footnote{Some aspects of
such type calculus  in a particular case related to a QG
were previously considered in \cite{WZ}.}. Namely, let $W$ be a vector space such that in
$W^{\ot 2}$ there are defined two
operators $Q$ and $Q'$ satisfying the above conditions. Then the
partial derivatives in the algebra $\Sym_q(W)=T(W)/ \langle \Im(\Id-Q) \rangle$
can be defined as follows
$$
\partial^i x_j=\de _j^i,\quad \partial^i (x_j\, x_k)=\partial^i_1
(\Id+Q')(x_j\, x_k),
\quad
\partial^i (x_j\, x_k\, x_l)=\partial^i_1 (\Id+Q'_{12}+ Q'_{12} Q'_{23})
(x_j\, x_k\, x_l)
$$
and so on. Here, $\{x_i\}$ is a basis of $W$ and  $\partial^i_1$ stands for the derivative
in the generator $x_i$ applied to
the first factor.

Observe that  we do not assume the braidings $Q$ and $Q'$
to be skew-invertible. We do not need this property since we do
not transpose
the derivatives and the generators $x_i$. So, this way of
proceeding differs from the usual Leibnitz rule.

Now, by assuming $W=\LL$ apply these partial derivatives to a degree $n$ homogeneous
element $f\in \LL^{(n)}(R_q)$ written in the canonical form.
Then we have
$$
(\Id+Q'_{12}+ Q'_{12} Q'_{23}+...+ Q'_{12} Q'_{23}...Q'_{n-1\,
n})(f)=n\, f.
$$
This follows from the fact that
$$
Q'_{k-1\, k}P_+^{(n)}=P_+^{(n)}Q'_{k-1\, k}=P_+^{(n)},\qquad
2\le k\le n,
$$
where $P_+^{(n)}$ is the "symmetrization" projector. Note that this
relation follows from the minimal polynomial
for the operator $Q'$ (which can be found similarly to that for $Q$)
and  the fact that a bigger projector "absorb" a smaller one:
$$
(P_+^{(2)})_{k-1\, k}P_+^{(n)}=P_+^{(n)}(P_+^{(2)})_{k-1\, k}=P_+^{(n)},
\qquad 2\leq k\leq n.
$$
This implies that our definition of partial derivatives and that from
\cite{Me1} are equivalent.

Let us go again to a general space $W$ endowed with operators $Q$ and $Q'$ and  define
a $q$-analog of the de Rham differential in the same manner. Namely, we put
$$
d(x_i\, x_j)= (d\,x_i)\, x_j+ (d\ot \Id) (Q'(x_i\, x_j))
$$
and so on. So, the space of one-differentials is realized as a right
module over the algebra $\Sym_q(W)$.
Introducing the skew-symmetric algebra
$$
\bw_q(W)=T(W)/ \langle \Im(\Id+Q') \rangle
$$
we treat the space $\bw_q(W)\ot \Sym_q(W)$ as the space
of all differentials on $W$. An analog of the de Rham operator on this
algebra can be introduced in the same way as above.

We complete this section with the following comment. The operators
$Q,\, Q',\, R_{\End}$ appeared in the early 90's (see \cite{MMe} and
the references therein). They play a very important role in the $q$-analysis.
In particular, each of the operators $Q$ and $Q'$ enables us to introduce
the subspaces $I_{\pm}(\LL)\subset \LL^{\ot 2}$ and therefore those
$I_{\pm}(\SL)\subset \SL^{\ot 2}$. However, it is not clear whether the
latter subspaces can be defined via the operator $R_{\End}$.

\section{$q$-Wave operators on $q$-Minkowski space algebra}
\label{sec:8}

As we noticed in Introduction the $q$-Minkowski space algebra was
treated in some papers to be a particular case of the REA. However,
we want to slightly modify its definition in order to make it more similar to
its classical counterpart.

There exists a family of quadratic, $\uq$-covariant algebras with central
time variable $t=\ell$ which are deformations of $\krr$. We get such algebras by
replacing the factor $q-\qq$ of the central element $\ell$ in the system (\ref{sys2})
(with $\h=0$) by an arbitrary multiplier $\alpha$. A similar change can be done in
the algebra $\lq$ related to any skew-invertible Hecke symmetry.

Let us now set $\alpha=0$ and denote the corresponding algebra  $\tlq$.
We call this algebra the {\em truncated REA}. In a similar way we define the
{\em truncated mREA} $\tlhq$. Explicitly, these algebras are defined as
follows
$$
\tlhq=\slhq\ot \K[\ell]\quad {\rm and}\quad \tlq=\slq\ot \K[\ell].
$$
In the sequel we denote the algebra $\tlq$ in the $\uq$ case
as $\krrq$ and consider it to be the $q$-Minkowski space algebra.
Also, we put $\krq=\krrq/\langle \ell \rangle$.

Thus, we have $\tlq=T(\LL)/\langle \tim \rangle$ where
$$
\tim =I_-\oplus {\rm span}(\ell\ot b - b\ot \ell ,\; \ell\ot
h-h\ot \ell,\; \ell\ot c-c\ot \ell)\subset \LL^{\ot 2}
$$
and $I_-=I_-(q)$ is the subspace spanned by the left hand side of
(\ref{sys3}). More explicitly, $\krrq$ is generated by four
generators $\{\ell, h, b, c\}$ subject to the relations
$$
q^2hb-bh=0,\quad q^2c h-hc=0,\quad 2_q q (bc-cb)+(q^2-1) h^2=0,
\quad \ell f=f\ell,\quad \forall\,f\in \span(b,h,c).
$$

Besides, consider the subspace
$$
\tip=I_+\oplus {\rm span}(\ell\ot b + b\ot \ell,\; \ell\ot
h+h\ot \ell,\; \ell\ot c+c\ot \ell,\; \ell\ot \ell)\subset
\LL^{\ot 2}
$$
where $I_+=I_+(q)=V_0\oplus V_2$ (see (\ref{Iplus})).

Note that regularity of the couple $(I_-,\,I_+)$ entails the same
property for $(\tim,\, \tip)$. Therefore, we can apply the scheme of the previous
section to the algebras $\krrq$ and $\twl=T(\LL)/\langle \tip \rangle$ and define
the partial derivatives $\partial_i^j$ on the algebra $\krrq$.

Also, we need the $R$-invariant pairing  on the space $\LL$.
Such a pairing can be defined by the formula
\be
\langle\,,\,\rangle:\LL^{\ot 2}\to \K,\quad \langle L_i^{\,j}, L_k^{\,l} \rangle=
\de_i^l\, B_k^j
\label{spariv}
\ee
where the operator $B$ was introduced in (\ref{def:B}).
So, the space $\LL$ can be identified with its dual.

However, such a pairing is not unique. It becomes  unique (up to a factor)
being restricted to the space $\SL$. In the standard case the pairing $\SL\ot \SL \to \K$
can be chosen with the following  normalization
\be
\langle b,c \rangle=\qq,\quad \langle h, h \rangle=2_q,\quad
\langle c, b \rangle=q
\label{q-pair}
\ee
(all other terms are trivial). Let us extend this pairing to the space $\LL$ by putting
\be
\langle \ell, \ell \rangle=\epsilon^{-1},\quad \langle \ell,
f\rangle=\langle f, \ell \rangle=0\,\,\,\,\, \forall f\in \SL,\quad
\epsilon\in \K, \epsilon\not=0
\label{pair}
\ee
This is the most general form of a nondegenerated $\uq$-invariant pairing on the
space $\LL$.

Then one can extend this pairing to the spaces $\SL^{\ot k}\ot \SL^{\ot k}$,
$k\ge2$ as was explained in the previous section. The following proposition can
be proved by straightforward computations.

\begin{proposition}
\label{pro}
The spaces $I_-$ and $I_+$ are orthogonal to each other with respect to the pairing
{\rm (\ref{q-pair})}. It is also true for the spaces $\tim$ and $\tip$ and the pairing
{\rm (\ref{q-pair})--(\ref{pair})}.
\end{proposition}

For the extended pairing we keep the same notation.

\begin{definition} We say that a triple $(I_-,\, I_+ \subset \vv, \,
\langle \,\,,\,\rangle:\vv\to \K)$ is regular if the couple
$(I_-,\, I_+)$ is regular and the subspaces $I_-$ and $I_+$ are
orthogonal to each other with respect to the pairing $\langle
\,\,,\,\rangle$.
\end{definition}

So, the proposition \ref{pro} states that the above triples are regular.
This property enables us to identify the algebra $A_-$ and that $A^*_-$ in notations of the
previous section and consequently to consider the differentials $\ddk$ and $\de$ as acting on the
same terms ${A_-}^{(m)}\ot A_+^{(n)}$. Thus, the operator
$$
\de\, \ddk: {A_-}^{(m)}\ot A_+^{(n)}\to {A_-}^{(m)}\ot A_+^{(n-2)}
$$
is well defined.

Also, this property of the subspaces $\tim$ and $\tip$ enables us to
compute the commutation relations between the partial derivatives
$\partial_{\ell},\, \partial_h,\, \partial_b,\, \partial_c$
acting on the algebra $\krrq$. To this end we consider another basis
$\{D_{\ell},\, D_c,\, D_h,\, D_b\}$ in the space spanned by
these derivatives. Namely, we put
$$
D_{\ell}=\partial_{\ell},\quad D_c=q\,\partial_b, \quad D_h=
2_q\,\partial_h, \quad D_b=\qq\,\partial_c.
$$

Note that this basis is more convenient to deal with since
the map $b\to D_b, h\to D_h,\, c\to D_c, \, \ell\to D_{\ell}$ is R-invariant whereas
that $b\to \partial_b,...$ is not. Explicitly, these operators are defined via the pairing :
$$
D_b(f)=n\,\langle b, f\rangle_{12},\quad {\rm if } \quad f\in \LL^{(n)}(R_q) \quad {\rm and \,\, so\,\, on}$$
where the subscript means that the pairing with $b$
is applied to the first factor of an element $f$.

The next proposition is a straightforward corollary of
proposition \ref{pro}.

\begin{proposition}
\label{prop:20}
The derivative $D_{\ell}$ commutes with
$D_c,\, D_h,\, D_b$. The derivatives $D_c$, $ D_h$ and $D_b$
satisfy the relations
{\rm
\begin{eqnarray}
&q^2D_h\, D_b-D_b\, D_h=0&\nonumber\\
& 2_q\, q (D_b\, D_c-D_c\, D_b)+(q^2-1)D_h\, D_h=0&\label{der-
alg}\\
&q^2 D_c\, D_h-D_h\, D_c=0.&\nonumber
\end{eqnarray}}
Otherwise stated, the map $\tau: b\to D_b,\, h\to D_h,\, c\to
D_c,\, \ell \to D_{\ell}$ is a representation of the algebra
$\krrq$.
\end{proposition}

This proposition entails that the operators
\be
\De_{\krq}= \qq
D_b\, D_c+{ D_h^2 \over 2_q} +q D_c\,
D_b=\qq\,\partial_c\,\partial_b+2_q\,\partial_h^2+q\,\partial_b\,
\partial_c,
\label{Lap1}
\ee
\be
\De_{\krrq}=\epsilon\,
D_{\ell}^2+ \qq D_b\, D_c+{ D_h^2 \over 2_q} +q D_c\, D_b=
\epsilon\, \partial_{\ell}^2+
\qq\,\partial_c\,\partial_b+2_q\,\partial_h^2+q\,\partial_b\,
\partial_c,
\label{Lap2}
\ee
are central. Here, $\epsilon$ is an
arbitrary non-trivial factor. We call them $q$-Laplace operators on
the algebras $\krq$ and $\krrq$ respectively. Thus, our $q$-Laplace operator
on the algebra $\krrq$ depends on $\epsilon$.

Note that these operators are obtained from the quadratic central elements in the
algebra $\krq$ and $\krrq$ respectively in which we replaced the generators by
the corresponding derivatives: $b\to D_b$ and so on.

Now, we pass to constructing $q$-Dirac operators on the algebras
$\krq$, $\krrq$ and the $q$-hyperboloid algebra $\khq$. All of these
operators can be defined via a universal scheme making use of the
{\em split} Casimir elements.
Note that this scheme is similar to that from \cite{GLS2} where a
way of relating the REA Cayley-Hamilton identities and the
split Casimir elements was suggested.

We begin with classical algebras. First, consider the $sl(2)$ {\em
split} Casimir element
$$
b\ot c + {h\ot h\over 2}+c\ot b\in sl(2)\ot sl(2).
$$
Let $\pi$ be the spin ${1 \over 2}$ representation of the algebra
$sl(2)$ in the standard basis. Applying this representation to the
right factors of the split Casimir and the map $\tau$ (with $q=1$)
defined in proposition \ref{prop:20} to the left factors we convert
the split Casimir element into the $2\times 2$ matrix. Its transposed
matrix reads
\be
D_b \left(\begin{array}{cc}
0 & 1\\
0& 0
\end{array}\right)+{1\over 2} D_h\left(\begin{array}{cc}
1 & 0\\
0& -1
\end{array}\right)+D_c \left(\begin{array}{cc}
0 & 0\\
1& 0
\end{array}\right)= \left(\begin{array}{cc}
{D_h\over 2} & D_b\\
D_c& -{D_h\over 2}
\end{array}\right).
\label{six}
\ee
It can be treated as the Dirac operator on the space $\R^3\cong
sl(2)^*$. It is easy to see that the square of this matrix is
$$({D_h^2\over 4}+D_b\, D_c)\,\Id=
(\partial^2_h+\partial_b\partial_c)\, \Id.$$

In order to get a similar operator on the space $\R^3\cong
so(3)^*$ equipped with the Euclidean coordinates
$$x={\ii (b+c) \over 2},\quad y={c-b \over 2},\quad z={\ii (a-d)\over
2}$$
we proceed in a similar way with the $so(3)$ split Casimir
element $-2\ii(x\ot x+y\ot y+ z\ot z$) (all numerical factors are
introduced
for our convenience). Then we arrive to the operator \be
\partial_x \left(\begin{array}{cc}
0 & 1\\
1& 0
\end{array}\right)+\partial_y \left(\begin{array}{cc}
0 & -\ii\\
\ii& 0
\end{array}\right)+\partial_z \left(\begin{array}{cc}
1 & 0\\
0& -1
\end{array}\right)=\left(\begin{array}{cc}
\partial_z & \partial_x-\ii \partial_y\\
\partial_x+\ii \partial_y& -\partial_z \end{array}\right). \label{lab} \ee
Being squared it equals
$(\partial_x^2+\partial_y^2+\partial_z^2)\, \Id$.

A passage to the 4-dimensional Minkowski space can be done in
the usual way (see formula (\ref{dir}) below).

To obtain the Dirac operator on the sphere
$S^2$ of the radius $r=1$ we replace the partial
derivatives in (\ref{lab}) by the infinitesimal rotations
$$
X=z\,\partial_y-y\,\partial_z, \;Y=x\,\partial_z-z\,\partial_x, \;Z=
y\,\partial_x-x\,\partial_y.
$$
Thus, we get the operator
$$
{\Dir}_{\K[S^2]}=\left(\begin{array}{cc}
Z & X-\ii Y\\
X+\ii Y& -Z \end{array}\right).
$$
It satisfies the equation
$$
{\Dir}_{\K[S^2]}^2=\ii\,{\Dir}_{\K[S^2]}+(X^2+Y^2+Z^2)\Id.
$$

In a similar way, on replacing the derivatives $D_b,\, D_h,\, D_c$ in
(\ref{six}) by
the hyperbolic infinitesimal rotations
$$
B=-2b\, \partial_h+h\,\partial_c,\quad H=2b\,
\partial_b-2c\,\partial_c,\quad C=-h\partial_b+2c\, \partial_h
$$
we can introduce the Dirac operator on a hyperboloid (see formula
(\ref{trinad}) below).

Now, we are going to use the same scheme for the braided split Casimir
$$
\qq\, b\ot c +{h\ot h\over 2_q}+q\, c\ot b\in \SL\ot \SL.
$$
Applying the representation (\ref{taut}) to the right factors
in this split Casimir we get a matrix whose transposed form coincides with
the matrix (\ref{Lmat}) (up to a numerical multiplier). So, by definition, the $q$-Dirac
operator on the algebra $\krq$ is the matrix (\ref{Lmat}) where we
assume that $\h=0$ and replace the left factors $b,\,h,\, c$ by
$D_b,\,D_h,\, D_c$ respectively (i.e. we apply the representation
$\tau$). Finally, we have
\be
{\Dir}_{\krq}=D_b \sigma_b+D_h
\sigma_h+D_c \sigma_c=\left(\begin{array}{cc}
{q\,D_h\over 2_q} & D_b\\
D_c& -{D_h \over 2_q\,q} \end{array}\right).
\label{dirr}
\ee
Here the matrices $\sigma_b,\, \sigma_h,\, \sigma_c$ are
respectively the multipliers of $b,\, h,\,c$ in formula
(\ref{Lmat}).

Using the Cayley-Hamilton identity (\ref{CHid}) for the matrix
(\ref{Lmat}) (where we put $\h=0$) we find that
$$
{\Dir}_{\krq}^2={1\over 2_q}(\qq \, D_b\, D_c+{1\over 2_q}
D_h^2+q\, D_c\, D_b) \, \Id.
$$

Now, define the $q$-Dirac operator on the algebra $\krrq$ in the
usual way by setting \be {\Dir}_{\krrq}=\epsilon D_{\ell}
\left(\begin{array}{cc}
0_2 & {{\rm I}}_2\\
{{\rm I}}_2& 0_2
\end{array}\right)+D_b \left(\begin{array}{cc}
\sigma_b & 0_2\\
0_2& -\sigma_b
\end{array}\right)+D_h\left(\begin{array}{cc}
\sigma_h & 0_2\\
0_2& -\sigma_h
\end{array}\right)+D_c \left(\begin{array}{cc}
\sigma_c & 0_2\\
0_2& -\sigma_c
\end{array}\right) \label{Dir} \ee
where $0_2$ and ${{\rm I}}_2$ are respectively the trivial and
unit $2\times 2$ matrices.

It is easy to see that \be {\Dir}_{\krrq}^2=(\epsilon^2\,
D_{\ell}^2+{1\over 2_q}(\qq \, D_b\, D_c+{1\over 2_q}D_h^2+q\,
D_c\, D_b))\Id \label{Dirr} \ee

In order to introduce a $q$-Dirac operator on the $q$-hyperboloid we
need braided analogs of the hyperbolic infinitesimal rotations.
They are defined via a braided analog of the Lie algebra $sl(2)$.
This analog can be introduced in frameworks of the general scheme
discussed in section \ref{sec:5}.  However, since the $sl(2)$-module $sl(2)^{\ot 2}$
is multiplicity free, this braided analog can be defined in a more
simple way. Let $\SL$ be the space $sl(2)$ endowed with the
action of the QG $\uq$ deforming the classical adjoint one.
There exists a unique (up to a nontrivial factor) $\uq$-morphism
$$
[\,\,,\,]: \SL\ot \SL\to \SL.
$$
Explicitly, it is given by the following multiplication table
\be
\begin{array}{c}
[b,b]=0,\quad [b,h]=-w \,b,\quad [b,c]=w \,{q \over 2_q}\,h,\quad
[h,b]=w\,q^2\, b,\\
\rule{0pt}{5mm} [h,h]=w\,(q^2-1) h,\quad [h,c]=-w\, c, \quad
[c,b]=-w\,{q\over 2_q}\, h,\quad [c,h]=w\, q^2\, c,\quad [c,c]=0
\end{array}
\label{qLie}
\ee
where $w\in \K,\, w\not=0$ is an arbitrary factor.

The corresponding adjoint action
$$
\ad\,x (y)=[x,y] \quad \forall\,\,x,y\in \SL
$$
gives rise to three operators
\be
B_q=\ad\,b,\qquad  H_q=\ad\,h, \qquad C_q=\ad\,c.
\label{bhc-q}
\ee
In the in the basis $\{b, h,c\}$ they are represented
by  the following matrices
\be
B_q=w\left(\begin{array}{ccc}
0&-1&0\\
0&0&{q\over 2_q}\\
0&0&0
\end{array}\right)\quad
H_q=w\left(\begin{array}{ccc}
q^2&0&0\\
0&q^2-1&0\\
0&0&-1
\end{array}\right)
\quad C_q=w\left(\begin{array}{ccc}
0&0&0\\
-{q\over 2_q}&0&0\\
0&q^2&0
\end{array}\right)\quad
\label{q-ad}
\ee
Note that these operators satisfy the relations (\ref{sys3}) with
$\h={w(q^4-q^2+1) \over 2_q} $. By specializing $q=1$ and $w=2$
we get the adjoint representation of the Lie algebra $sl(2)$.

Now, we want to extend the action of operators $B_q$, $H_q$ and
$C_q$ to the higher components of the algebra $\krq$. Such an extension
can be constructed via the coproduct described in section \ref{sec:4}.
However, we use another method which is similar to that used in definition of
the partial derivatives. Presenting a degree $n$ homogeneous
element $f\in \krq$ in the canonical form we define the action of the
extended operators as follows
$$
B_q(f)=n\, (B_q)_1 (f),\quad H_q(f)=n\, (H_q)_1 (f),\quad
C_q(f)=n\, (C_q)_1 (f)
$$
where as usual the subscript means that these operators are
applied to the first factors. (We keep the same notation for the
extended operators.)

Note that on each homogeneous component of the algebra $\krq$
the operators $B_q,\,H_q,\,C_q$ realize a representation of the
algebra (\ref{sys3}) but (in contrast  with the classical case)
with different factors $\h$.  Let $\h(n)$ be the value of
the factor $\h$ on the homogeneous component of the degree
$n$. A proof of this fact and the computation of the values
$\h(n)$ can be found in \cite{DGR}.

Now, upon replacing the operators $D_b,\,D_h,\, D_c$ in
(\ref{dirr}) by the operators $B_q,\,H_q,\,C_q$ respectively we
get the $q$-Dirac operator on the $q$-hyperboloid (more precisely,
a free $\K_q[H^2]$-module).
 Namely, we have
\be
{\Dir}_{\K_q[H^2]}=B_q \,\sigma_b+H_q\,
\sigma_h+C_q\,\sigma_c=\left(\begin{array}{cc}
{q\,H_q\over 2_q} & B_q\\
C_q& - {H_q \over 2_q\,q} \end{array}\right).
\label{trinad}
\ee
In virtue of (\ref{CHid})  this operator on the degree $n$ homogeneous
component satisfies the relation
$$
{\Dir}_{\K_q[H^2]}^2=q\, \h(n){\Dir}_{\K_q[H^2]}+{1\over 2_q}
(\qq\, B_qC_q+{H_q^2\over 2_q}+C_q B_q)\Id.
$$

Now, we discuss a way of definition of a $q$-analog of the Maxwell
operators on the quantum algebras in question.
 First, consider the classical Maxwell operator
$$
\Mw(\omega)=\partial\,d\,\omega = \De(\omega)-d\,
\partial(\omega), \,\,\,{\rm where}\,\,\,\om \in \Om^1,\,\,\,
\partial=*^{-1}\,d \,*,
$$
$d$ is the de Rham operator, and $*$ is the Hodge one.
The Hodge operator is introduced via a metric which
is assumed to be (quasi)Euclidian. Identifying the differential
forms
$$
\om=dt\,\al +dx\, \beta +dy\, \gamma +dz\, \de  \in \Om(\R^4),\quad \al,\,
\beta,\, \gamma,\,\delta \in \K[\R^4]
$$
with the columns $(\al,\, \beta,\, \gamma,\,\delta)^{\sf T}$ we can
present the Maxwell operator on the Minkowski space as follows
$$
\Mw_{\K[\R^4]}\left(\begin{array}{c}
\al\\
\beta\\
\gamma\\
\delta
\end{array}\right)=
\left(\begin{array}{c}
\De_{\K[\R^4]}(\al)\\
\De_{\K[\R^4]}(\beta)\\
\De_{\K[\R^4]}(\gamma)\\
\De_{\K[\R^4]}(\delta)
\end{array}\right)-\left(\begin{array}{c}
\partial_t\\
\partial_x\\
\partial_y\\
\partial_z
\end{array}\right)(\partial_t\,-\partial_x,\,-\partial_y,\,-\partial_z)
\left(\begin{array}{c}
\al\\
\beta\\
\gamma\\
\delta
\end{array}\right).
$$
In a similar manner we can realize the Maxwell operator on the
Euclidian spaces $\R^3\cong so(3)^*$ and $\R^3\cong sl(2)^*$
(see \cite{DG} for detail).

Now, introduce $q$-analogs of the Maxwell operators on the
algebras $\krq$ and $\krrq$ respectively by the relations
$$
\Mw_{\K_q[\R^3]}\left(\begin{array}{c}
\al\\
\beta\\
\gamma
\end{array}\right)=
\left(\begin{array}{c}
\De_{\K_q[\R^3]}(\al)\\
\De_{\K_q[\R^3]}(\beta)\\
\De_{\K_q[\R^3]}(\gamma)
\end{array}\right)-
\left(\begin{array}{c}
\partial_b\\
\partial_h\\
\partial_c
\end{array}\right)
(q^{-1}\partial_c,\,2_q\partial_h,\,q\,
\partial_b)\left(\begin{array}{c}
\al\\
\beta\\
\gamma
\end{array}\right)
$$

\be
\Mw_{\K_q[\R^4]}\left(\begin{array}{c}
\al\\
\beta\\
\gamma\\
\delta
\end{array}\right)=
\left(\begin{array}{c}
\De_{\K_q[\R^4]}(\al)\\
\De_{\K_q[\R^4]}(\beta)\\
\De_{\K_q[\R^4]}(\gamma)\\
\De_{\K_q[\R^4]}(\delta)
\end{array}\right)-
\left(\begin{array}{c}
\partial_b\\
\partial_h\\
\partial_c\\
\partial_{\ell}
\end{array}\right)
(q^{-1}\partial_c,\,2_q\partial_h,\,q\, \partial_b,\, \epsilon
\partial_{\ell})\left(\begin{array}{c}
\al\\
\beta\\
\gamma\\
\delta
\end{array}\right)
\label{dir}
\ee
where $\al,\,\beta,\,\gamma\in \K_q[\R^3]$
(resp., $\al,\,\beta,\,\gamma,\delta\in \K_q[\R^4]$) and the
operator $\De_{\K_q[\R^3]}$ (resp., $\De_{\K_q[\R^4]}$) is defined
by formula (\ref{Lap1}) (resp., (\ref{Lap2})).

The basic property of these operators is that their kernels are
similar to those of the classical Maxwell operators. Namely, the
kernel of the operator $ \Mw_{\K_q[\R^3]}$ (resp.,
$\Mw_{\K_q[\R^4]}$) contains all columns of the form $(\partial_b
\varphi,\,\partial_h \varphi,\,\partial_c \varphi)^{{\sf T}}$
(resp., $(\partial_b \varphi,\,\partial_h \varphi,\,\partial_c
\varphi,\, \partial_{\ell} \varphi)^{{\sf T}}$). This property is
a consequence of the fact that the $q$-Laplace operators are central
in the corresponding algebras.

In order to get the $q$-Maxwell operator on the $q$-hyperboloid algebra we
present the partial derivatives in a form similar to that expressing the classical
operators in the spherical coordinates. We do not know any quantum analogs
of the angle variables, but instead we use the tangent braided vector fields
which are analogs of the infinitesimal hyperbolic rotations. First, consider the
space $\R^3\cong so(3)^*$. Let us present the partial derivatives in the Euclidean
variables $x,y,z$ as follows
\be
\partial_x={y\, Z-z\, Y \over r^2}+{x\over r}\partial_r =
{y\, Z-z\, Y \over \rho}+2\, x\, \partial_\rho\,, \quad (c.p.\rightarrow \partial_y,\;\partial_z)
\label{part}
\ee
where $\rho = r^2=(x^2+y^2+z^2)$ and {\it c.p.} stands for the cyclic permutations
$x\rightarrow y\rightarrow z\rightarrow x$.

Here, instead of the derivatives in the angles we use the vector
fields $X,\, Y,\, Z$ tangent to all spheres $S^2_r=\{(x,y,z)\,|\,
x^2+y^2+z^2=r^2\}$.
These tangent fields are bound by the relation
\be
x\, X+y\, Y+z\, Z=0.
\label{modpro}
\ee
Besides, they commute with the derivatives $\partial_r$ and $\partial_\rho$.
Note that the derivative $\partial_\rho$ acts on the Cartesian
variables $x,\, y,\,z$ as follows
\be
\partial_{\rho}\, x={x\over 2\rho}, \quad \partial_{\rho}\, y={y\over 2\rho},
\quad \partial_{\rho}\, z={z\over 2\rho}.
\label{derrho}
\ee

In terms of the vector fields $X,\, Y,\, Z$ and the derivative
$\parho$ the Laplace operator on $\R^3$ takes the following form
\be
\Delta_{\kr}=\partial^2_x+\partial^2_y+\partial^2_z=
{X^2+Y^2+Z^2\over \rho}+6\partial_\rho+4 \rho\partial_\rho^2.
\label{Lap}
\ee

Taking in consideration this formula it is reasonable to extend the algebra
$\krr$ by the element $\rho^{-1}$  (and to proceed in a similar way with the
other algebras considered below). For detail we refer the reader to \cite{GS3}.

In a similar way we can proceed while dealing with $\R^3\cong
sl(2)^*$ endowed with a $sl(2)$-invariant metric. Then, using
the hyperbolic infinitesimal rotations
$$
B=-2b\, \partial_h+h\,\partial_c,\quad H=2b\,
\partial_b-2c\,\partial_c,\quad C=-h\partial_b+2c\, \partial_h
$$
we can present the derivatives $\partial_b$,$\partial_h$ and $\partial_c$
in a form similar to (\ref{part}) (see formula (\ref{psds}) below for $q=1$).

Note that the operators $B,\, H,\,C$ are not independent but are
bound by a relation analogous to (\ref{modpro})
$$
b\,C+{h\,H\over 2}+ c\, B=0.
$$

A similar statement is valid for the braided vector fields $B_q$,
$H_q$ and $C_q$ acting on the algebra $\krq$.

\begin{proposition}
The operators $B_q$, $H_q$ and $C_q$ obey the relation
{\rm
\be
\qq b\,C_q+{h\,H_q\over 2_q}+q c\, B_q=0.
\label{re}
\ee}
\end{proposition}

{\bf Proof.} It suffices to check this relation on the generators
of the algebra $\krq$. Taking in consideration our method of
prolongation of the operators $B_q$, $H_q$ and $C_q$ to the
higher components we can conclude that it remains true on the whole
algebra $\krq$.\hfill\rule{6.5pt}{6.5pt}
\medskip

Note that the relation (\ref{re}) is independent on the normalization
of the operators $B_q$, $C_q$ and $H_q$ on the higher homogeneous
components of the algebra $\krq$.

All combinations $\al\,B_q+\beta\,H_q+\gamma\,C_q$ where
$\al,\, \beta,\, \gamma\in \krq$ are called {\em the braided tangent vector
fields}. The space of such fields is a left $\krq$-module
$$
M=\krq^{\oplus 3}/\overline{M}\quad {\rm where}\quad
\overline{M}=\{\varphi (\qq b\,C_q+{h\,H_q\over 2_q}+q c\, B_q)\,|\,
\forall\, \varphi\in \krq\}.
$$
On any $q$-hyperboloid algebra $\khq$ the modules $\overline{M}$
and $M$ are projective. The module $M$ is similar to $\bw_q^1(H^2)$
though one of them is right module and the other is left one.

Also, we need the following $q$-analog of the variable $\rho$.
We put
$$
\rho_q={1\over 2_q}\left(\qq\,b\,c+{h^2\over 2_q}+q\, c\,b\right).
$$
In analogy with the classical formulae, we introduce the derivative
$\partial_{\rhoq}$ setting by definition
$$
\partial_{\rhoq} b={b\over 2\, \rhoq},\quad \partial_{\rhoq} h={h\over
2\, \rhoq},
\quad \partial_{\rhoq} c={c\over 2\, \rhoq}.
$$
In addition, we assume that the derivative $\partial_{\rhoq}$ is
subject to the usual Leibnitz rule. It is easy to see that this
way of introducing the derivative $\partial_{\rhoq}$ is compatible
with the defining relations of the algebra $\krq$. Note that we
impose the Leibnitz rule only on the derivative in a central
element. However, here we can also do without this rule by
applying the method above. Namely, assuming $f\in \krq$ to be
a homogeneous element of degree $n$ we apply this derivative to
its first factor and multiply the result by $n$.

\begin{proposition}
The following operator equalities are valid on
the algebra $\krq$
{\rm
\be
D_b=\frac{q^{-2}}{2_q\rhoq}\, \BBB+ \frac{2b}{2_q}
\,\partial_{\rhoq}, \quad
D_h=\frac{q^{-2}}{2_q\rhoq}\,\HHH+\frac{2h}{2_q}\,\partial_{\rhoq},
\quad D_c=\frac{q^{-2}}{2_q\rhoq}\,\CCC+
\frac{2c}{2_q}\,\partial_{\rhoq},
\label{psds}
\ee}
where $\BBB=q^2 hB_q-b H_q$, $\HHH=q2_q(bC_q-cB_q)+(q^2-1)hH_q$,
$\CCC= q^2cH_q- hC_q$ and in the definition {\rm (\ref{bhc-q}) -- (\ref{qLie})}
of the operators $B_q,\, H_q,\, C_q$ we put $w=1$.
\end{proposition}

Analogously to the previous proposition it suffices to check these
relation on the generators of the algebra $\krq$.
We call the relations (\ref{psds}) the {\em pseudospherical form} of
the derivatives $D_b$, $D_h$ and $D_c$.

Now, define the $q$-Laplace operator on the $q$-hyperboloid with
$\rho_q=q^{-2} 2_q$ by
$$
\De_{\K_q[H^2]}=\qq\,\BBB\,\CCC+{1\over 2_q} \HHH^2+
q\,\CCC\,\BBB
$$
and the $q$-Maxwell operator by
$$
\Mw_{\K_q[H^2]}\left(\begin{array}{c}
\al\\
\beta\\
\gamma
\end{array}\right)=e'\,\left(
\left(\begin{array}{c}
\De_{\K_q[H^2]}(\al)\\
\De_{\K_q[H^2]}(\beta)\\
\De_{\K_q[H^2]}(\gamma)
\end{array}\right)-
\left(\begin{array}{c}
\qq\CCC\\
{\HHH\over 2_q}\\
q\BBB
\end{array}\right)
(\BBB,\,\HHH,\,\CCC)\left(\begin{array}{c}
\al\\
\beta\\
\gamma
\end{array}\right)\right),
$$
where $(\al,\,\beta,\,\gamma)^{\sf{T}} \in e'\,\K_q[H]^{\oplus 3}$,
$e'=1-{\overline e}'$, and ${\overline e}'$ is defined by
(\ref{eprim}).

This method of defining  the $q$-Laplace and the $q$-Maxwell operators on
the $q$-hyperboloid algebra was suggested in \cite{DG}. It is motivated by the
classical case, where on a subvariety of an affine space these operators are
restrictions of similar operators on the ambient space (see \cite{DG} for detail).

\section*{Acknowledgement} The work of one of the authors (P.S.) was partially supported by
the RFBR grant 08-01-00392-a and the joint RFBR and DFG grant 08-01-91953. The work of D.G. and
P.S. was partially supported by the joint RFBR and CNRS grant 09-01-93107.

\end{document}